\documentclass[11pt]{article}

\usepackage[margin=1in]{geometry}  
\usepackage{graphicx}           
\usepackage{amsmath, accents}              
\usepackage{dsfont} 
\usepackage{mathrsfs}
\usepackage{amsthm}    
\usepackage{amssymb}           
\usepackage{soul}
\usepackage{xcolor}
\usepackage{longtable}
\usepackage{float}
\usepackage{algorithm}
\usepackage{algpseudocode}
\usepackage{geometry}
\usepackage{amsfonts,mathtools}
\usepackage{amsxtra}
\usepackage{amstext}
\usepackage{amsbsy}
\usepackage{amscd}
\usepackage{float}
\usepackage{bm}
\usepackage[bottom]{footmisc}
\usepackage{multicol}
\usepackage{multirow}

\newtheorem{thm}{Theorem}[section]

\newtheorem{defi}[thm]{Definition}

\newtheorem{remark}[thm]{Remark}
\newtheorem{prop}[thm]{Proposition}

\theoremstyle{remark}
\newtheorem{algo}[thm]{Algorithm}

\newcommand{\Real}{\mathbb{R}}

\begin{document}

\nocite{*}

\title{A BDF2-Semismooth Newton Algorithm for the Numerical Solution of the Bingham Flow with Temperature Dependent Parameters\thanks{Supported in part by the Escuela Polit\'ecnica Nacional del Ecuador, under Project PIJ 16-05 ``Optimización Matemática del Flujo Dependiente de la Temperatura de Fluidos Viscoplásticos y Elasto-Viscoplásticos''.}}

\author{Sergio~Gonz\'alez-Andrade \\ \footnotesize Research Center on Mathematical Modeling (MODEMAT) and \\\footnotesize Departamento de Matem\'atica -  Escuela Polit\'ecnica Nacional\\\footnotesize Ladr\'on de
Guevara E11-253, Quito 170413, Ecuador\\\footnotesize \tt sergio.gonzalez@epn.edu.ec} 
\date{\today}
\maketitle

\begin{abstract}
This paper is devoted to the numerical solution of the non-isothermal instationary Bingham flow with temperature dependent parameters by semismooth Newton methods. We discuss the main theoretical aspects regarding this problem. Mainly, we discuss the existence of solutions for the problema, and focus on a multiplier formulation which leads us to a coupled system of PDEs involving a Navier-Stokes type equation and a parabolic energy PDE. Further, we propose a Huber regularization for this coupled system of partial differential equations, and we briefly discuss the well posedness of these regularized problems. A detailed finite element discretization, based on the so called (cross-grid $\mathbb{P}_1$) - $\mathbb{Q}_0$ elements, is proposed for the space variable, involving weighted stiffness and mass matrices. After discretization in space, a second order BDF method is used as a time advancing technique, leading, in each time iteration, to a nonsmooth system of equations, which is suitable to be solved by a semismooth Newton algorithm. Therefore, we propose and discuss the main properties of a SSN algorithm, including the convergence properties. The paper finishes with two computational experiments that exhibit the main properties of the numerical approach.
  \vskip .2in

\vspace{0.2cm}
\end{abstract}

\section{Introduction}
The numerical simulation of complex fluids has received an increasing amount of attention because of the wide variety of fields where these materials play central roles. In particular, viscoplastic materials are widely used in food industry (production of sauces and pastes), geophysics and other important fields of application. Recently, the focus on the heat transfer in these materials and its effects on the different material regions are of great interest. The development of theoretical and numerical tools to understand and simulate the effect of temperature fields on the evolution of rigid zones of the material has a particular impact in the understanding of several processes like transportation and exploitation of waxy-crude oils, simulation of lava flows, optimization of processes in industrial bakery, production and packing of jellies and honey-based products, etc. From the modelling perspective, this study is carried out by coupling the viscoplastic Bingham flow model with energy equations for the temperature fields. 

In this interesting context, the analysis and simulation of viscoplastic flows with temperature dependent parameters represent a challenging and interesting field of research. From the analytical point of view, the viscoplastic model needs to be coupled with energy equations which are nonlinear and complex to analyze. Further, the temperature field should be inserted in the differential operators of the flow model, as functional weights. These models allow us to directly affect the mechanical properties of the material by changing the viscosity and the yield stress due to the action of temperature. This constitutes a key first step in the design of a control strategy for viscoplastic materials which are subject to heating by means of controlled heat sources.

There are several contribution on the theoretical analysis of the problem, starting with the classical work \cite{DuvChaleur}, where the problem was introduced and studied from the variational point of view. Recently, the steady flow of the Bingham model with temperature dependent parameters has been deeply analyzed in \cite{ConRod}. Besides a detailed discussion regarding the existence and uniqueness of solutions for the problem, the authors analyze the asymptotic behavior of the velocity field when the thermal conductivity of the material tends to infinity. This regime can be understood as a superfluid model, where the parameters stop depending on the geometry (non-local parameters). Following a similar perspective, in \cite{Messelmi} the authors extend this analysis to the steady flow of the Herschell-Bulkley model for a specific shear-thinning setting.

In contrast with the analytical work developed, the numerical solution of these models is, to the best of our knowledge, scarce. In \cite{Vinay}, the authors propose a numerical solution of this problem with the background of waxy crude oil flows. The authors consider that this kind of crude oils behaves as a viscoplastic material (Bingham model), and they numerically simulate the flow in a geometry representing a pipe. They propose a discretization based on a finite volume method. However, they dismiss the effect of the dissipation term in the energy equation, following mechanical properties of the flow under analysis. Finally, they proposed an Augmented Lagrangian method with an Uzawa type algorithm to solve the discretized system. 

In \cite{GloWachs}, the authors introduce the so called Houska model, which represents a generalized approach to this kind of coupled problems. In this model, the energy equations is replaced by a transport equation and the flow is coupled with a structural parameter function which can be seen as temperature, concentration of certain components in mixtures-type materials, etc.

In this article, we study the Bingham flow with temperature dependent parameters, by analyzing the instationary Bingham model coupled with a parabolic partial differential equation for the temperature field. We first briefly discuss the variational formulation of this coupled fluid-heat system. Next, we discuss the existence of solutions for the coupled model and propose a multiplier formulation, which yields a system of equations formed by a Navier-Stokes type system (the Bingham model) and a parabolic PDE, coupled through the convection terms in both equations and a dissipation term in the energy equation. The analytical part of the article closes with the introduction of a smoothing procedure for the coupled multiplier system, based on the Huber regularization approach. We propose such approximation to the problem since it has proven to be an efficient way to obtain fast convergent algorithms, maintaining an accurate approximation of the mechanical behavior of the flow (see \cite{DlRG1}).

For the numerical solution of the problem, we aim at extending the approach developed in [6] to the case of the coupled heat-fluid flow under analysis. We start by discussing in detail a finite element discretization for the space variable. This discussion includes the construction of weighted stiffness and mass matrices for both the flow and the energy partial differential equations. Further, we discuss the time advancing scheme, considering a second order BDF time method, combined with a lag operator, as proposed in \cite{baker, DlRG1}. In such a way, we obtain a fully discretized coupled fluid-heat system to be solved. The main characteristic of this system is that, because of the discretization approach, the system matrices are convective free. This means that the convective matrices can be constructed only with information of the velocity and temperature fields in previous time steps. Further, the coupled system involves, at each time step, a semismooth system of equations related to the regularized Bingham model. The solution of this nonsmooth system is carried out with a semismooth Newton method, wich we propose and analyse. The main result is that the superlinear convergence of the algorithm holds in this case.

The paper is organized as follows. In Section \ref{sec:problem} we introduce the problem given by a coupled system of nonlinear PDEs. Next, we briefly discuss the variational formulation of the PDEs system, which yields a coupled system of a variational inequality of the second kind for the velocity field and a parabolic nonlinear energy PDE for the temperature field. Finally, we characterize the solutions for this coupled system by the introduction of a tensor multiplier and propose a local smoothing procedure on the multiplier system based on the Huber approach. In Section \ref{sec:numerics}, we discuss the fully discretization of the regularized system. First, we propose a stable finite element discretization for the space variable based on the so called (cross-grid $\mathbb{P}_1$)-$\mathbb{Q}_0$ elements. Next, we discuss the time advancing technique based on the application of a semi-implicit BDF2 method. Section \ref{sec:algo} is devoted to the construction of the BDF2-SSN algorithm to numerically solve the coupled problem. Also, we construct the semismooth Newton inner algorithm to solve the flow equations. Finally, in Section \ref{sec:comput}, two detailed numerical experiments are carried out in order to verify the theoretical properties of the proposed approach.

\section{Problem Statement}\label{sec:problem}
Let $\Omega\subset\Real^d$, $d=2,3$, be an open and bounded set with Lipschitz boundary $\partial \Omega$. Let us assume that there exist $\Gamma,\Gamma_0\subset\partial\Omega$, such that $\partial \Omega = \Gamma\cup\Gamma_0,\,\,\Gamma\cap\Gamma_0=\emptyset$ and $|\Gamma|>0$. Given a final time $T_f>0$, we define the following space-time sets $Q:=\Omega\times(0,T_f)$ and $\Sigma:=\partial\Omega\times(0,T_f)$.

We are concerned with the non-isothermal flow of a Bingham fluid, considering that both the viscosity and the plasticity threshold depend on temperature. We assume the existence of volume forces acting on the fluid. The constitutive equations for this phenomenon are given by the following problem: find a velocity field $\mathbf{u}:Q\rightarrow \Real^d$, a pressure field $p:Q\rightarrow\Real$ and a temperature field $\theta:Q\rightarrow \Real$ such that

\begin{equation}\label{floweq}\tag{$\mathcal{B}$}
\begin{array}{cll}
\rho\left[\partial_t \mathbf{u}+ (\mathbf{u}\cdot\nabla)\mathbf{u}\right]- \boldsymbol{\nabla}\cdot\boldsymbol{\tau}(\theta) + \nabla p= \mathbf{f},&\mbox{in $\Omega$}\vspace{0.2cm}\\ \nabla\cdot \mathbf{u} =0,&\mbox{in $\Omega$}\vspace{0.2cm}\\ 
\boldsymbol{\tau}(\theta)=\mu(\theta)\mathcal{E}\mathbf{u} + g(\theta)\frac{\mathcal{E}\mathbf{u}}{\|\mathcal{E}\mathbf{u}\|}&\mbox{if $\mathcal{E}\mathbf{u}\neq 0$,}\vspace{0.2cm}\\ \|\boldsymbol{\tau}(\theta)\|\leq g(\theta)&\mbox{if $\mathcal{E}\mathbf{u}= 0$,}\vspace{0.2cm}\\\mathbf{u}=0,&\mbox{on $\Sigma$}, \vspace{0.2cm}\\ \mathbf{u}(x,0)=\mathbf{u}_0,&\mbox{in $\Omega$}.
\end{array}
\end{equation}

\begin{equation}\label{eneq}\tag{$\mathcal{E}$}
\begin{array}{cll}
\rho C_p\left[\partial_t\theta+\mathbf{u}\cdot \nabla\theta\right] -\kappa\Delta\theta =\boldsymbol{\tau}(\theta):\mathcal{E}\mathbf{u}- \alpha\theta,&\mbox{in $\Omega$.}\vspace{0.2cm}\\\frac{\partial\theta}{\partial\mathbf{n}}=0,&\mbox{on $\Gamma_0\times(0,T)$}\vspace{0.2cm}\\ \kappa\frac{\partial\theta}{\partial\mathbf{n}}+\beta\, C_p\theta =0,&\mbox{on $\Gamma\times(0,T)$}	\vspace{0.2cm}\\\theta(x,0)=\theta_0,&\mbox{in $\Omega$}.
\end{array}
\end{equation}

System \eqref{floweq} represents the Bingham flow model. However, in this case, both the viscosity and the yield stress depend on  temperature, which is provided as the solution of the energy equation \eqref{eneq}. The analysis and numerical solution of the coupled system \eqref{floweq}-\eqref{eneq} is the main goal of this work. Here, $\rho$ denotes the density of the fluid, $C_p$ is the heat capacity, $\kappa$ stands for the thermal conductivity and $\mathbf{f}$ for the volume forces. $\alpha\geq0$ and $\beta\geq 0$ are given parameters and, as usual, $\mathbf{n}$ stands for the unit outward normal to $\partial\Omega$. Notice that, in addition to the action provided by the dissipation energy factor $\boldsymbol{\tau}(\theta):\mathcal{E}\mathbf{u}$, we admit the effect of a possible external heat sink proportional to the temperature, if $\alpha>0$. We assume only nonslip boundary conditions for the flow velocity. For the energy equation, we impose a Robin boundary condition in $\Gamma$, with given coefficient $\beta$, and Neumann homogeneous boundary conditions in $\Gamma_0$. Further, we assume given initial conditions $\mathbf{u}_0$ and $\theta_0$ for the velocity and temperature, respectively. Finally, by following \cite{DuvChaleur}, we assume that $\rho=1$.

In the following section we analyze the variational formulation of the coupled system \eqref{floweq}-\eqref{eneq}. Further, we discuss existence of solutions and propose a mixed formulation, based on a regularization approach, for this system.

\subsection[Variational Formulation]{Variational Formulation}
We consider the Bingham incompressible flow model. Therefore, we introduce the  divergence free spaces $V:=\{\mathbf{v}\in \mathbf{H}_0^1(\Omega)\,:\, \nabla\cdot \mathbf{v}=0\}$ and $H:=\{\mathbf{w}\in \mathbf{L}^2(\Omega)\,:\, \nabla\cdot \mathbf{w}=0\,\mbox{ and }\, \frac{\partial\mathbf{w}}{\partial\mathbf{n}}|_{\partial\Omega}=0\}$. Next, the variational formulation of the flow model is given by the following problem: find $\mathbf{u}(t)\in V$, a.e. in $(0,T)$ such that

\begin{equation}\label{flowvar}\tag{$\mathcal{VB}$}
\begin{array}{cll}
\left(\partial_t\mathbf{u}(t),\mathbf{v}-\mathbf{u}\right) + \int_\Omega\langle (\mathbf{u}(t)\cdot\nabla)\mathbf{u}(t)\,,\,\mathbf{v}-\mathbf{u}(t)\rangle\,dx+ \int_{\Omega}\mu(\theta)\,\mathcal{E}\mathbf{u}(t)\,:\,\mathcal{E}(\mathbf{v}-\mathbf{u}(t))\,dx \vspace{0.2cm}\\ \hspace{3cm}+ \int_\Omega g(\theta)\,\|\mathcal{E}\mathbf{v}\|\,dx-  \int_\Omega g(\theta)\,\|\mathcal{E}\mathbf{u}(t)\|\,dx\geq \int_\Omega(\mathbf{f}\,,\,\mathbf{v}-\mathbf{u}(t))\,dx,\quad \forall\,\mathbf{v}\in V\vspace{0.2cm}\\ \mathbf{u}(0)=\mathbf{u}_0
\end{array}
\end{equation}

One important issue regarding this formulation is the definition of the functions $\mu$ and $g$, which measure the dependence of the viscosity and the plasticity threshold with temperature, respectively. These functions must be defined in such a way that the integrals in \eqref{flowvar} are well posed. Usually, they are required to satisfy the following hypotheses \cite{ConRod,DuvChaleur,Messelmi}: $\mu\in C(\Real)$, $g\in C(\Real)$ and there exist $\mu_0,\mu_1>0$, and $g_0\geq0$ and $g_1>0$ such that
\begin{equation}\label{muandg}
\begin{array}{lll}
\mu_0\leq\mu(\theta)\leq \mu_1<+\infty,&\forall \theta\in \Real\vspace{0.2cm}\\
g_0\leq g(\theta)\leq g_1<+\infty,&\forall \theta\in \Real.
\end{array}
\end{equation}

In our approach, we propose to follow a similar setting as the one given by the Houska model (see \cite{GloWachs}). This model is a generalization of the coupled system studied here, where the  flow parameters depend on a structure parameter, which can be seen as several magnitudes, including temperature. The main difference is that the Houska model does not include a diffusion term for the energy equation, which is replaced by a transport equation. In that model, the viscosity and the plasticity threshold are supposed to be affine functions of the  structure parameter.  Summarizing, by following ideas in the Houska model, we propose to analyse the system \eqref{floweq}-\eqref{eneq} considering the following affine functions of the temperature parameter $\theta$
\begin{equation}\label{gmu}
g(\theta):=g_0 + \delta_g\,\theta\,\,\,\mbox{ and }\,\,\, \mu(\theta):= \mu_0 + \delta_\mu\,\theta.
\end{equation}
We assume that $\theta$ is a function of $x$ and $t$ which takes its values in the closed interval $[0,1]$ (see \cite[Section 3]{GloWachs}). Moreover, $\delta_g$ and $\delta_\mu$ stand for the expected variation, due to the action of heating, in the yield stress and the viscosity, respectively. Clearly, with this consideration these functions satisfy the condition \eqref{muandg}. Further, it is accurate for us to suppose that $g(\theta)\geq 0$, since the fluid is supposed to be a Bingham material which will be altered by the action of temperature. Finally, we can state that all the integrals in \eqref{flowvar} are well posed, so the variational formulation for the flow system holds.

Now, let us focus on the variational formulation of the energy equation \eqref{eneq}. Let $1<q<\frac{d}{d-1}$. Then, the variational formulation of the energy equation is given by: find $\theta(t)\in W^{1,q}(\Omega)$ a.e. in $(0,T_f)$ such that
\begin{equation}\label{varen}\tag{$\mathcal{VE}$}
\begin{array}{cll}
C_p\int_\Omega \partial_t \theta(t)\, \phi\,dx + C_p\int_\Omega (\mathbf{u}(t)\cdot \nabla\theta)\phi\,dx +\kappa \int_\Omega (\nabla \theta (t),\nabla\phi)\,dx + \beta\,C_p \int_\Gamma \theta(t)\,\phi\,dx \vspace{0.2cm}\\ \hspace{2cm}= \int_\Omega \left[\mu(\theta(t))\|\mathcal{E}\mathbf{u}(t)\|^2 + g(\theta(t))\|\mathcal{E}\mathbf{u}(t)\|\right]\phi\,dx -\alpha\int_\Omega \theta\,\phi\,dx,\,\,\forall\, \phi\in W^{1,q'}(\Omega),\vspace{0.2cm}\\\theta(0)=\theta_0,
\end{array}
\end{equation}
where $1/q+1/q'=1$. Note that the dissipation term is only nonzero in regions where  $\mathcal{E}\mathbf{u}(t)\neq 0$, a.e. in $(0,T_f)$. Therefore, we use the corresponding form of $\boldsymbol{\tau}(\theta)$ in the energy equation. Furthermore, the associated integral term in the variational formulation is well posed since, for $1<q<d/(d-1)$, we have that $q'>d$, which implies that $W^{1,q'}(\Omega)$ is continuously embedded in $L^\infty(\Omega)$.

\begin{remark}\label{rem:exist}
Existence and uniqueness of solutions for the coupled system \eqref{flowvar}-\eqref{varen}, to the best of our knowledge, constitutes an open problem. The instationary case has been deeply analyzed in \cite{ConRod}, where existence of solutions has been proved and the regularity $(\mathbf{u},\theta)\in V\times W^{1,q}(\Omega)$, $1<q<\frac{d}{d-1}$, has been established. However, for the instationary case there are only partial results. For instance, if we neglect the advection term in \eqref{varen} and consider a heat source/sink term which is not dependent on $\theta$, we can state, by following \cite[Th. 3.1]{DuvChaleur}, that there exist $\mathbf{u}\in L^2(V)$ such that $\partial_t\mathbf{u}\in L^2(V')$, and $\theta\in L^{q}(Q)$ solutions of \eqref{flowvar}-\eqref{varen}, for $1<q<\frac{d}{d-1}$, and considering that $\mathbf{f}\in L^2(V')$, $\mathbf{u}_0\in H$ and $\theta_0\in L^1(\Omega)$. The development of a general existence result needs further research in the PDEs theory. These results are out of the scope of this paper, thus, in the following, we assume the existence of solutions for \eqref{flowvar}-\eqref{varen}.
\end{remark}

Let us recall that we use the notation $L^p(W)$ for the spaces $L^p(0,T;W):=\{f:[0,T]\rightarrow W\,:\,\int_0^T\|f(t)\|^p\,dx<\infty\}$.

\subsection[Multiplier Approach]{Multiplier Approach}
The use of tensor multiplier-type functions provide a versatile characterization for the solutions of problems involving variational inequalities. This approach has been used in previous contributions focused on Bingham flow (see \cite{DlRG1,DlRG2}). With such a characterization, a partial differential equation involving a multiplier, together with additional complementarity relations, is obtained. Further, in the case of non-isothermal flow, \cite[Th. 3.2]{DuvChaleur} has established the existence of such a multiplier $\mathbf{q}$ for \eqref{flowvar}. The resulting system, in variational form, reads as follows: find $\mathbf{u}(t)\in \mathbf{H}^1_0(\Omega)$, $p(t)\in L_0^2(\Omega)$ and $\mathbf{q}(t)\in \mathbf{L}^{2\times 2}(\Omega)$ a.e. in $(0,T_f)$ such that.

\begin{equation}\label{multflow}\tag{$\mathcal{MB}$}
\begin{array}{cll}
\int_\Omega (\partial_t\mathbf{u}(t),\mathbf{v})\,dx +\int_\Omega\langle (\mathbf{u}(t)\cdot\nabla)\mathbf{u}(t)\,,\,\mathbf{v}\rangle\,dx+\int_\Omega \mu(\theta)\, (\mathcal{E}\mathbf{u}:\mathcal{E}\mathbf{v})\,dx \vspace{0.2cm}\\\hspace{3cm} + \int_\Omega g(\theta)\, (\mathbf{q}:\mathcal{E}\mathbf{v})\,dx -\int_\Omega p(t)\nabla\cdot\mathbf{v}\,dx = \int_\Omega (\mathbf{f}(t),\mathbf{v})\,dx,\,\,\forall \mathbf{v}\in \mathbf{H}^1_0(\Omega)\vspace{0.2cm}\\ \int_\Omega r\,\nabla \cdot\mathbf{u}(t)=0\,dx,\,\,\, \forall r\in L^2_0(\Omega)\vspace{0.2cm}\\ \|\mathbf{q}(x,t)\|\leq g(\theta),\,\,\mbox{a.e. in $Q$}\vspace{0.2cm}\\ (\mathbf{q}(x,t):\mathcal{E}\mathbf{u}(x,t)) = g(\theta)\|\mathcal{E}\mathbf{u}(x,t)\|,\,\,\mbox{a.e. in $Q$}\vspace{0.2cm}\\\mathbf{u}(x,0)=\mathbf{u}_0.
\end{array}
\end{equation}

In this system, the pressure function $p$ is recovered by a direct application of the de Rahm's Theorem, so the variational formulation is constructed by using the Sobolev space $\mathbf{H}^1_0(\Omega)$, as the test space (see \cite[Rem. 1.9, pp. 14.]{Temam}). 

The active and inactive sets for the flow are defined, respectively, by
\[
\mathcal{A}:=\{(x,t)\in Q\,:\, \|\mathcal{E}\mathbf{u}(x,t)\|\neq 0\}\,\,\mbox{ and }\,\,\mathcal{I}:=Q\setminus\mathcal{A}.
\]

Further, since the multiplier $\mathbf{q}$ is undetermined in the solid regions and not necessarily unique (see \cite{DuvChaleur}), the computational approach depends on projection techniques (see \cite{GloWachs, Vinay}). Moreover, instability in the numerical methods may occur due to the nonuniqueness of the multiplier. In this contribution, we propose to extend the approach based on local regularization techniques (see \cite{DlRG2,DlRG1}) which leads us to superlinear convergent methods of semismooth Newton type. 

\subsection[Huber Regularization Approach]{Huber Regularization Approach}
In this section, we introduce a family of regularized problems to approximate the solution of \eqref{multflow}-\eqref{varen}. This regularization approach is based on the so called Huber local smoothing of the stress tensor and it has proved to be efficient in the numerical solution of the Bingham flow in the stationary as well as in the instationary cases (see \cite{DlRG0,DlRG1}). Further, we have used this approach in the numerical solution of the convective flow of Bingham fluids in the Bousinessq paradigm \cite{DlRG2}. The Huber regularization provides an equivalent formulation as the one obtained with the bi-viscosity approach (see \cite{DlRG0,DlRG2}). Here, the active and inactive sets are approximated by sets depending on a regularization parameter. As the parameter grows, the regularized inactive set tends to the actual inactive set, which makes that the viscosity of the flow, in the solid regions, tends to infinity, which is the expected mechanical behavior. 

For a parameter $\gamma>0$, the regularized problem consists in: find $\mathbf{u}_\gamma(t)\in \mathbf{H}^1_0(\Omega)$, $p_\gamma(t)\in L_0^2(\Omega)$, $\mathbf{q}_\gamma(t)\in \mathbf{L}^{2\times 2}(\Omega)$ and $\theta_\gamma(t)\in W^{1,q}(\Omega)$ a.e. in $(0,T_f)$ such that.

\begin{equation}\label{regflow}\tag{$\mathcal{RB}$}
\begin{array}{cll}
\int_\Omega (\partial_t\mathbf{u}_\gamma(t),\mathbf{v})\,dx + \int_\Omega\langle (\mathbf{u}_\gamma(t)\cdot\nabla)\mathbf{u}_\gamma(t)\,,\,\mathbf{v}\rangle\,dx+\int_\Omega \mu(\theta_\gamma)\, (\mathcal{E}\mathbf{u}_\gamma:\mathcal{E}\mathbf{v})\,dx \vspace{0.2cm}\\\hspace{3cm} + \int_\Omega g(\theta_\gamma)\, (\mathbf{q}_\gamma:\mathcal{E}\mathbf{v})\,dx -\int_\Omega p_\gamma(t)\nabla\cdot\mathbf{v}\,dx = \int_\Omega (\mathbf{f}(t),\mathbf{v})\,dx,\,\,\forall \mathbf{v}\in \mathbf{H}^1_0(\Omega)\vspace{0.2cm}\\ \int_\Omega r\,\nabla \cdot\mathbf{u}_\gamma(t)=0\,dx,\,\,\, \forall r\in L^2_0(\Omega)\vspace{0.2cm}\\ \mathbf{q}_\gamma:=\left\{\begin{array}{lll}
g(\theta_\gamma(x,t))\frac{\mathcal{E}\mathbf{u}_\gamma(x,t)}{\|\mathcal{E}\mathbf{u}_\gamma(x,t)\|},\,\mbox{a.e. in $\mathcal{A}_\gamma$}\vspace{0.2cm}\\ \gamma\mathcal{E}\mathbf{u}_\gamma(x,t),\,\mbox{a.e. in $\mathcal{I}_\gamma$}
\end{array}\right.\vspace{0.2cm}\\\mathbf{u}_\gamma(x,0)=\mathbf{u}_0\vspace{0.2cm}\\
C_p\int_\Omega \partial_t \theta_\gamma(t)\, \phi\,dx + C_p\int_\Omega (\mathbf{u}_\gamma(t)\cdot \nabla\theta_\gamma)\phi\,dx+\kappa \int_\Omega (\nabla \theta_\gamma (t),\nabla \phi)\,dx + \beta\,C_p \int_\Gamma \theta_\gamma(t)\,\phi\,dx \vspace{0.2cm}\\ \hspace{1.5cm}= \int_\Omega \left[\mu(\theta_\gamma(t))\|\mathcal{E}\mathbf{u}_\gamma(t)\|^2 + g(\theta_\gamma(t))\|\mathcal{E}\mathbf{u}_\gamma(t)\|\right]\phi \,dx -\alpha\int_\Omega \theta_\gamma\,\phi\,dx,\,\,\forall\, \phi\in W^{1,q'}(\Omega),\vspace{0.2cm}\\\theta_\gamma(0)=\theta_0,
\end{array}
\end{equation}
Here the regularized active and inactive sets are defined by
\[
\mathcal{A}_\gamma:=\{(x,t)\in Q\,:\, \gamma\|\mathcal{E}\mathbf{u}(x,t)\|\geq g(\theta_\gamma(x,t))\}\,\,\mbox{ and }\,\,\mathcal{I}_\gamma:=Q\setminus\mathcal{A}_\gamma.
\]
Note that, for given $\theta_\gamma(x,t)$ a.e. in $Q$, we have that $\|\mathbf{q}_\gamma(x,t)\|\leq g(\theta_\gamma(x,t))$ a.e. in $Q$.

Given $\theta_\gamma$, we usually rewrite the equation for $\mathbf{q}_\gamma$ as follows
\[
\max(g(\theta_\gamma(x,t)), \gamma\|\mathcal{E}\mathbf{u}(x,t)\|)\mathbf{q}_\gamma(x,t) = \gamma g(\theta_\gamma(x,t))\mathcal{E}\mathbf{u}_\gamma(x,t),\,\,\mbox{a.e. in $Q$}.
\]

\begin{remark}
Considering the discussion in Remark \ref{rem:exist}, we assume the existence of solutions for the energy equation \eqref{varen}. Taking this fact into account, the existence of solutions for \eqref{regflow} follows from \cite[Th. 3.1]{DlRG1}. Furthermore, \cite[Th. 3.3]{DlRG1} guarantee that for a given $\theta_\gamma$, $\mathbf{u}_\gamma\rightarrow \mathbf{u}$ strongly in $L^2(V)$ and $\mathbf{q}_\gamma\rightarrow\mathbf{q}$ weakly in $L^2(\mathbf{L}^{2\times 2})$. This convergence result shows that our approach is consistent and produce reliable approximations for the non-isothermal flow under study.
\end{remark}

\section[Space-Time Discretization]{Space-Time Discretization}\label{sec:numerics}
In this section, we propose a space-time discretization scheme for both the flow equation and the energy equation in the coupled system \eqref{regflow}. For the flow equation, we use the scheme proposed in \cite{DlRG1}, which is based on a combination of a first-order finite element approximation for the space variable with (cross-grid $\mathbb{P}_1$)-$\mathbb{Q}_0$ elements and the semi-implicit BDF2 scheme for the time variable. This class of finite elements allows us to use the same test functions for the velocity gradient and the dual variable. In such a way, a direct relation between these two variables is obtained, and an accurate determination of active and inactive sets is achieved. On the other hand, the semi-implicit BDF2 is an advancing scheme that leads us to convection-independent systems in each time step. This is particularly useful to reduce the computational cost in every step of the Newton iteration.

For the energy equation we propose a first-order finite element method for the space variable and a BDF2 method for the time variable. In such a way, we also obtain a convection-independent system. In this case, the main advantage is that the fully discretized equation becomes linear and, consequently, computationally cheap to be solved. The temperature $\theta$ is discretized in the same nodes as the velocity, and a simple restriction method (weighting) is used to obtain the value of the temperature at each triangle. 

\begin{figure}
\begin{center}
\includegraphics[width=50mm, height=45mm]{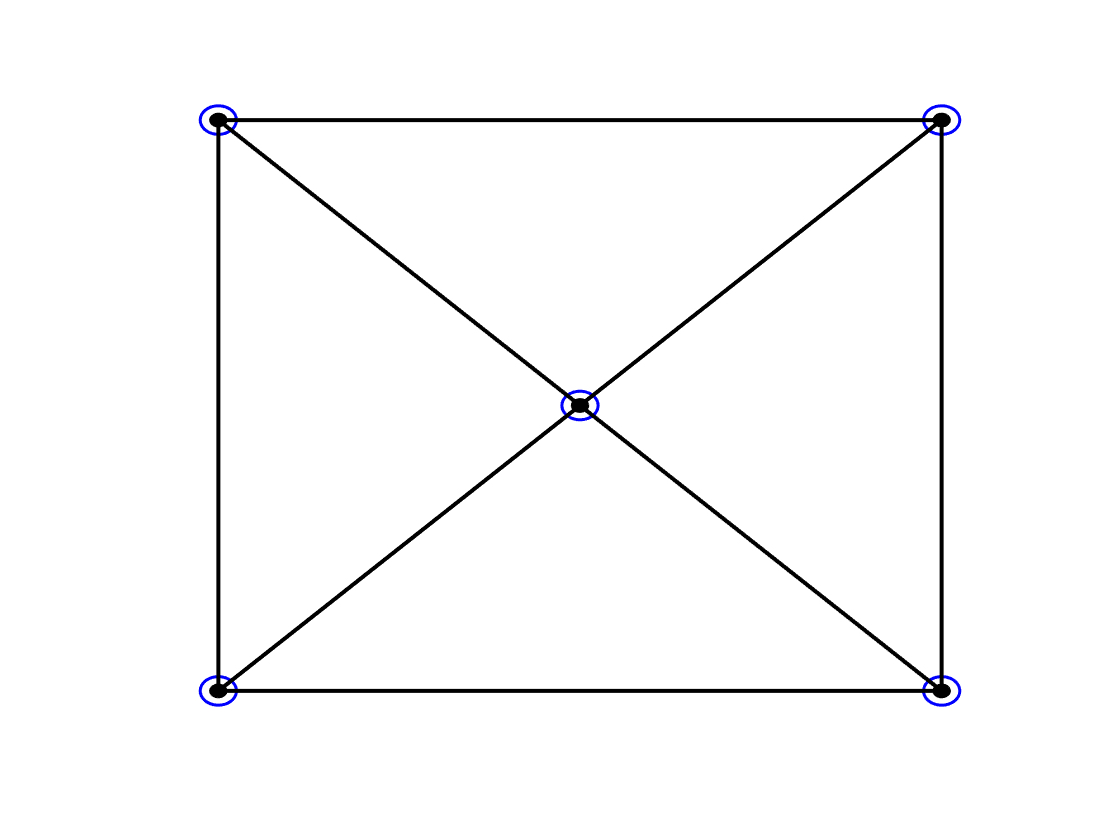}\includegraphics[width=50mm, height=45mm]{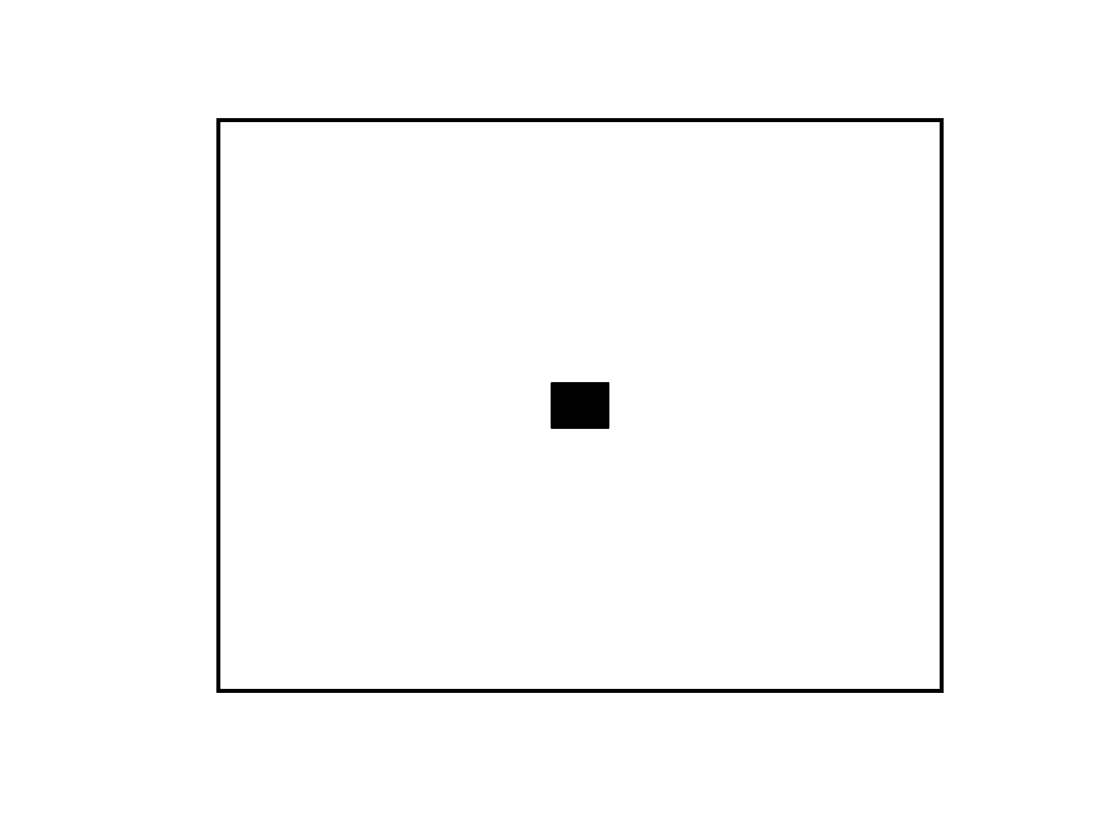}\includegraphics[width=50mm, height=45mm]{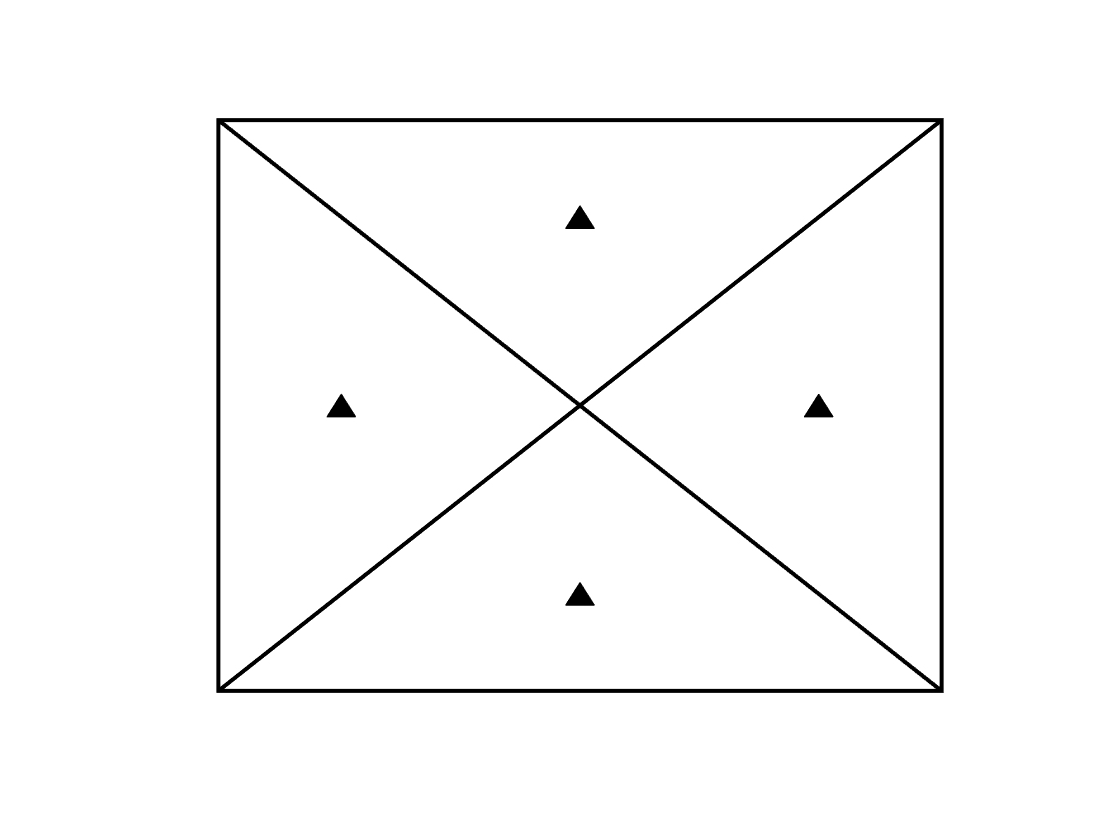}
\end{center}
\caption{\small Finite Elements for the coupled system. (Cross-grid $\mathbb{P}_1$)-$\mathbb{Q}_0$ macroelements: $\bullet$ are the nodes for the velocity, $\blacksquare$ are the nodes for the pressure and $\blacktriangle$ are the nodes for the multiplier. For the energy equation: $\bigcirc$ are the nodes for the temperature.}\label{fig:Fem}
\end{figure}
\subsection[Finite element discretization]{Finite element discretization}
By following \cite[Sec. 9.3]{QuarVall}, we start the discussion of the space discretization by introducing the following finite dimensional spaces 
\[
\begin{array}{lll}
\mathcal{Y}_\mathcal{T}^h&:=&\{\varphi^h\in L^2(\Omega)\,:\, \varphi^h|_T\in \Pi_1,\, \mbox{for all $T\in \mathcal{T}^h$}\}\vspace{0.2cm}\\
\mathcal{Y}_\mathcal{Q}^h&:=&\{\varphi^h\in L^2(\Omega)\,:\, \varphi^h|_T\in \Pi_1,\, \mbox{for all $Q\in \mathcal{Q}^h$}\}\vspace{0.2cm}\\
\mathcal{X}^h&:=&\mathcal{Y}_\mathcal{T}^h\cap C(\overline{\Omega}),
\end{array}
\]
where $\mathcal{Q}^h$ is a regular quadrangulation of $\Omega$, and $\mathcal{T}^h$ is the regular triangulation obtained by dividing any square in $\mathcal{Q}^h$ by using its two main diagonals \cite[Sec. 6]{QuarVall}.

Next, we define the discrete velocity and pressure spaces for the so called (cross-grid $\mathbb{P}_1$)-$\mathbb{Q}_0$ finite elements, as follows.
\begin{equation*}
\mathbf{V}^h:=(\mathcal{X}^h\cap H_0^{1}(\Omega))^2\,\,\,\mbox{and}\,\,\, U^h:=\mathcal{Y}_\mathcal{Q}^h\cap L^2_0(\Omega).
\end{equation*}
Further, we define the following discrete space for the multiplier
\[
\mathbf{W}^h:=\{(q_1^h,q_2^h,q_3^h,q_4^h)^\top\in (L^2(\Omega))^4\,:\, q_j^h |_T\in \Pi_0,\mbox{ for all }\, T\in \mathcal{T}^h\}.
\]
Clearly, we have that  $\mathbf{V}^h\subset\mathbf{H}^1(\Omega)$, $U^h\subset L^2(\Omega)$ and $\mathbf{W}^h\subset (L^2(\Omega))^4$. Further, we write that $\dim \mathbf{V}^h=2n$, $\dim\mathbf{W}^h=4m$ and $\dim U^h=\ell$ with $n,m,\ell\in \mathbb{N}$, respectively. To simplify the analysis, we assume that $\Omega$ has a polygonal boundary. The degrees of freedom associated to these elements are depicted in Figure \ref{fig:Fem}. 

It is remarkable that the (cross-grid $\mathbb{P}_1$)-$\mathbb{Q}_0$ elements satisfy the Ladyzhenskaya-Babuska-Brezzi (LBB) or inf-sup condition and, therefore, lead to a stable approximation of the Navier-Stokes-like systems, such as the Bingham model (see \cite[p. 435]{QuarVall}).

Further, for the energy equation, we discretize the temperature in the following finite dimensional spaces 
\[
X^h=\mathcal{X}^h\cap W^{1,q}(\Omega)\,\,\,\mbox{  and  }\,\,\,X^{h'}=\mathcal{X}^h\cap W^{1,q'}(\Omega).
\] 
The degrees of freedom are shown in Figure \ref{fig:Fem}, left.

By using the classical Galerkin approach, we obtain the following semi-discrete approximation for the coupled system \eqref{regflow}.

\begin{equation}\label{femflow}
\begin{array}{ccc}
\mathbf{M}^h\partial_t \vec{\mathbf{u}}(t) + \mathbf{C}^h(\vec{\mathbf{u}}(t))\vec{\mathbf{u}}(t) + \mathbf{A}^h_\mu (\vec{\theta}(t))\vec{\mathbf{u}}(t) + B^h \vec{p}(t)+ \mathbf{Q}^h_g (\vec{\theta}(t))\vec{\mathbf{q}}(t) = \vec{\mathbf{f}}(t)\vspace{0.2cm}\\ -(B^h)^\top \vec{\mathbf{u}}(t)=0,\vspace{0.2cm}\\  \max\left[G_T(\vec{\theta}(t)), \gamma N^h(\mathcal{E}^h\vec{\mathbf{u}}(t))\right]\star\vec{\mathbf{q}}(t) = \gamma \textrm{diag}(G_T(\vec{\theta}(t)))\mathcal{E}^h\vec{\mathbf{u}}(t),\,\,\mbox{a.e. in $[0,T_f]$ and $\forall \,T\in\mathcal{T}^h$}\vspace{0.2cm}\\ \vec{\mathbf{u}}(0)= \vec{\mathbf{u}}_0.
\end{array}
\end{equation}
\begin{equation}\label{femen}
\begin{array}{lcc}
C_p [M^h \partial_t \vec{\theta}(t) +C^h(\vec{\mathbf{u}}(t))\vec{\theta}(t)] +\kappa A^h \vec{\theta}(t)\vspace{0.2cm}\\\hspace{4cm} + \alpha M^h \vec{\theta}(t)+C_p \beta M_\Gamma^h \vec{\theta}(t) = [\mathbf{K}_\mu^h(\vec{\mathbf{u}}(t)) + \mathbf{K}_g^h(\vec{\mathbf{u}}(t)) ]\vec{\theta}(t) \vspace{0.2cm}\\\hspace{5cm} \vec{\theta}(0)=\vec{\theta}_0.
\end{array}
\end{equation}
Here, $G_T(\vec{\theta}(t))$ denotes a vector whose components are the values of $g(t)$ in the center of gravity of each $T\in\mathcal{T}^h$. $\vec{\mathbf{u}}(t)\in\Real^{2d}$, $\vec{p}(t)\in\Real^\ell$, $\vec{\mathbf{q}}(t) \in\Real^{4m}$, $\vec{\mathbf{\theta}}(t)\in\Real^d$, $\vec{\mathbf{u}}_0\in\Real^{2d}$ and $\vec{\theta}_0\in\Real^{d}$ are the time-dependent vectors of coefficients in the finite element representation of the 4-tuple $(\mathbf{u}^h(t),p^h(t),\mathbf{q}^h(t),\theta^h(t))\in \mathbf{V}^h\times U^h\times \mathbf{W}^h\times X^h$, and the initial conditions $\vec{\mathbf{u}}(0)$ and $\vec{\theta}(0)$, respectively. $\mathbf{M}^h$ and $M^h$ are the mass matrices for $\mathbf{V}^h$ and $X^h-X^{h'}$, respectively, while $A^h$ stands for the stiffness matrix associated to $X^h-X^{h'}$. $M^h_\Gamma$ is the boundary mass matrix constructed for the Robin boundary condition of the energy equation (see \cite[Sec. 4.6.2]{Larson}). Matrix $B^h$ is obtained in the usual way from the bilinear form $-(\cdot,\nabla\cdot (\cdot))_{L^2}$. 

The discretized convective matrices $\mathbf{C}^h(\vec{\mathbf{u}}(t))$ and $C^h(\vec{\mathbf{u}}(t))$ are given by
\[
\mathbf{C}^h(\vec{\mathbf{u}}(t))_{ij}:=\sum_{k=1}^{2d} \mathbf{u}_k \int_\Omega\langle (\boldsymbol{\varphi}_k\cdot\nabla)\boldsymbol{\varphi}_j\,,\,\boldsymbol{\varphi_i}\rangle\,dx\,\,\mbox{ and }\,\, C^h(\vec{\mathbf{u}}(t))_{ij}:=\sum_{k=1}^{2d} \mathbf{u}_k \int_\Omega (\boldsymbol{\varphi}_k\cdot\nabla \phi_j)\phi'_i\,dx,
\]
where $\boldsymbol{\varphi}_j$, $j=1,\ldots,2d$ are the basis functions of $\mathbf{V}^h$, $\phi_j$, $j=1,\ldots,d$ are the basis functions of $X^h$ and $\phi'_j$, $j=1,\ldots,d$ are the basis functions of $X^{h'}$, respectively. The discrete approximation of the deformation tensor $\mathcal{E}^h$ and the right hand side $\vec{\mathbf{f}}(t)$ are constructed by using the basis functions $\boldsymbol{\varphi}_j$, $j=1,\ldots,2d$ (see \cite[Sec. 5]{DlRG0}). Finally, the function $N^h:\Real^{4m}\rightarrow\Real^{4m}$ is defined by 
\[
N^h(q)_i=N^h(q)_{i+m}=\cdots=N^h(q)_{i+4m}:= |(q_i,q_{i+m},...,q_{i+4m})|,
\]
for $q\in\Real^{4m}$ and $i=1,\ldots ,m$. The values of $\mathcal{E}^h\vec{\mathbf{u}}$ and $N^h(\mathcal{E}^h\vec{\mathbf{u}})$ are given in the gravity centers of each $T\in \mathcal{T}^h$ (see Figure \ref{fig:Fem}).

\begin{figure}
\begin{center}
\includegraphics[width=75mm, height=60mm]{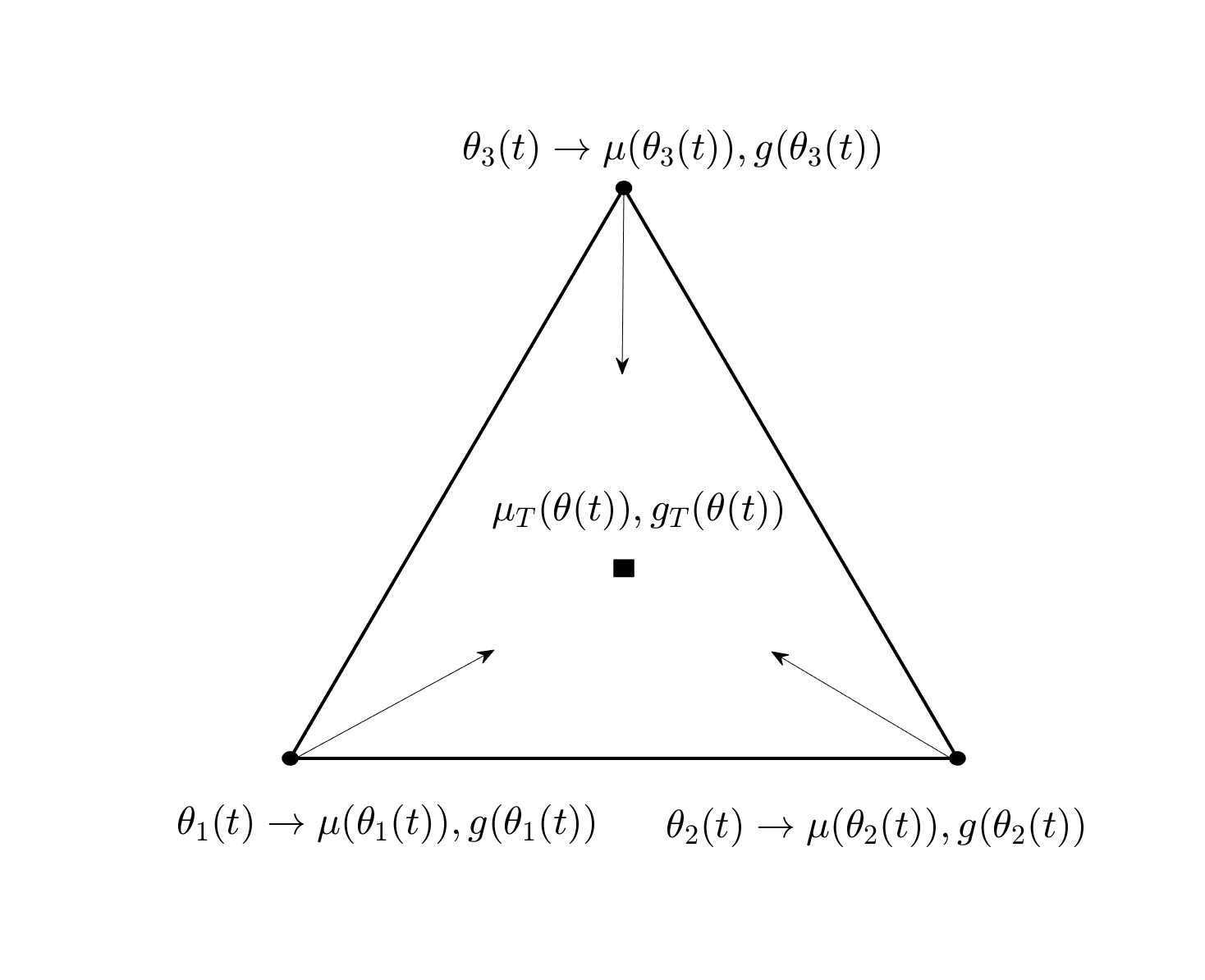}
\caption{\small Weighting technique for $\mu$ and $g$ in $T\in \mathcal{T}^h$. Here, $\theta_i$ represents the value of variable $\theta$ at each vertex $i$ of $T$ and $\mu(\theta_i)$ and $g(\theta_i)$ are the values of the functions $\mu$ and $g$ in $\theta_i$, respectively.}\label{fig:muTgT}
\end{center}
\end{figure}

Let us now explain the matrices $\mathbf{A}_\mu^h(\vec{\theta}(t))$ and $\mathbf{Q}_g^h(\vec{\theta}(t))$. These matrices are defined as follows
\[
\mathbf{A}_\mu^h(\vec{\theta})_{ij}:= \sum_{T\in\mathcal{T}^h} \int_T \mu_T(\vec{\theta}(t)) (\mathcal{E}\boldsymbol{\varphi}_i : \mathcal{E}\boldsymbol{\varphi}_j)\,dx\,\,\mbox{  and  }\,\, \mathbf{Q}_g^h(\vec{\theta})_{ij}:= \sum_{T\in\mathcal{T}^h} \int_T g_T(\vec{\theta}(t)) (\boldsymbol{\psi}_i : \mathcal{E}\boldsymbol{\varphi}_j)\,dx, 
\]
where, as before,  $\boldsymbol{\varphi}_j$, $j=1,\ldots,2d$ are the basis functions of $\mathbf{V}^h$ and $\boldsymbol{\psi}_j$, $j=1,\ldots,4m$, are the basis functions of $\mathbf{W}^h$. Here, $\mu_T(\vec{\theta}(t))$ represents the application of a simple weighting operator to obtain the value of the function $\mu$ in the gravity center of each triangle $T\in \mathcal{T}^h$ from the values of $\mu(\vec{\theta}(t))$ at each vertex of $T$, a.e. in $(0,T_f)$. Since we are working in a uniform mesh, in this contribution we use a simple average operator. The same applies for $g(\vec{\theta}(t))$. This is depicted in Figure \ref{fig:muTgT}.

Finally, let us discuss the space discretization for the dissipation term 
$$\int_\Omega \left[\mu(\theta(t))\|\mathcal{E}\mathbf{u}_\gamma(t)\|^2 + g(\theta(t))\|\mathcal{E}\mathbf{u}_\gamma(t)\|\right]\phi\,dx.$$
Let us recall that both $g$ and $\mu$ are affine functions on $\theta$ (see \eqref{gmu}), and we consider that $\theta$ is a function with values in the interval $[0,1]$. Therefore, we decompose the two terms in the integral above as follows.
\[
\int_\Omega \mu(\theta(t))\|\mathcal{E}\mathbf{u}_\gamma(t)\|^2 \phi\,dx=\delta_\mu \int_\Omega  \|\mathcal{E}\mathbf{u}_\gamma(t)\|^2\, \theta(t) \phi\,dx + \mu_0 \int_\Omega \|\mathcal{E}\mathbf{u}_\gamma(t)\|^2 \phi\,dx,
\]
and
\[
\int_\Omega g(\theta(t))\|\mathcal{E}\mathbf{u}_\gamma(t)\| \phi\,dx=\delta_g \int_\Omega  \|\mathcal{E}\mathbf{u}_\gamma(t)\|\, \theta(t) \phi\,dx + g_0 \int_\Omega \|\mathcal{E}\mathbf{u}_\gamma(t)\| \phi\,dx.
\]
Therefore, by following the classical Galerkin method, the discretization of this terms reads as follows
\[
\mathbf{K}_\mu^h(\vec{\mathbf{u}}(t)):=\delta_\mu M^{h,2}(\vec{\mathbf{u}}(t)) \vec{\theta}(t) + \mu_0\Theta^{h,2}(\vec{\mathbf{u}}(t)),
\]
where
\[
M^{h,2}(\vec{\mathbf{u}}(t))_{i,j}:= \sum_{T\in\mathcal{T}^h} \int_T [N^h(\mathcal{E}^h\vec{\mathbf{u}}(t))]_T^2\, (\phi_i\,\phi'_j)\,dx,\,\mbox{ and }\,\, (\Theta^{h,2})_i := \frac{1}{6}|T| [N^h(\mathcal{E}^h\vec{\mathbf{u}}(t))]_T^2,
\]
for all $i,j=1,\ldots,d$. Here $[N^h(\mathcal{E}^h\vec{\mathbf{u}}(t))]_T^2$ stands for approximated value of $\|\mathcal{E}\mathbf{u}\|^2$ at each $T\in \mathcal{T}^h$, and $|T|$ is the measure of the given triangle. The components of vector $\Theta^{h,2}$ are constructed by following the approximation given in \cite[Sec. 6]{cc50}.

By using the same argumentation, we have that
\[
\mathbf{K}_g^h(\vec{\mathbf{u}}(t)):= \delta_g M^{h,1}(\vec{\mathbf{u}}(t)) \vec{\theta}(t) + g_0\Theta^{h,1}(\vec{\mathbf{u}}(t)),
\]
where
\[
M^{h,1}(\vec{\mathbf{u}}(t))_{i,j}:= \sum_{T\in\mathcal{T}^h} \int_T [N^h(\mathcal{E}^h\vec{\mathbf{u}}(t))]_T\, (\phi_i\,\phi'_j)\,dx,\,\mbox{ and }\,\, (\Theta^{h,1})_i := \frac{1}{6}|T| [N^h(\mathcal{E}^h\vec{\mathbf{u}}(t))]_T,
\]
for all $i,j=1,\ldots,d$. Here $[N^h(\mathcal{E}^h\vec{\mathbf{u}}(t))]_T$ stands for approximated value of $\|\mathcal{E}\mathbf{u}\|$ at each $T\in \mathcal{T}^h$.

\begin{remark}
By construction, matrices $\mathbf{A}^h_\mu(\vec{\mathbf{u}}(t))$,  $\mathbf{Q}^h_g(\vec{\mathbf{u}}(t))$, $M^{h,2}_\mu(\vec{\mathbf{u}}(t))$ and $M^{h,1}_\mu(\vec{\mathbf{u}}(t))$ are weighted stiffness and mass matrices, respectively. In the particular case of the dissipation term, it is possible to obtain the respective weighted mass matrices because of the specific form of the functions $\mu$ and $g$. Other class of functions need a further analysis. 
\end{remark}

\subsection{Time advancing method}
Usual time advancing techniques, such as the one-step $\theta$-method (see \cite{QuarVall}), applied to Navier–Stokes-type equations lead to the numerical solution of nonlinear and convective systems of algebraic equations, which change in every time step. This fact provokes an increase of the computational cost. One approach to avoid this issue is to use operator splitting techniques (see \cite{Sanchez}). However, the use of such methods needs the solution of several other systems to construct the solution of the problem. Another approach is the use of semi implicit methods. One important characteristic shared by a class of these methods is that they lead to Stokes-type matrices with no convective term active in every time step (see \cite{baker,DlRG1,DlRG2}). In this article, we focus on a semi-implicit method proposed in \cite{DlRG1} for the Bingham flow, based on the second-order backward differentiation formulae (BDF2) and on the introduction of a lag-operator. This approximation enjoys the same kind of property: a system whose associated matrix is a Stokes-type one and does not change in every time step.

We discuss separately the time advancing for the flow and the energy equation. Lets us start by the flow equation. By following \cite[Sec. 4.B.]{DlRG1}, we formulate the BDF2 approximation for \eqref{femflow}, as follows: given a vector $\vec{\theta}\in\Real^{d}$, at each time level $t_{k+1}=(k+1)\delta_t$, for $k=0,\ldots, N-1$, solve the system
\begin{equation}\label{flowfd}
\begin{array}{ccc}
\left(\frac{3}{2 \delta_t}\mathbf{M}^h + \mathbf{A}_\mu^h(\vec{\theta})\right)\vec{\mathbf{u}}_{k+1} + B^h \,\vec{p}_{k+2}+ \mathbf{Q}_g^h(\vec{\theta})\,\vec{\mathbf{q}}_{k+2} = \vec{\mathbf{F}}_{k+2}\vspace{0.2cm}\\ -(B^h)^\top \,\vec{\mathbf{u}}_{k+2}=0\vspace{0.2cm}\\  \max\left[G_T(\vec{\theta}), \gamma N^h(\mathcal{E}^h\vec{\mathbf{u}}_{k+2} )\right]\star\vec{\mathbf{q}}_{k+2} = \gamma \textrm{diag}(G_T (\vec{\theta}))\mathcal{E}^h\vec{\mathbf{u}}_{k+2},
\end{array}
\end{equation}
where $\vec{\mathbf{u}}_k$ represents the approximation of $\vec{\mathbf{u}}(t_k)$. The right hand side is given by
\[
\vec{\mathbf{F}}:= \vec{\mathbf{f}}_{k+2} - \mathbf{C}^h(\Lambda(\vec{\mathbf{u}_k}))\Lambda(\vec{\mathbf{u}_k}) + \frac{2}{\delta_t}\mathbf{M}^h\vec{\mathbf{u}}_{k+1}- \frac{1}{2 \delta_t}\mathbf{M}^h\vec{\mathbf{u}}_k.
\]
Here, $\Lambda(\vec{\mathbf{u}}_k)$ stands for the lag operator and it is defined by
\[
\Lambda(\vec{\mathbf{u}}_k):= 2\vec{\mathbf{u}}_{k+1} - \vec{\mathbf{u}}_k.
\]
The BDF2 scheme combined with the lag operator allows us to approximate the convection matrix with information given by the function in the two previous time steps, which implies that this matrix can be moved to the right hand side of the system. 

The initialization of this scheme is performed as follows: having the discretized initial condition $\vec{\mathbf{u}}_0$, we calculate two intermediate steps $\vec{\mathbf{u}}_{2/3}$ and $\vec{\mathbf{u}}_{4/3}$, by applying  consecutively  the backward Euler method. This process leads us to the following systems (see \cite[p. 370]{baker} and \cite[p. 13]{DlRG1})
\begin{equation*}
\begin{array}{ccc}
\left(\frac{3}{2 \delta_t}\mathbf{M}^h + \mathbf{A}_\mu^h(\vec{\theta})\right)\vec{\mathbf{u}}_{2/3} + B^h \,\vec{p}_{2/3}+ \mathbf{Q}_g^h(\vec{\theta})\,\vec{\mathbf{q}}_{2/3} =\vec{\mathbf{f}}_{2/3} - \mathbf{C}^h(\Lambda(\vec{\mathbf{u}_0}))\Lambda(\vec{\mathbf{u}_0}) + \mathbf{M}^h\vec{\mathbf{u}}_{0}\vspace{0.2cm}\\ -(B^h)^\top \,\vec{\mathbf{u}}_{2/3}=0\vspace{0.2cm}\\  \max\left[G_T(\vec{\theta}), \gamma N^h(\mathcal{E}^h\vec{\mathbf{u}}_{2/3} )\right]\star\vec{\mathbf{q}}_{2/3} = \gamma \textrm{diag}(G_T (\vec{\theta}))\mathcal{E}^h\vec{\mathbf{u}}_{2/3},
\end{array}
\end{equation*}
and
\begin{equation*}
\begin{array}{ccc}
\left(\frac{3}{2 \delta_t}\mathbf{M}^h + \mathbf{A}_\mu^h(\vec{\theta})\right)\vec{\mathbf{u}}_{4/3} + B^h \,\vec{p}_{4/3}+ \mathbf{Q}_g^h(\vec{\theta})\,\vec{\mathbf{q}}_{4/3} =\vec{\mathbf{f}}_{4/3} - \mathbf{C}^h(\Lambda(\vec{\mathbf{u}_0}))\Lambda(\vec{\mathbf{u}_0}) + \mathbf{M}^h\vec{\mathbf{u}}_{2/3}\vspace{0.2cm}\\ -(B^h)^\top \,\vec{\mathbf{u}}_{4/3}=0\vspace{0.2cm}\\  \max\left[G_T(\vec{\theta}), \gamma N^h(\mathcal{E}^h\vec{\mathbf{u}}_{4/3} )\right]\star\vec{\mathbf{q}}_{4/3} = \gamma \textrm{diag}(G_T (\vec{\theta}))\mathcal{E}^h\vec{\mathbf{u}}_{4/3}.
\end{array}
\end{equation*}
Finally, we set $\vec{\mathbf{u}}_1:= \frac{1}{2}(\vec{\mathbf{u}}_{2/3}+\vec{\mathbf{u}}_{4/3})$. 

\begin{remark}\label{rem:dt}
In  \cite[p. 371]{baker} it is proved that this time-advancing method guarantees a second order approximation in time, if the algorithm is initialized with the method explained above and if $\delta_t\leq C h^{4/5}$, where $h$ is the size of the space mesh and $C>0$ is a given arbitrary constant. This result holds for sufficiently slow flows, \textit{i.e.}, for flows in which the velocity is bounded (see \cite[(2.35)]{baker}). Consequently, the method is useful and efficient for the kind of flows that we are numerically analysing in this contribution. However, for flows in which the advective terms dominate over the diffusive ones, the uniform boundedness of the velocity can be lost, provoking a deterioration of the second order approximation. This fact will be studied in future contributions.
\end{remark}

Let us now focus on the energy equation. Given $\vec{\mathbf{u}}\in \Real^{2d}$, we use the same approach based on the BDF2 method and the lag operator to obtain the following system at each time level $t_{k+1}=(k+1)\delta_t$, for $k=0,\ldots, N-1$.
\begin{equation}\label{enfd}
\begin{array}{lcc}
\left(\frac{3C_p}{2\delta_t}M^h + \kappa A^h\right)\theta_{k+2} + \alpha M^h\theta_{k+2} +C_p \beta M^h_\Gamma \theta_{k+2} - \left[\delta_\mu M^{h,2}(\vec{\mathbf{u}})+ \delta_g M^{h,1}(\vec{\mathbf{u}}) \right] \vec{\theta}_{k+2}= \vec{F}_{k+2},
\end{array}
\end{equation}
where $\vec{\theta}_k$ represents the approximation of $\vec{\theta}(t_k)$, and the right hand side is given by
\[
\vec{F}_{k+2}:=C^h(\Lambda(\vec{\mathbf{u}}_k))\Lambda(\vec{\theta}_k) +\frac{2C_p}{\delta_t} M^h\theta_{k+1} - \frac{C_p}{2\delta_t} M^h\vec{\theta}_k+\mu_0\Theta^{h,2}(\vec{\mathbf{u}}) + g_0\Theta^{h,1}(\vec{\mathbf{u}}).
\]
The initialization follows from the discretized initial condition and the calculation of $\theta_1$ by using a similar approximation as the one used for the flow equation. Due to the use of the lag operator, this discretized equation is a linear equation, which does not need the application of Newton type methods.

\section[Combined BDF2-SSN Algorithm]{Combined BDF2-SSN Algorithm}\label{sec:algo}
In this section, we discuss the combined BDF2-Semismooth Newton Algorithm to solve the system \eqref{flowfd}-\eqref{enfd}. We propose a sequential algorithm, which means that we solve the flow equation with a given temperature field, and then we update the temperature with the resulting velocity field. Due to the discretization scheme, the fully discretized energy equation \eqref{enfd} is now a linear equation, which implies that its solution depends only on the nonsigularity of the system matrix. On the other hand, the flow equation requires a nonlinear algorithm to be solved. As stated before, we propose a semismooth Newton algorithm to find the numerical solution of system \eqref{flowfd}. 

Let us show the proposed combined BDF2-SSN algorithm for the non-isothermal Bingham flow with temperature dependent parameters in the time interval $[0,T_f]$.

\begin{algo}(non-isothermal BDF2-SSN)\label{algo:full}
\begin{enumerate}
\item \textit{Initialization}: Given $\vec{\theta}_0$ and $\vec{\mathbf{u}}_0$, calculate $\vec{\mathbf{u}}_1:=\frac{1}{2}(\vec{\mathbf{u}}_{2/3}+\vec{\mathbf{u}}_{4/3})$. Introduce $\vec{\mathbf{u}}_1$ in \eqref{enfd}, calculate $\vec{\theta}_1$, and set $k:=0$.
\item For $k=0,\ldots \mathcal{N}-2$ do
\begin{enumerate}
\item \textit{Flow update}: Given $\vec{\theta}_{k+1}$, obtain $\vec{\mathbf{u}}_{k+2}$ by applying the SSN algorithm to the following system
\begin{equation}\label{flowalgo}
\begin{array}{ccc}
\left(\frac{3}{2 \delta_t}\mathbf{M}^h + \mathbf{A}_\mu^h(\vec{\theta}_{k+1})\right)\vec{\mathbf{u}}_{k+2} + B^h \,\vec{p}_{k+2}+ \mathbf{Q}_g^h(\vec{\theta}_{k+1})\,\vec{\mathbf{q}}_{k+2} = \vec{\mathbf{F}}_{k+2}\vspace{0.2cm}\\ -(B^h)^\top \,\vec{\mathbf{u}}_{k+2}=0\vspace{0.2cm}\\  \max\left[G_T(\vec{\theta}_{k+1}), \gamma N^h(\mathcal{E}^h\vec{\mathbf{u}}_{k+2} )\right]\star\vec{\mathbf{q}}_{k+2} = \gamma \textrm{diag}(G_T (\vec{\theta}_{k+1}))\mathcal{E}^h\vec{\mathbf{u}}_{k+2},
\end{array}
\end{equation}

\item \textit{Temperature update}: Use the calculated $\vec{\mathbf{u}}_{k+2}$ to obtain $\vec{\theta}_{k+2}$ by solving
\begin{equation}\label{enalgo}
\begin{array}{lcc}
\left(\frac{3C_p}{2\delta_t}M^h + \kappa A^h + \alpha M^h+C_p \beta M^h_\Gamma \right)\theta_{k+2}  \vspace{0.2cm}\\\hspace{3cm}- \left[\delta_\mu M^{h,2}(\vec{\mathbf{u}}_{k+2})+ \delta_g M^{h,1}(\vec{\mathbf{u}}_{k+2}) \right] \vec{\theta}_{k+2}= \vec{F}_{k+2},
\end{array}
\end{equation}
\end{enumerate}
\end{enumerate}
\end{algo}

\begin{remark}\label{rem:conven}
Note that \eqref{enalgo} stands for the fully discretized version of \eqref{enfd}, rewritten in order to emphasize the fact that the solution of this system depends only on the non singularity of the system matrices. This fact directly follows from the positive definiteness of mass and stiffness matrices \cite[p.148]{QuarVall}. Further, this property is shared by the weighted matrices $M^{h,1}(\vec{\mathbf{u}})$ and $M^{h,2}(\vec{\mathbf{u}})$, since the weights are positive numbers which are constant in every $T\in\mathcal{T}^h$.
\end{remark}

Note that, the Remark \ref{rem:conven} implies that the success of Algorithm \ref{algo:full} depends only on the convergence of the inner SSN algorithm developed to numerically solve \eqref{flowalgo}.

\subsection{Semismooth Newton Algorithm}
The main difficulty regarding system \eqref{flowalgo} is that the functions involved are not necessarily differentiable in the classical sense (Fréchet or Gateaux differentiable). In this section, we propose a semismooth Newton algorithm to numerically solve the system, following ideas from [5, 6]. For the sake of readability, we start by introducing the definition of Newton or slantly differentiable functions.
\begin{defi}
Let $D\subset \Real^\ell$ be an open subset. The mapping $F:D\rightarrow \Real^\ell$ is called Newton or slantly differentiable on the open subset $U\subset D$ if there exists a family of applications $G:U\rightarrow \mathcal{L}(\Real^\ell)$ such that
\[
\underset{h\rightarrow 0}{\lim} \frac{1}{\|h\|_{\Real^\ell}}\| F(x+h) - F(x) - G(x+h)h\|_{\Real^\ell}=0,
\]
for all $x\in U$.
\end{defi}
This concept generalizes the classical Fr\'echet differentiability, and allows us to calculate generalized derivatives of functions such as the norm or the $\max$ functions, which are involved in the system under analysis. In fact, in \cite{HIK}, for instance, it is established that the function $\Real^n \ni x\mapsto \max(g\vec{e},x)$, with $g\geq0$, is slantly differentiable with slantly derivative given by
\[
(m(x))_i:= \left\{
\begin{array}{ccc}
1 &\mbox{if $x_i\geq g$}\vspace{0.2cm}\\ 0 &\mbox{if $x_i< g$}.
\end{array}
\right.
\]
Further, we write the following convergence result from \cite{HIK}.
\begin{prop}
If $x^*$ is a solution of $F(x) = 0$, $F$ is Newton or slantly differentiable in an open neighborhood U containing $x^*$ with generalized derivative G. If $\{\|G(y) ^{-1}\|\,:\,y \in U\}$ is bounded, then the Newton iterations
\[
x_{k+1} =x_k -G(x_k)^{-1}F(x_k)
\]
converge superlinearly to $x^*$, provided that $\|x_0 - x^*\|$ is sufficiently small.
\end{prop}

Suppose that we are given a vector $\vec{\theta}$ and let us rewrite system \eqref{flowfd} as the following operator equation
\begin{equation*}
\begin{pmatrix}
\left(\frac{3}{2 \delta_t}\mathbf{M}^h + \mathbf{A}_\mu^h(\vec{\theta})\right)\vec{\mathbf{u}}_{k+1} + B^h \,\vec{p}_{k+2}+ \mathbf{Q}_g^h(\vec{\theta})\,\vec{\mathbf{q}}_{k+2}- \vec{\mathbf{F}}_{k+2}\vspace{0.2cm}\\ -(B^h)^\top \,\vec{\mathbf{u}}_{k+2}\vspace{0.2cm}\\   \textrm{diag}(\mathbf{m}(\vec{\mathbf{u}}_{k+2},\vec{\theta}))\vec{\mathbf{q}}_{k+2} -\gamma \textrm{diag}(G_T (\vec{\theta}))\mathcal{E}^h\vec{\mathbf{u}}_{k+2},
\end{pmatrix}=\mathbf{0},
\end{equation*}
where $\textrm{diag}(\mathbf{m}(\vec{\mathbf{u}}_{k+2},\vec{\theta}))$ stands for the diagonal matrix with entries given by the components of the vector $\max\left[G_T(\vec{\theta}_{k+1}), \gamma N^h(\mathcal{E}^h\vec{\mathbf{u}}_{k+2} )\right]$.

Next, by following \cite[Sec. 4.C.]{DlRG1} and \cite[Sec. 3.3]{DlRG2}, we calculate the following semismooth Newton step.
\begin{equation}\label{SSNmat}
\begin{array}{lll}
\begin{pmatrix}
\boldsymbol{\Xi}(\vec{\theta}) & B^h & \mathbf{Q}_g^h(\vec{\theta})\vspace{0.2cm}\\ -(B^h)^\top & \mathbf{0} & \mathbf{0}\vspace{0.2cm}\\  \mathcal{S}^h(\vec{\mathbf{u}}_{k+2},\vec{\theta}) & \mathbf{0} & \textrm{diag}(\mathbf{m}(\vec{\mathbf{u}}_{k+2},\vec{\theta}))
\end{pmatrix} \begin{pmatrix}
\delta_{\mathbf{u}}\vspace{0.2cm}\\\delta_{p}\vspace{0.2cm}\\\delta_{\mathbf{q}}
\end{pmatrix}\vspace{0.2cm}\\\hspace{4.5cm}=
\begin{pmatrix}
-\boldsymbol{\Xi}(\vec{\theta})\vec{\mathbf{u}}_{k+1} - B^h \,\vec{p}_{k+2}-\mathbf{Q}_g^h(\vec{\theta})\,\vec{\mathbf{q}}_{k+2}+ \vec{\mathbf{F}}_{k+2}\vspace{0.2cm}\\ (B^h)^\top \,\vec{\mathbf{u}}_{k+2}\vspace{0.2cm}\\  -\textrm{diag}(\mathbf{m}(\vec{\mathbf{u}}_{k+2},\vec{\theta}))\vec{\mathbf{q}}_{k+2} +\gamma \textrm{diag}(G_T (\vec{\theta}))\mathcal{E}^h\vec{\mathbf{u}}_{k+2},
\end{pmatrix},
\end{array}
\end{equation}
where $\boldsymbol{\Xi}:= \left(\frac{3}{2 \delta_t}\mathbf{M}^h + \mathbf{A}_\mu^h(\vec{\theta})\right)$, and $\mathcal{S}^h(\vec{\mathbf{u}}_{k+2},\vec{\theta})$ stands for the slantly derivative of the third equation, which is the nondifferentiable part, and is given by.
\[
\mathcal{S}^h(\vec{\mathbf{u}}_{k+2},\vec{\theta}):=\gamma\left(\chi_{{\mathcal{A}}_{k}}\textrm{diag}(\vec{\mathbf{q}}_{k+2})N^h_{\mathbf{u}}(\mathcal{E}^h\vec{\mathbf{u}}_{k+2})-\textrm{diag}(G_T (\vec{\theta}))\right)\mathcal{E}^h.
\]
Let us discuss this expression with more detail. The function $N^h_{\mathbf{u}}$ stands for the slantly derivative of the norm function $N$ and it is given by
\[
N^h_{\mathbf{w}}:=\textrm{diag}(N^h(\mathbf{w}))^{-1}\begin{pmatrix}
D(w_1)&D(w_2)&D(w_3)&D(w_4)\vspace{0.2cm}\\D(w_1)&D(w_2)&D(w_3)&D(w_4)\vspace{0.2cm}\\D(w_1)&D(w_2)&D(w_3)&D(w_4)\vspace{0.2cm}\\D(w_1)&D(w_2)&D(w_3)&D(w_4)
\end{pmatrix},
\]
for any $\mathbf{w}=(w_1,w_2,w_3,w_4)^\top \in \Real^{4m}$ and $w_i\in\Real^m$. Here, we use the notation $D(w_i):=\textrm{diag}(w_i)$, for all $i=1,\ldots,4$.

Now, let us explain the array $\chi_{{\mathcal{A}}_{k}}$. This expression represents the slantly derivative of the expression $\max\left[G_T(\vec{\theta}_{k+1}), \gamma N^h(\mathcal{E}^h\vec{\mathbf{u}}_{k+2} )\right]$, and it is given by $\chi_{\mathcal{A}_k}:=\textrm{diag}(\vec{\vartheta})$, where

\[
(\vec{\vartheta})_i:= \left\{
\begin{array}{cll}
1 &\mbox{if $(N^h(\mathcal{E}^h\vec{\mathbf{u}}_{k+2}))_i \geq \frac{(G_T(\vec{\theta}_{k+1}))_i}{\gamma}$}\vspace{0.2cm}\\ 0 &\mbox{otherwise}.
\end{array}
\right.
\]
This function plays an important role: it is the numerical characterization of the active set $\mathcal{A}_\gamma$. This set represents an approximation of the yielded regions, \textit{i.e.}, the regions where the material behaves as an incompressible fluid. Therefore, the set $\mathcal{I}_\gamma:=Q\setminus\mathcal{A}_\gamma$ corresponds to an approximation of the unyielded zones in the fluid, which are the rigid regions in the material.

The semismooth Newton method converges locally with superlinear rate \cite[Th. 6.5]{DlRG0}. In practice, however, in the most cases it is not possible to accurately estimate the convergence neighborhood, which can be very small. In such a case, it is mandatory to design a globalization strategy to have a convergent algorithm for arbitrary initialization values. We follow \cite{DlRG1}, where the system matrix in each iteration is slightly modified to guarantee that this matrix is always positive definite. For the system matrix in \eqref{SSNmat}, positive definiteness can only be guaranteed if $N^h(\vec{\mathbf{q}}_{k+2})_i\leq (G_T(\vec{\theta}_{k+1}))_i$, for all $i=1,\ldots,4m$, is satisfied (\cite[Prop. 6.1]{DlRG0}). Since this condition does not automatically hold in each iteration, we project the multiplier $\vec{\mathbf{q}}_{k+2}$ in the feasible set $\{\vec{\mathbf{w}}\in\Real^{4m}\,:\, N^h(\vec{\mathbf{w}})_i\leq(G_T(\vec{\theta}_{k+1}))_i,\,\forall i=1,\ldots,4m \}$, yielding the modified matrix $\widehat{\mathcal{S}}^h$, which is always positive definite. The SSN algorithm is given through the following steps.

\begin{algo}{SSN Algorithm}\label{SSN-Algorithm}
\begin{enumerate}
\item Initialization: Set the initial values $\vec{\mathbf{u}}^0_{k+2}$, $\vec{p}^{\,0}_{k+2}$ and $\vec{\mathbf{q}}^0_{k+2}$, and set $\ell=0$.
\item Active sets estimation: Determine $\chi_{\mathcal{A}_k}$.
\item SSN step: Solve the system \eqref{SSNmat} with $\widehat{\mathcal{S}}^h$.
\item Update: $\vec{\mathbf{u}}^{\ell+1}_{k+2}=\vec{\mathbf{u}}^{\ell}_{k+2}+\delta_\mathbf{u}$, $\vec{p}^{\,\ell+1}_{k+2}=\vec{p}^{\,\ell}_{k+2}+\delta_p$ and $\vec{\mathbf{q}}^{\ell+1}_{k+2}=\vec{\mathbf{q}}^{\ell}_{k+2}+\delta_\mathbf{q}$.
\item Stopping criteria: Verify if $\|\delta\|:=\|\delta_\mathbf{u}\|+\|\delta_p\|+\|\delta_\mathbf{q}\|\leq \epsilon<<1$. If so, stop, otherwise, go to step 2.
\end{enumerate}
\end{algo}
\begin{remark}
The system of equations can be solved with any direct or iterative method, and it is not computationally expensive to obtain its solution. In fact, in \cite[Sec. IV. C]{DlRG1} a decomposition of the matrix is explained, which leads us to the numerical solution of only a $2n\times 2n$ system of linear equations per iteration.
\end{remark}

\begin{remark}\label{rem:convflow}
Note that the matrices involved in system \eqref{SSNmat} depend on a given vector $\vec{\theta}$. However, in the Algorithm \ref{algo:full}, the temperature field is actually a constant vector for every time step, which implies that the SSN algorithm can process it as constant coefficient or given positive constant weights in the system matrices. Furthermore, this coefficients do not modify the properties of these matrices. Particularly, the matrices keep being positive definite and sparse. This fact allows us to conclude that the convergence of the SSN algorithm follows from \cite[Sec. IV.C]{DlRG1} and \cite[Sec. 6]{DlRG0}.
\end{remark}

\section[Computational Results]{Computational Results}\label{sec:comput}
In this section, we present two numerical experiments that show the behavior of the Algorithm non-isothermal BDF2-SSN. First, we analyze a regime in which the viscosity and the yield stress functions increase with temperature, and then a regime in which these functions decrease. We carry out the two experiments in the unit square $\Omega:=(0,1)\times(0,1)$ with $\Gamma:=(0,1)\times\{1\}$ and $\Gamma_0=\partial\Omega\setminus \Gamma$.

Let us first discuss the parameters of the algorithm. Let $\delta^h_k := \|\delta_\mathbf{u}\|_{H^{1,h}} + \|\delta_\mathbf{q}\|_{(L^{2,h})^4} + \|\delta_p\|_{L^{2,h}}$, where the upper index $k$ represents each time step and $H^{1,h} $, $(L^{2,h})^4$ and $L^{2,h}$ stand for the discrete versions of $\mathbf{H}^1(\Omega)$, $(L^2(\Omega))^4$ and $L^2(\Omega)$, respectively. We stop the inner algorithm SSN, at each time step $k$, as soon as $\delta^h_k$ is lower than $\sqrt{\epsilon}$, where $\epsilon$ denotes the machine accuracy ($\epsilon\approx 2.2204e-16$). We fix the regularization parameter $\gamma=10^3$, and  in both experiments we consider the action of a body force given by 
\[
\mathbf{f}(x_1, x_2) := 300(x_2 - 0.5, 0.5 - x_1).
\]
We consider uniform space meshes, whose components have all the same area, and we measure the size of these meshes by the constant radius of the inscribed circumferences of the triangles in the mesh, represented by $h$. Further, we define the time step size as $\delta_t := C(h^{4/5})$, $C>0$. This selection was explained in Remark \ref{rem:dt}.

\begin{figure}
\begin{center}
\includegraphics[width=80mm, height=60mm]{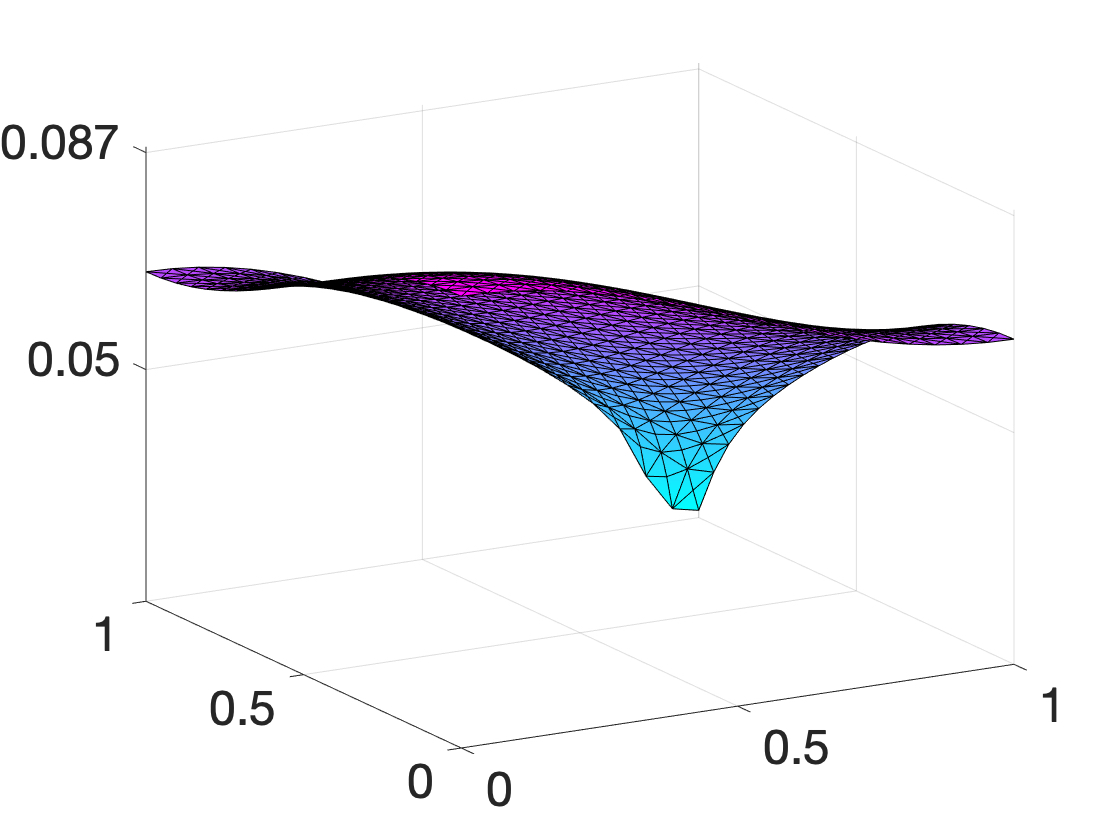}
\end{center}
\caption{\small Initial condition $\theta_0$, with $\beta=15$ and $C_p=1$.}\label{fig:th0exp1}
\end{figure}

The initial conditions in the two experiments are given as follows: for the flow equation we consider that $\mathbf{u}_0=\mathbf{0}$, while for  the energy equation, we consider the onset $\theta_0$ to be the solution of the following elliptic equation
\begin{equation}\label{th0}
\begin{array}{ccl}
-\Delta \theta_0  = \frac{x_1^2}{100}+\frac{x_2^2}{50}+\frac{1}{100},&\mbox{in $\Omega$}\vspace{0.2cm}\\
\frac{\partial \theta_0}{\partial\mathbf{n}}=0,&\mbox{on $\Gamma_0$}\vspace{0.2cm}\\
\kappa\frac{\partial \theta_0}{\partial\mathbf{n}}+C_p\beta\theta_0=0,&\mbox{on $\Gamma$}.
\end{array}
\end{equation}

\begin{figure}
\begin{center}
\includegraphics[width=40mm, height=40mm]{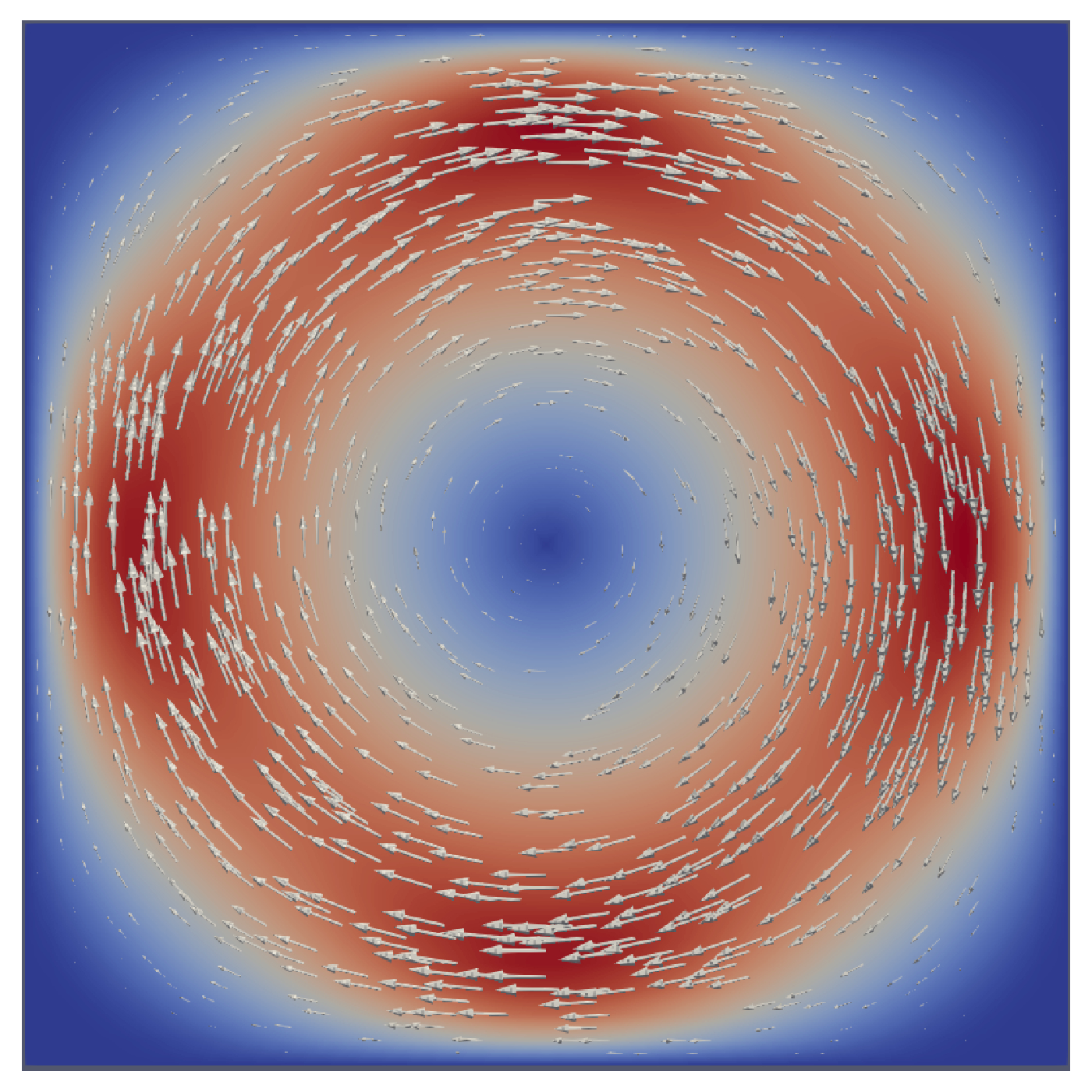}\includegraphics[width=40mm, height=40mm]{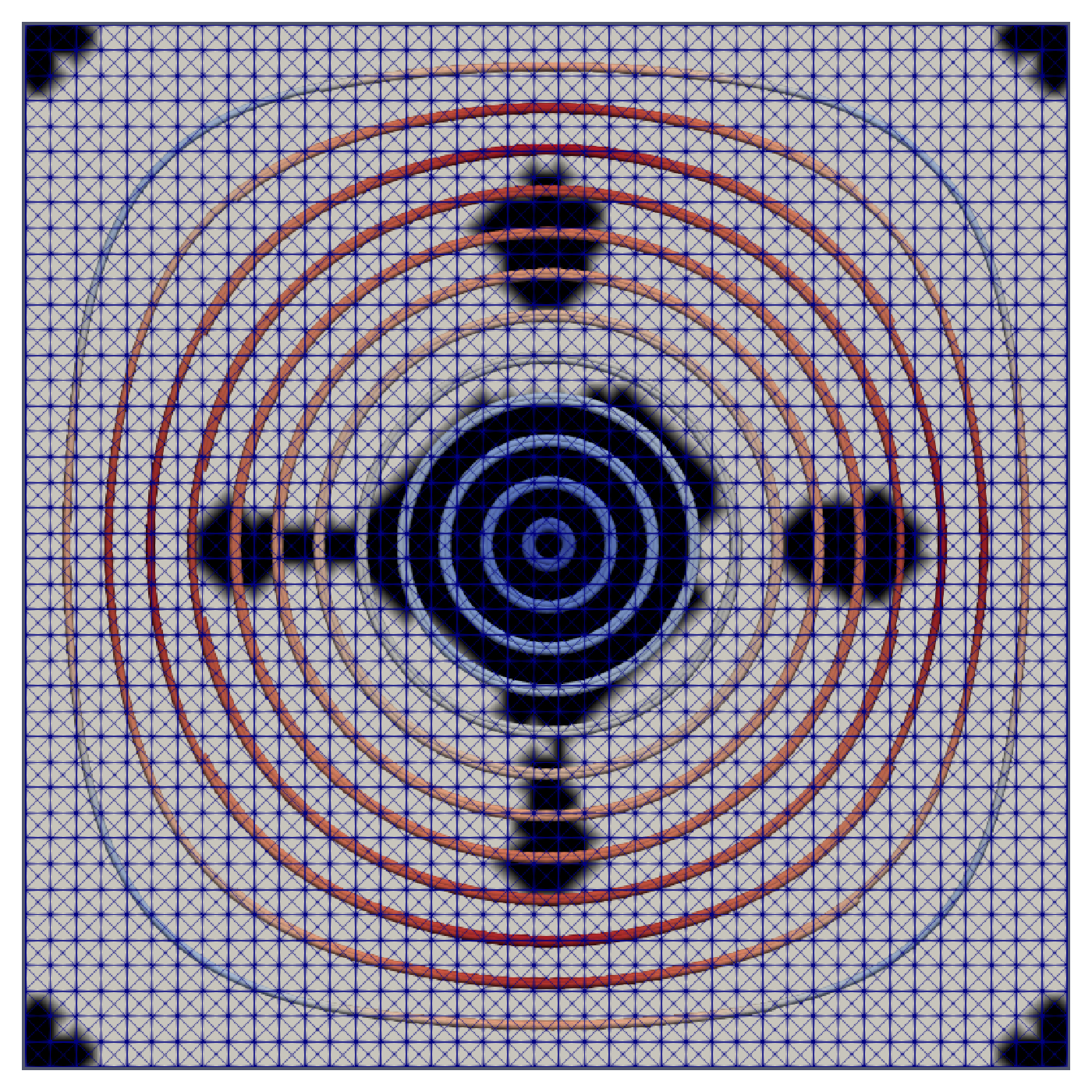}\includegraphics[width=50mm, height=40mm]{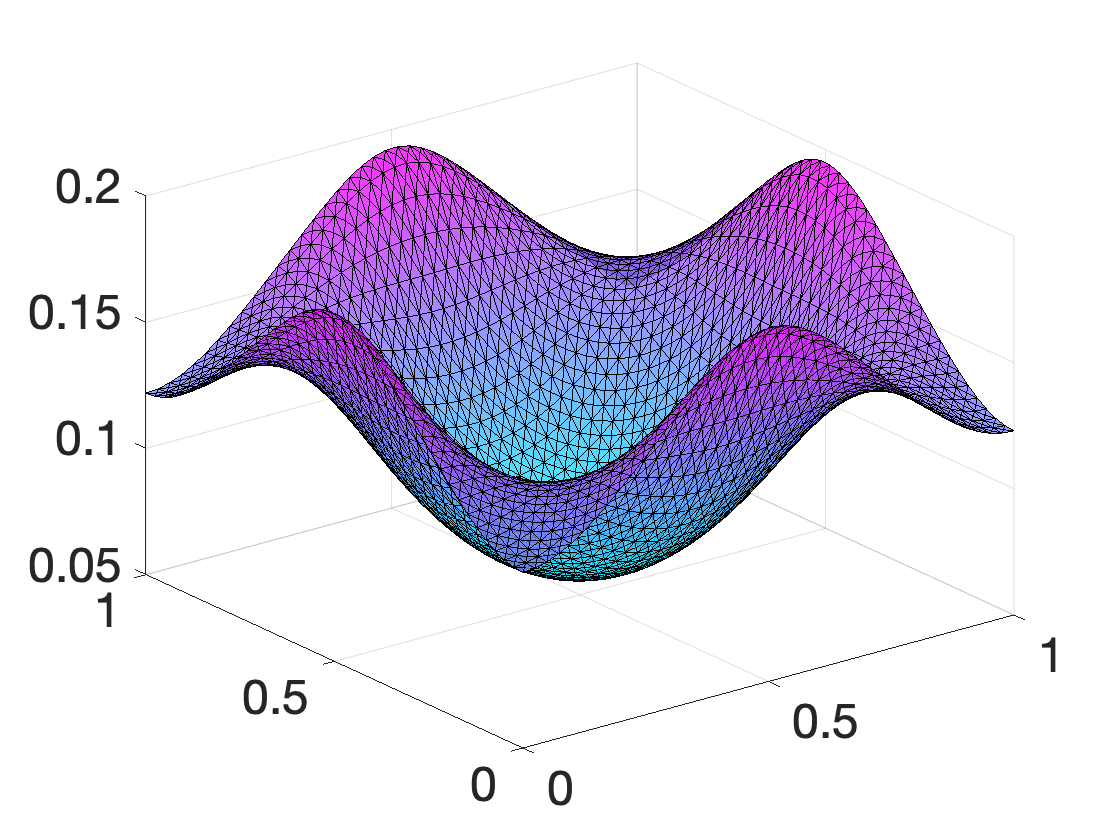}\\
\includegraphics[width=40mm, height=40mm]{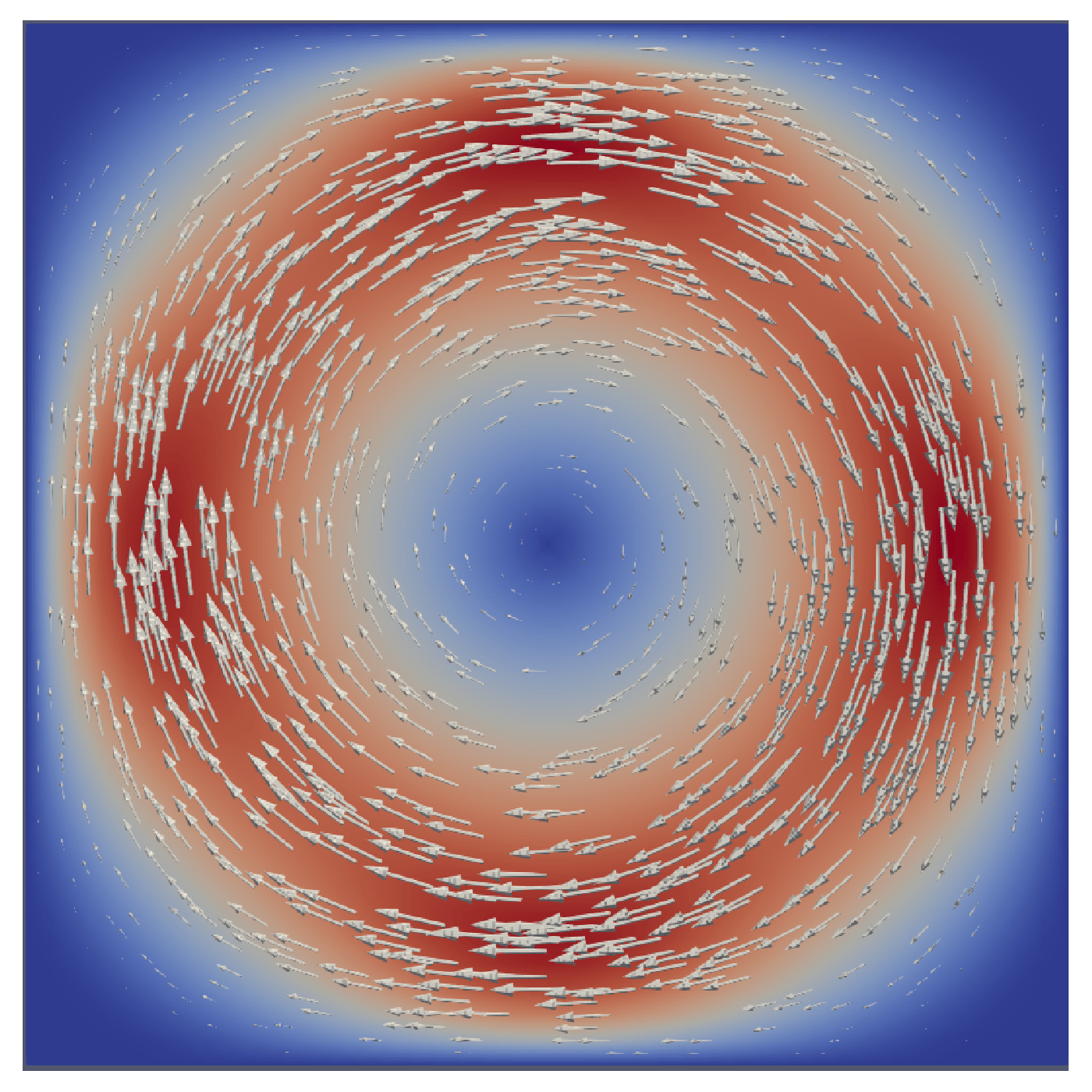}\includegraphics[width=40mm, height=40mm]{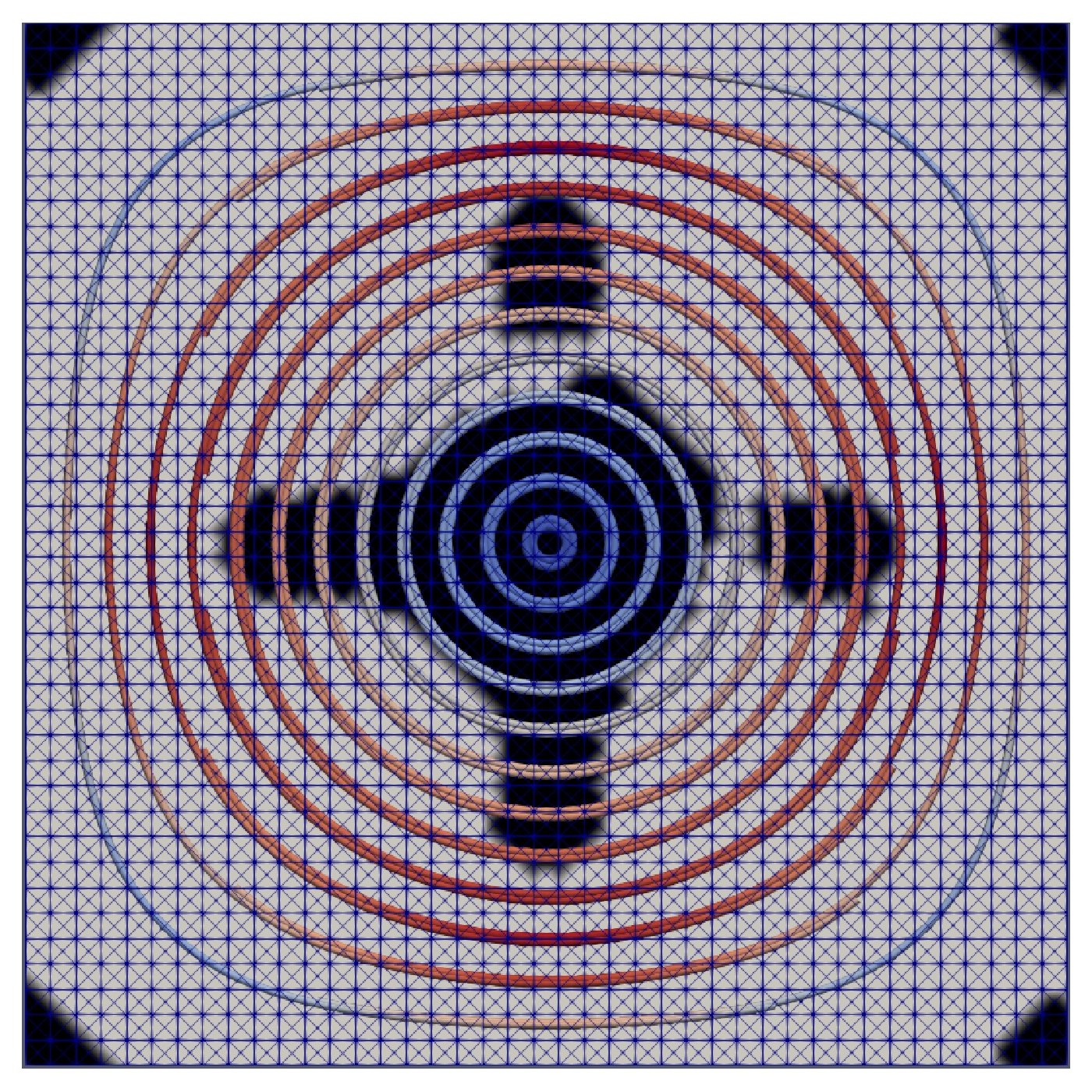}\includegraphics[width=50mm, height=40mm]{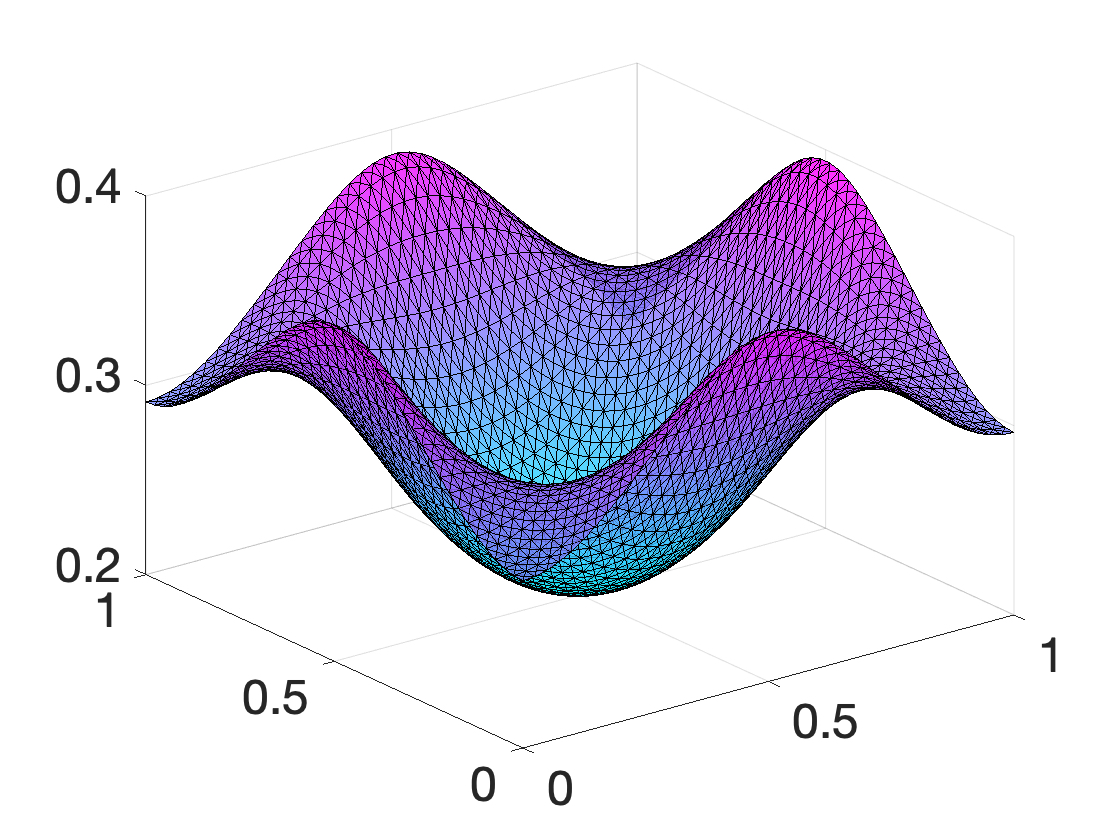}\\
\includegraphics[width=40mm, height=40mm]{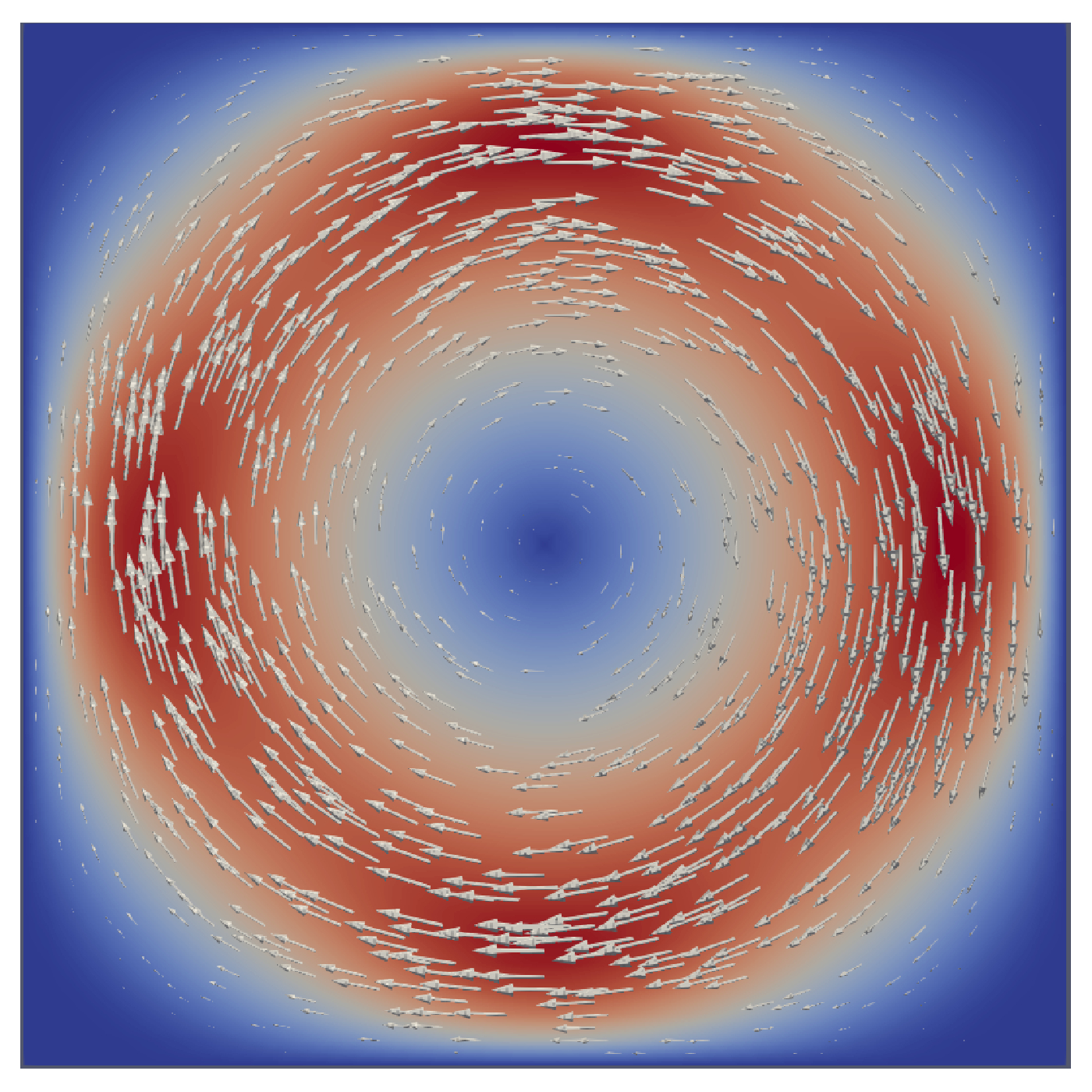}\includegraphics[width=40mm, height=40mm]{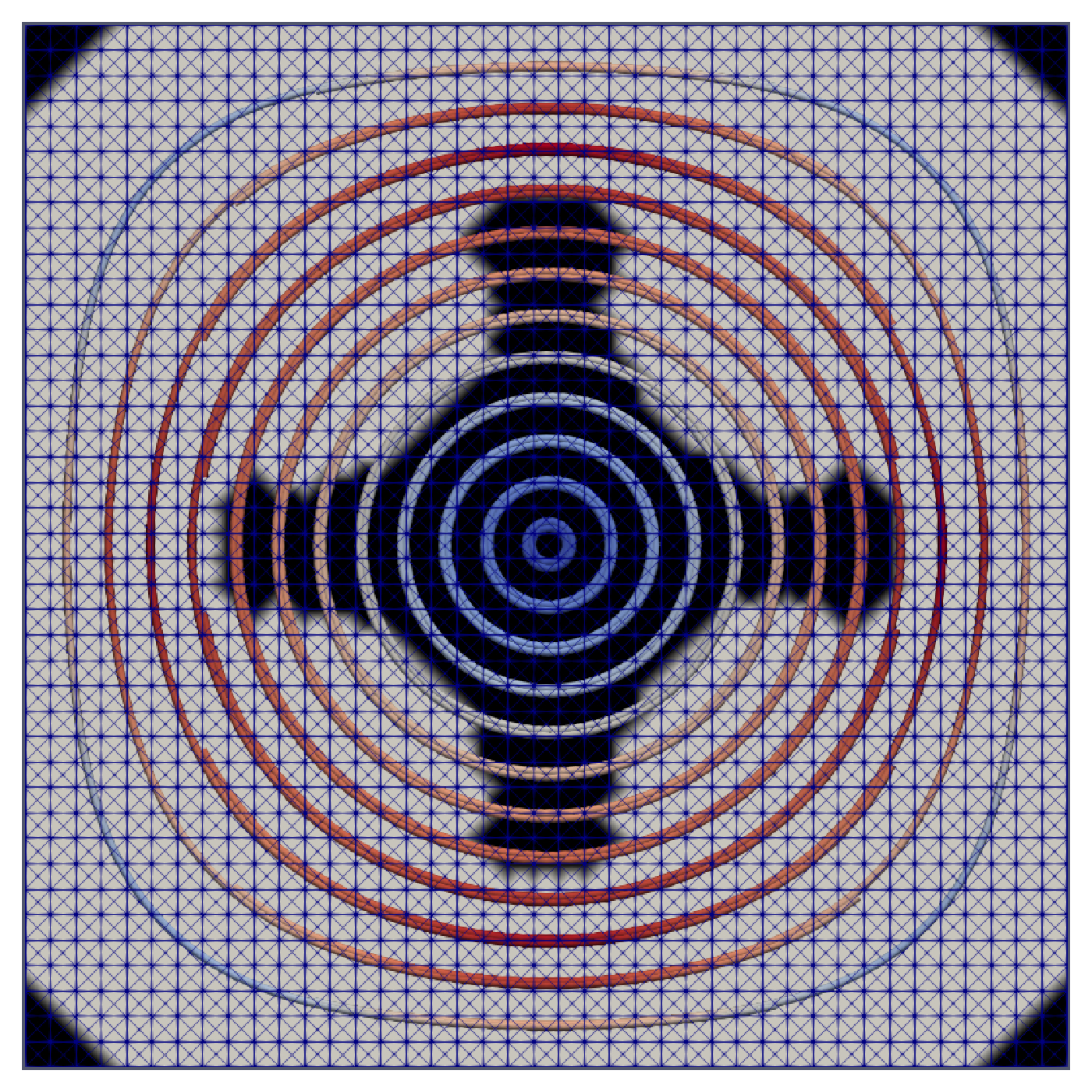}\includegraphics[width=50mm, height=40mm]{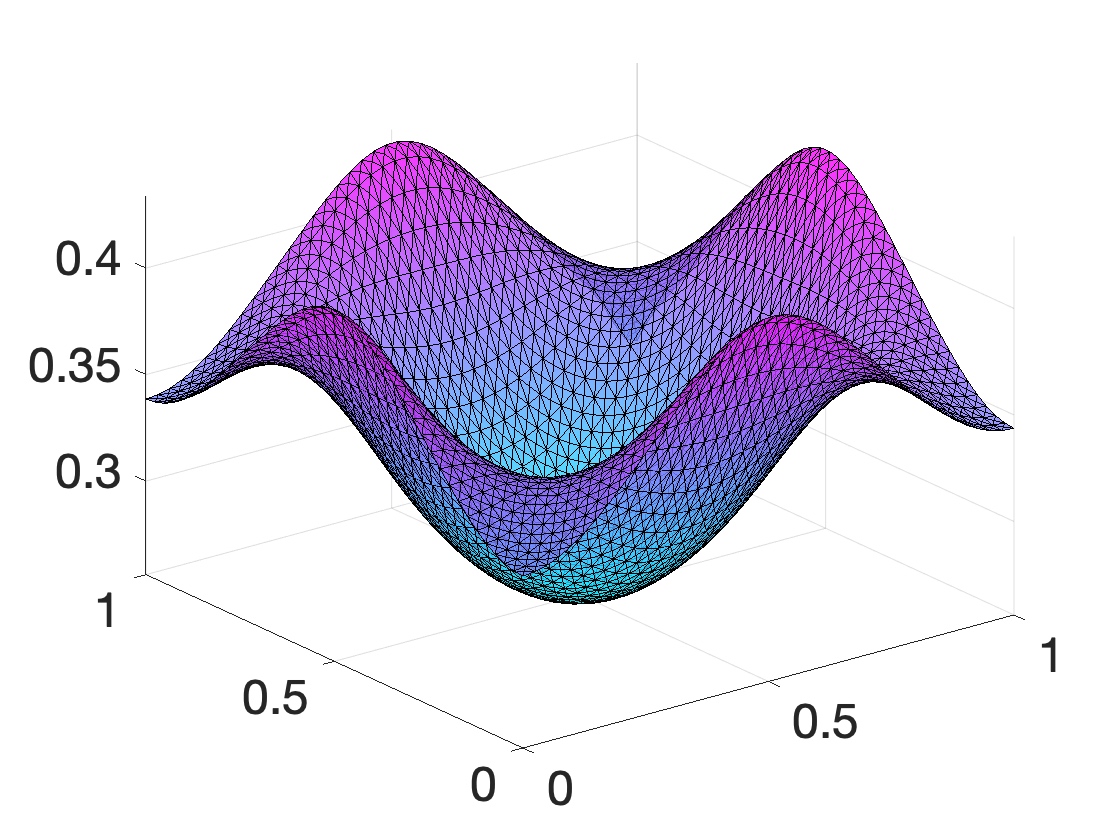}\\
\includegraphics[width=40mm, height=40mm]{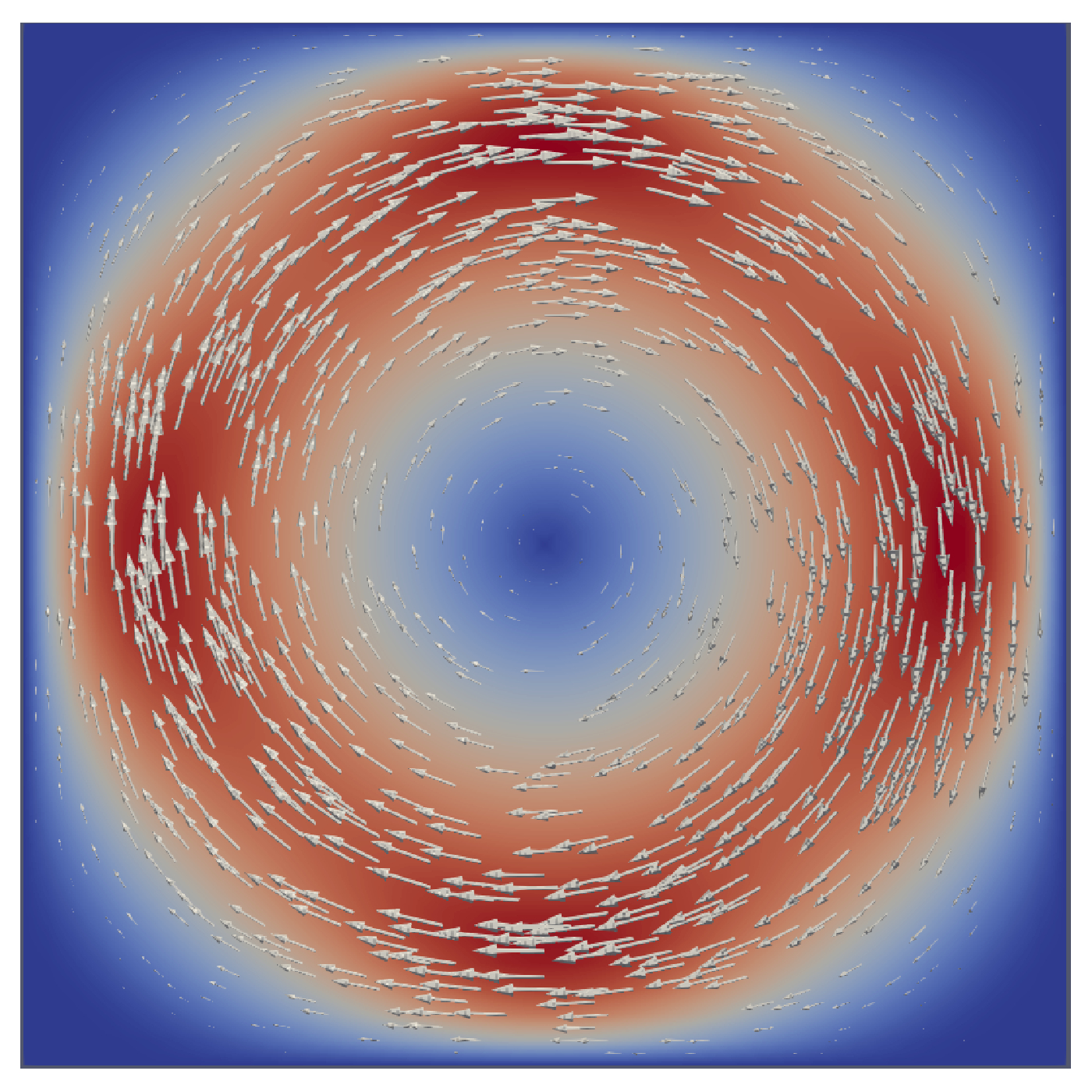}\includegraphics[width=40mm, height=40mm]{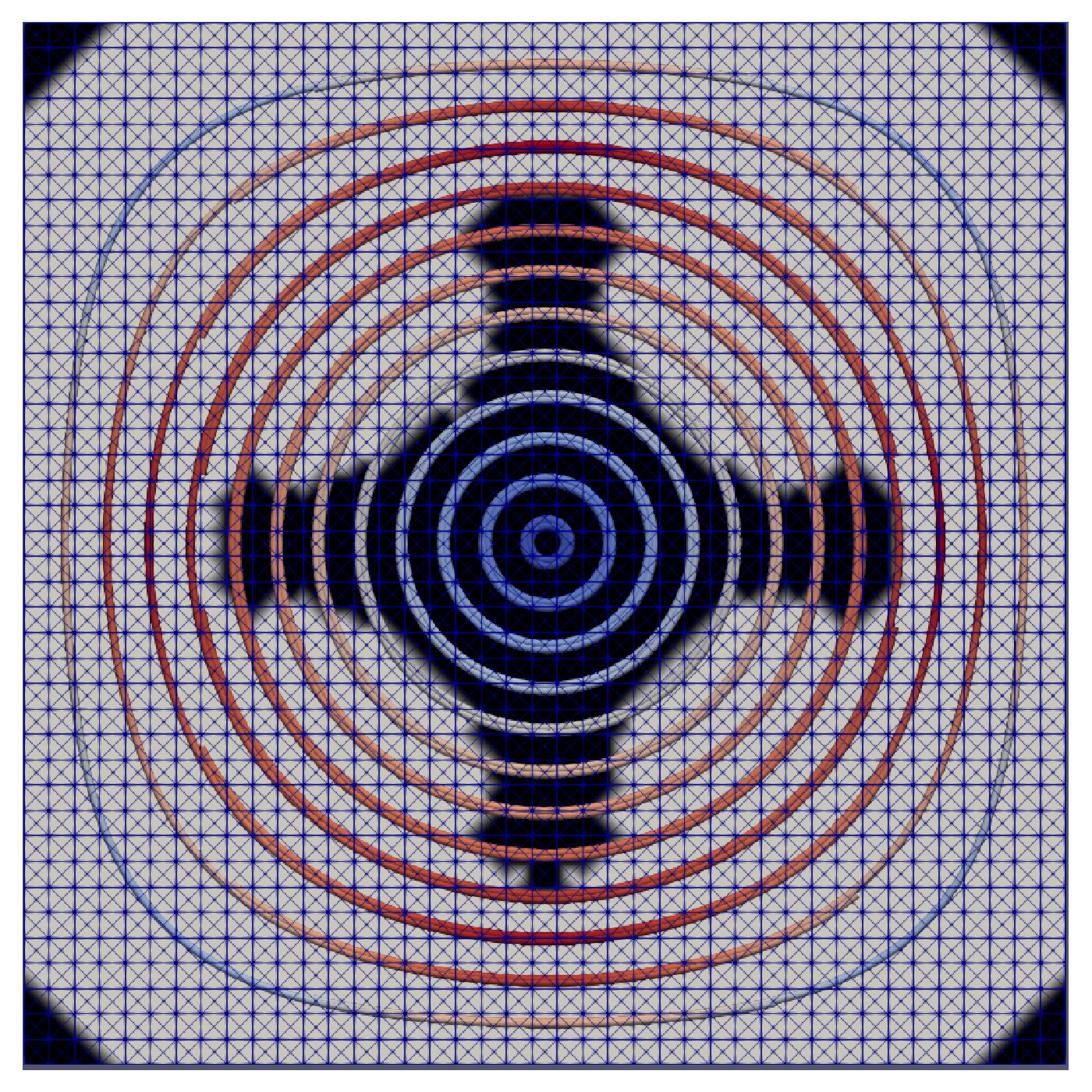}\includegraphics[width=50mm, height=40mm]{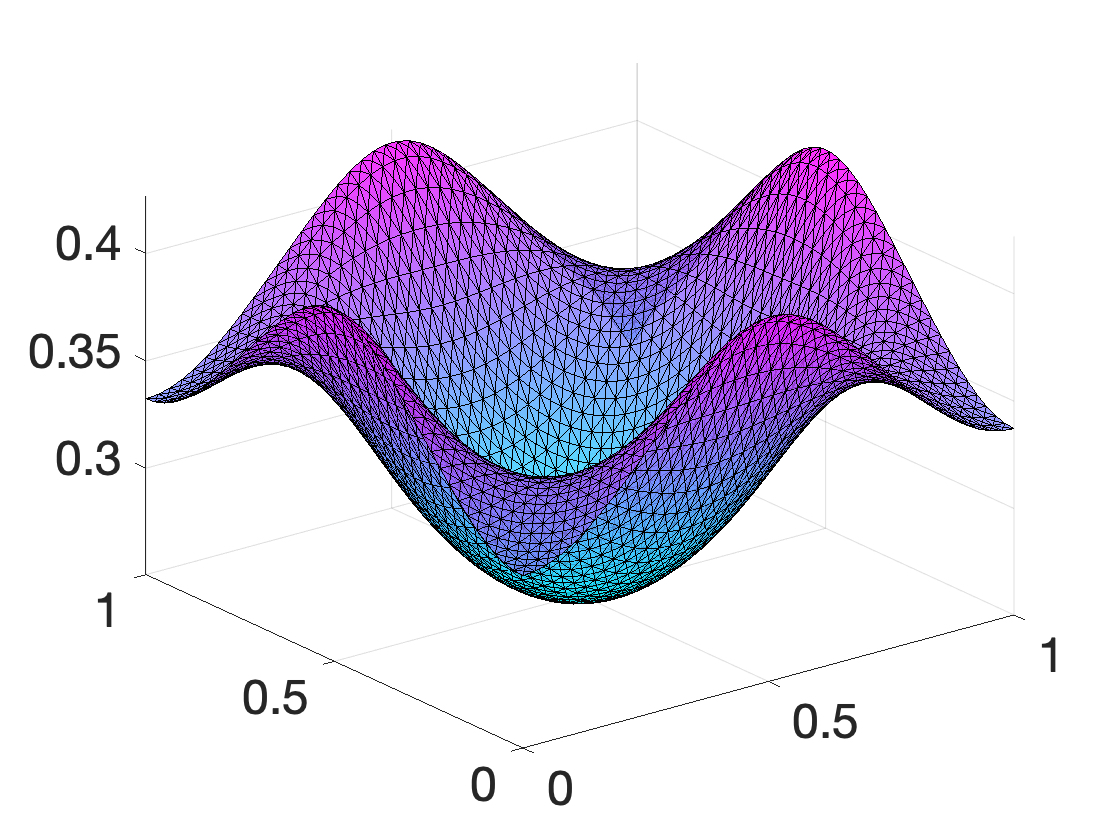}
\end{center}
\caption{\small Experiment 1. Left column: velocity flow. Central column: active and inactive sets with streamlines. Right column: temperature field. For $t=0.015$ (up), $t=0.030$, $t=0.060$ and $T_f=0.12$ (bottom).}\label{fig:flowexp1}
\end{figure}

\subsection{Experiment 1}
In this experiment we calculate the non-isothermal Bingham flow with the following structural parameters: $\mu_0=1$, $\delta_\mu=0.5$, $g_0=10$ and $\delta_g=8$. Thus, we are concerned with the flow in the regime where both the viscosity and the yield stress increase with temperature. We consider that $\alpha=100$, so we admit an external heat sink proportional to the temperature field. We have that $\beta=15$ for the Robin boundary condition imposed in the upper edge of the square geometry, and we assume that $\kappa=10$, $C_p=1$ and $\gamma=10^3$. We consider a mesh size of $h=0.0023$, and we set $\delta_t=0.1*h^{4/5}\approx 0.0018$. The initial condition for the flow and theta are given by $\mathbf{u}_0=\mathbf{0}$ and the solution of equation \eqref{th0}, respectively. Finally, we consider that $T_f=0.12$.

\begin{figure}
\begin{center}
\includegraphics[width=60mm, height=40mm]{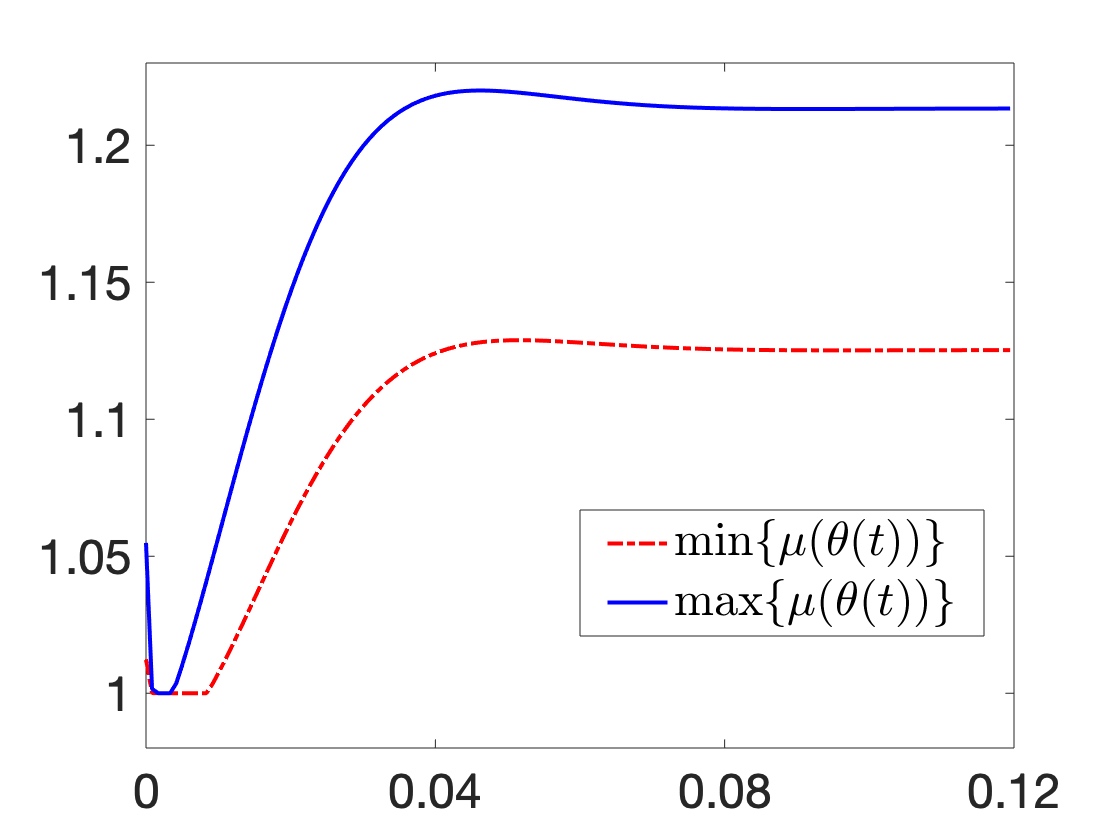}\hspace{1cm}\includegraphics[width=60mm, height=40mm]{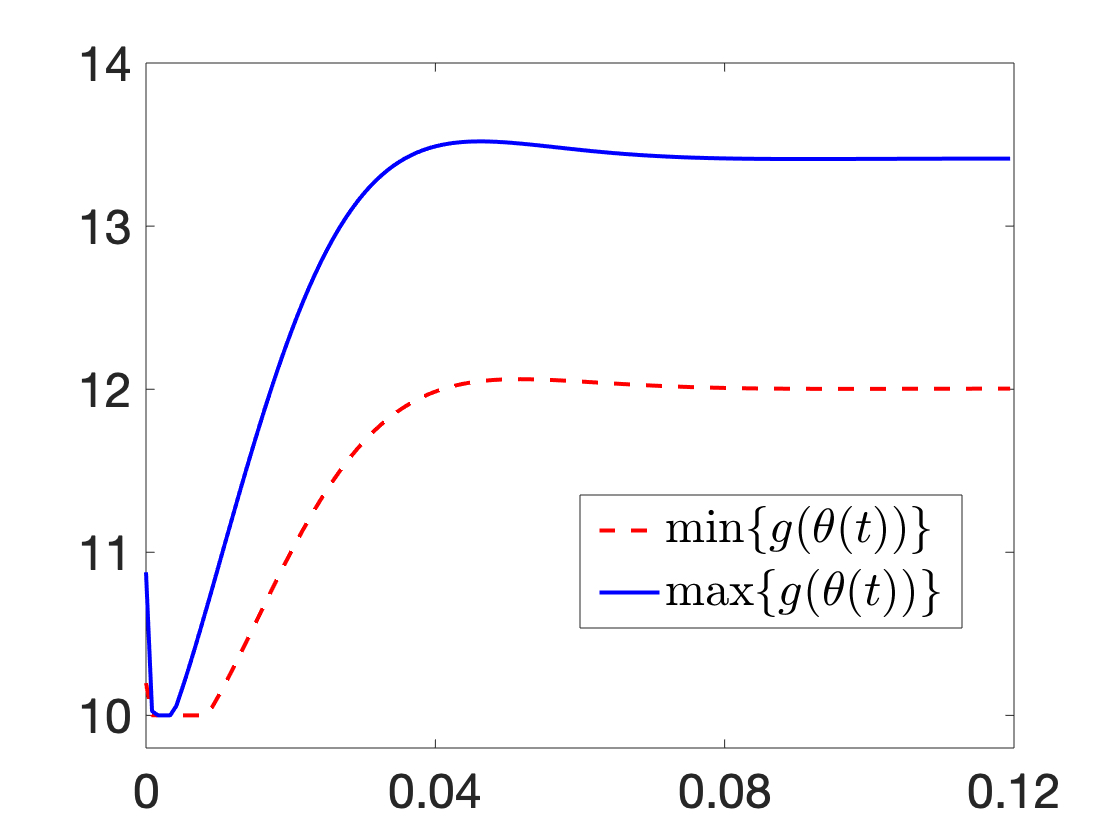}
\end{center}
\caption{\small Experiment 1. Higher and lower values reached by $\mu(t)$ (left) and $g(t)$ (right) in the time interval $[0,T_f]$.}\label{fig:maxmug}
\end{figure}

In Figure \ref{fig:flowexp1}, we show the calculated velocity field, the active and inactive sets, representing the yielded and unyielded
regions of the material, as well as the streamlines of the flow, and the calculated temperature fields, for several instants. The flow behaves as expected: the rigid regions change from the regions given by a constant value of $g$, close to $g_0$ for the initial values, and then evolves to  bigger regions. Further, in Figure \ref{fig:maxmug}, we show the minimum and maximum values reached by $\mu(t)$ and $g(t)$ all along the time interval $[0,T_f]$. In this case, we can observe that the functions $\mu$ and $g$ are bounded all along the time interval, and near $T_f$, the parameters tend to be constant.

\begin{figure}
\begin{center}
\includegraphics[width=60mm, height=40mm]{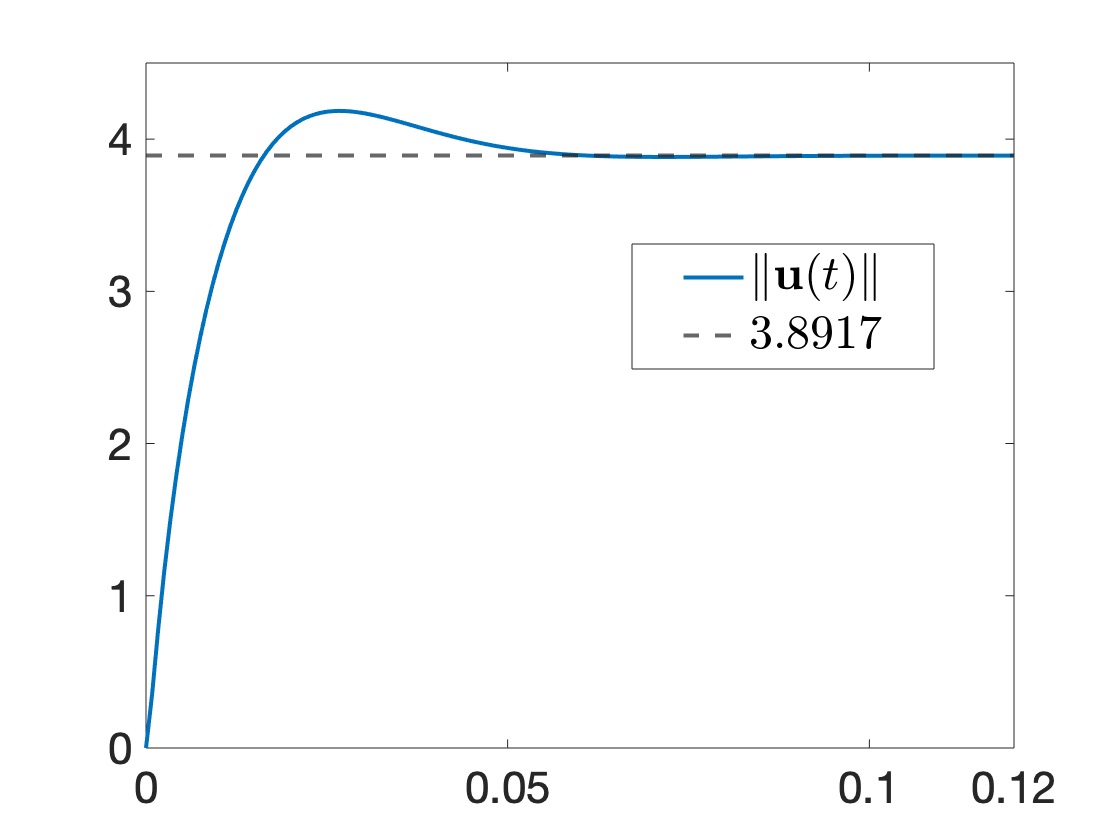}\hspace{1cm}\includegraphics[width=60mm, height=40mm]{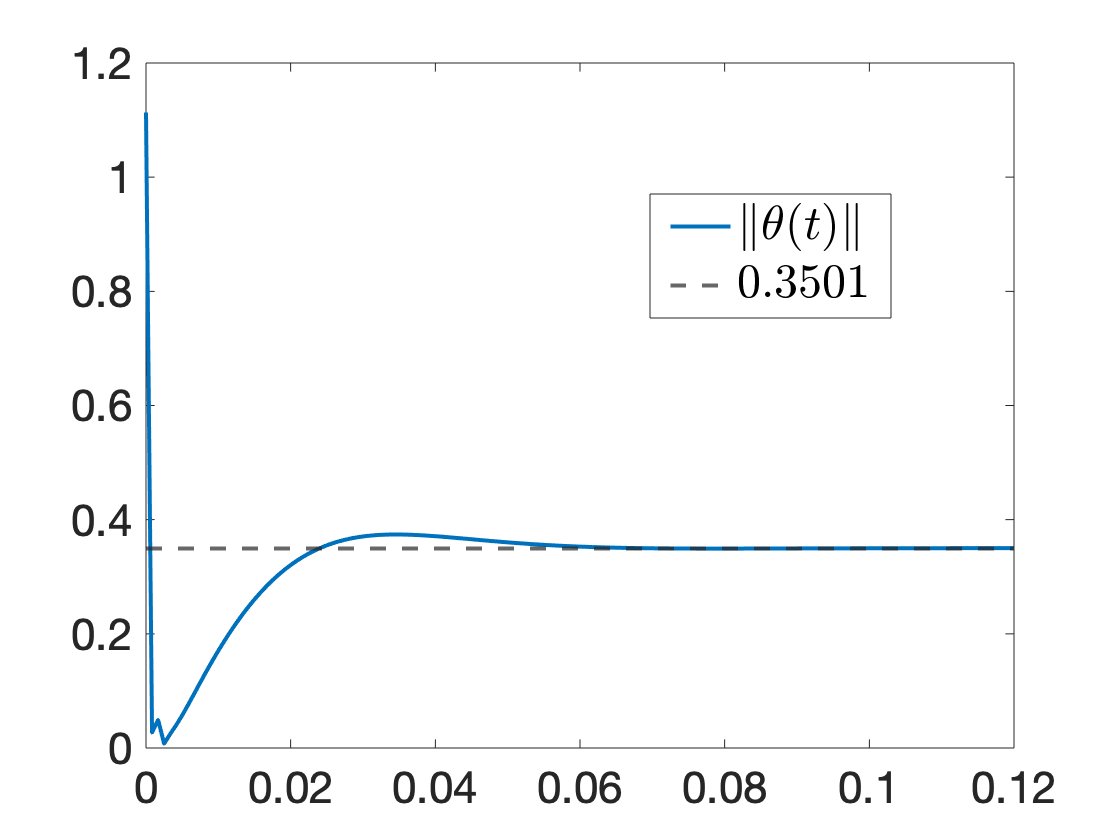}
\end{center}
\caption{\small Experiment 1. Evolution of $\|\mathbf{u}(t)\|_{H^{1,h}}$ (left) and $\|\theta(t)\|_{W^{q,h}}$ (right).}\label{fig:kine}
\end{figure}

\begin{table}[h]\label{tab:1}       
\begin{center}
\begin{tabular}{c c c c c}
\hline
$t$&0.015 & 0.030 & 0.060 & 0.12 \\
\hline
 & 1.1075e-4  & 2.5283e-5 & 0.0040& -\\
$\delta_k^h$ & 5.4368e-6 & 2.0906e-8 & 1.0234e-6 & 2.0329e-5\\
& 4.0606e-9 & 5.3680e-14 & 1.3834e-11 & 4.7995e-11\\\hline
\# it.& 5 & 4 & 3 & 2\\
\hline
\end{tabular}
\end{center}
\caption{\small Values of $\delta^h_k$ in the last three inner iterations of the algorithm SSN and the total number of inner iterations, for several time steps in $[0,T_f]$.}
\end{table}
In Figure \ref{fig:kine}, we show the evolution in time of $\|\mathbf{u}(t)\|_{H^{1,h}}$ and $\|\theta(t)\|_{W^{q,h}}$, in the time interval $[0, T_f]$. These pictures depict us a good insight into the variation of the kinetic energy of the flow and the heat transfer of the temperature field. It is possible to observe that both the norm of the velocity field $\|\mathbf{u}(t)\|_{H^{1,h}}$ and the energy field $\|\theta(t)\|_{W^{q,h}}$ tend toward constant limits as $t\rightarrow\infty$. This behaviour suggests that the kinetic energy of the flow tends to be constant asymptotically. A similar behavior is observed for the asymptotic behavior of the temperature field.

Regarding the behavior of Algorithm SSN, the average number of iterations is 3.27, which implies that the inner SSN algorithm needs to solve, in average, approximately 3 $2n \times 2n$-systems of linear equations per time iteration. Even if we do not apply the decomposition of the  system matrix, we need only to solve four Stokes-type linear systems per iteration, which is not computationally costly in comparison with first-order methods. This low computational cost is also a consequence of the superlinear convergence rate. In order to show this behavior, in Table \ref{tab:1}, we show the values of $\delta^h_k$ in the last three inner iterations of the algorithm SSN, for several time steps, as well as the total number of inner iterations in each one of these time steps. In this table it is possible to appreciate the fast decay of the residuum in the last iterations, and consequently, we can illustrate the local superlinear convergence rate of the inner algorithm in each time step.

\subsection{Experiment 2}
In this experiment we calculate the non-isothermal Bingham flow with the following structural parameters: $\mu_0=1.5$, $\delta_\mu=-0.5$, $g_0=18$ and $\delta_g=-8$. In this case, we focus on the flow in the regime where both the viscosity and the yield stress decrease with temperature. We consider that $\alpha=0$, so there is not an external heat sink proportional to the temperature field. We have that $\beta=15$ for the Robin boundary condition imposed in the upper edge of the square geometry. We consider that $\kappa=10$, $C_p=1.5$ and $\gamma=10^3$. We consider a mesh given by $h=0.0023$ and $\delta=0.1*h^{4/5}\approx 0.0018$. The initial condition for the flow and theta are given by $\mathbf{u}_0=\mathbf{0}$ and $\theta_0=0.0125$, respectively. Finally, we consider that $T_f=0.12$.

\begin{figure}
\begin{center}
\includegraphics[width=40mm, height=40mm]{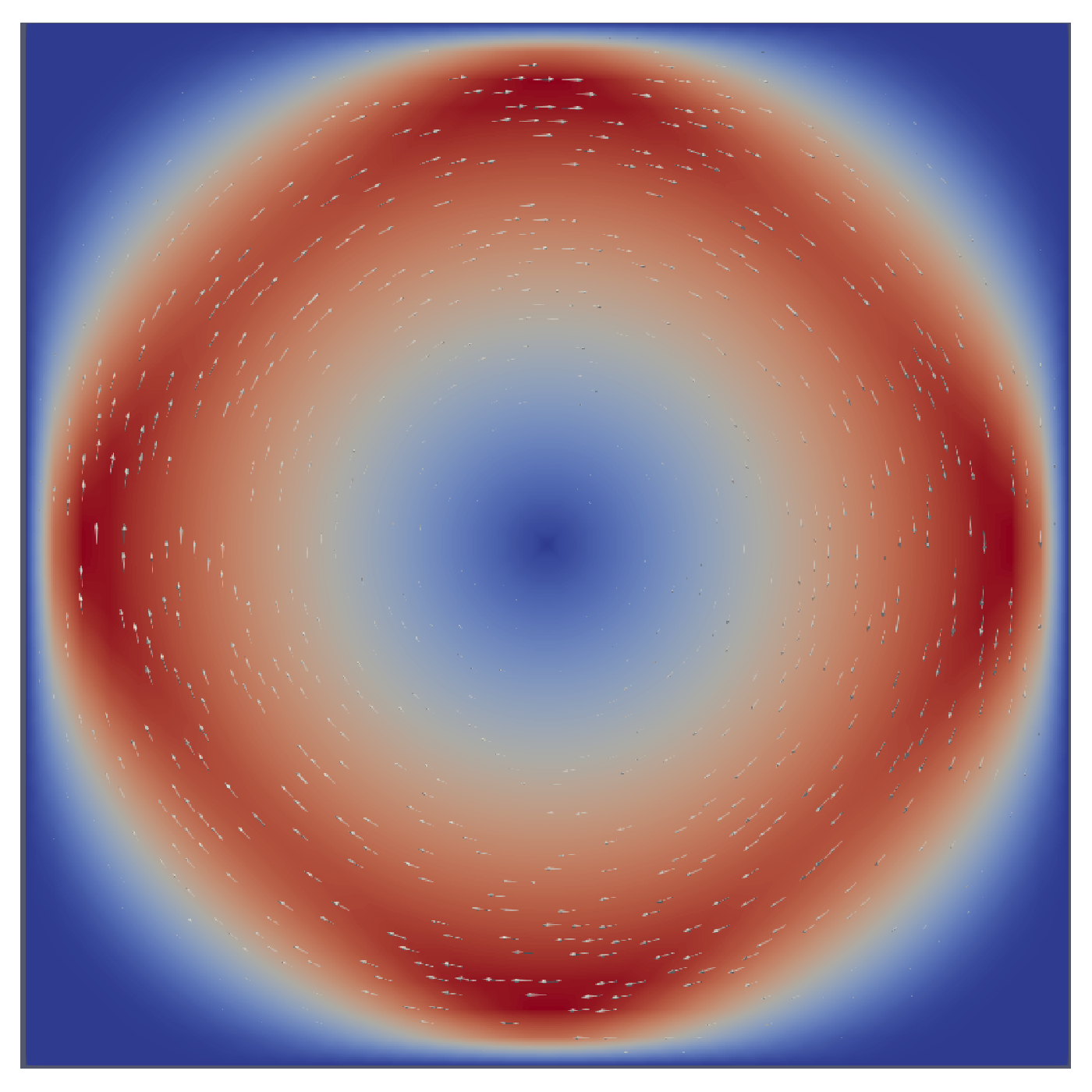}\includegraphics[width=40mm, height=40mm]{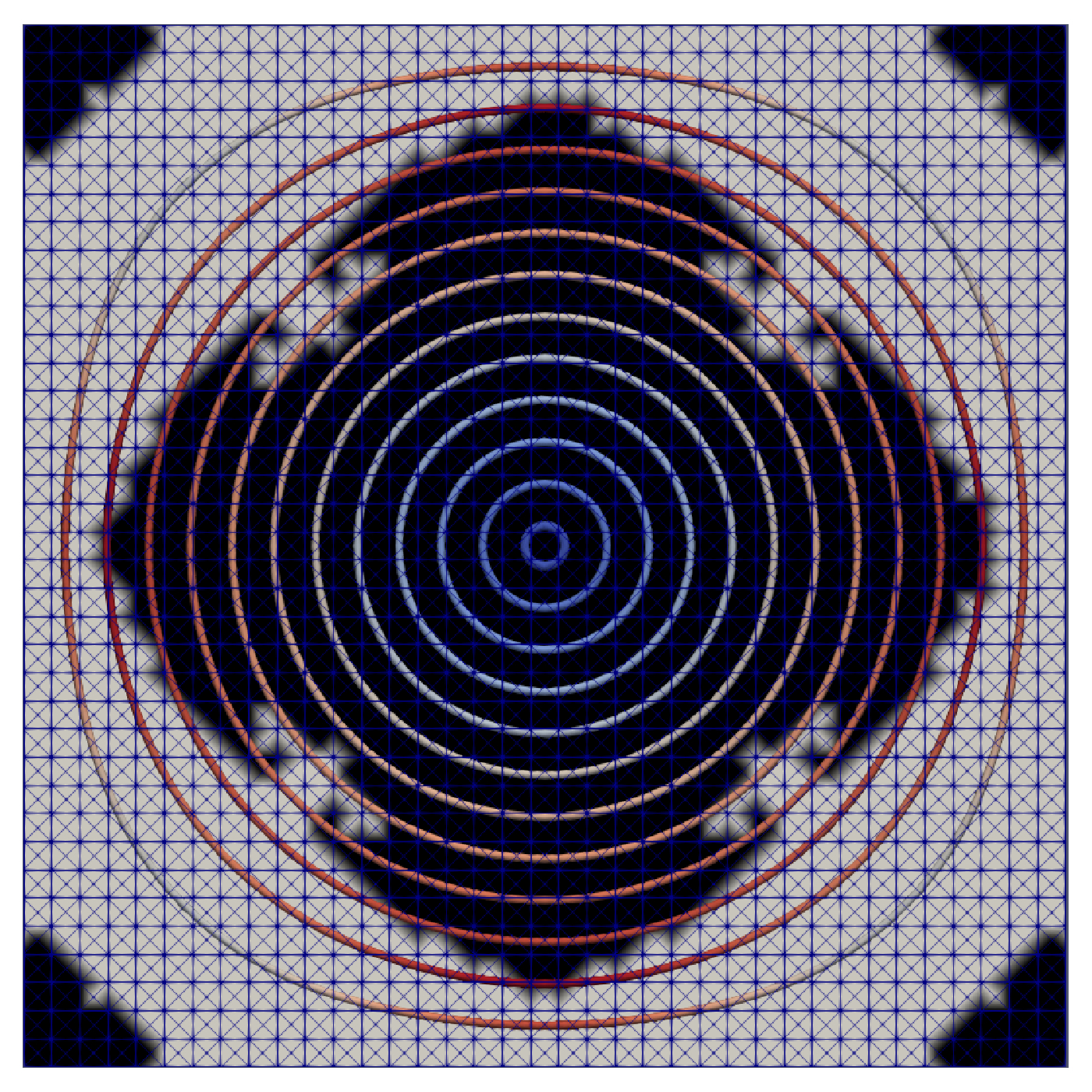}\includegraphics[width=50mm, height=40mm]{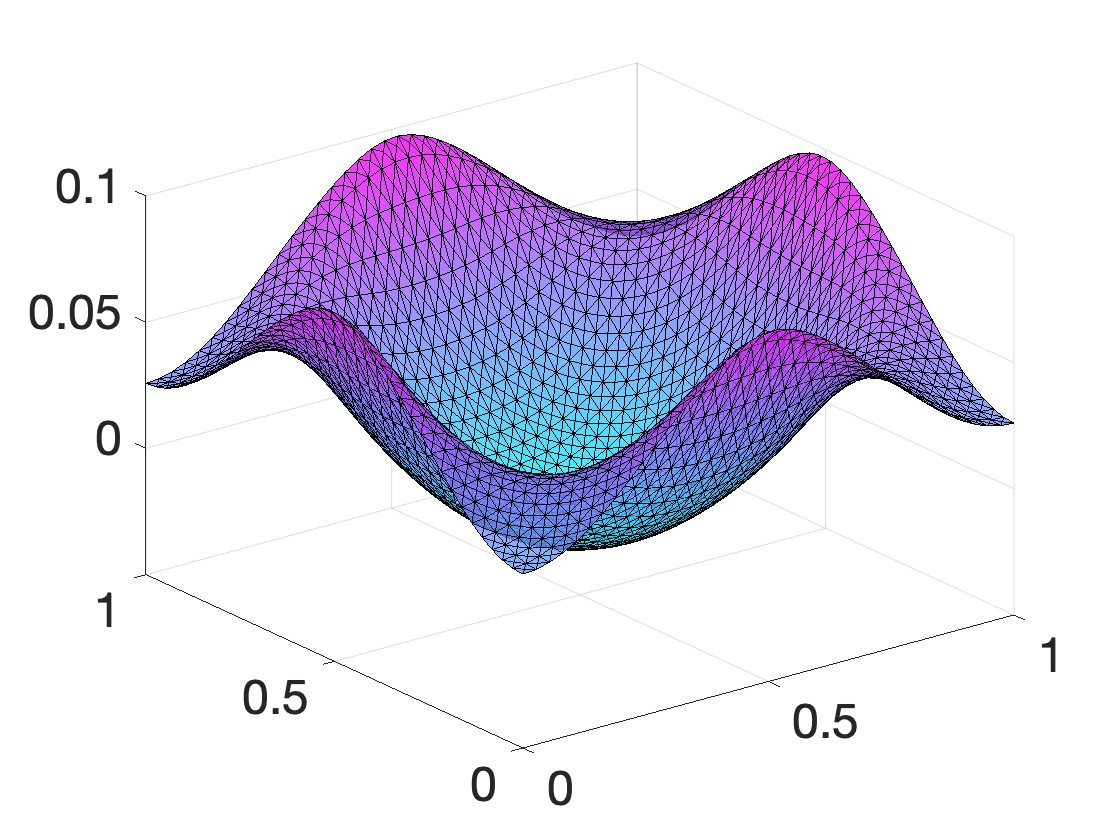}\\
\includegraphics[width=40mm, height=40mm]{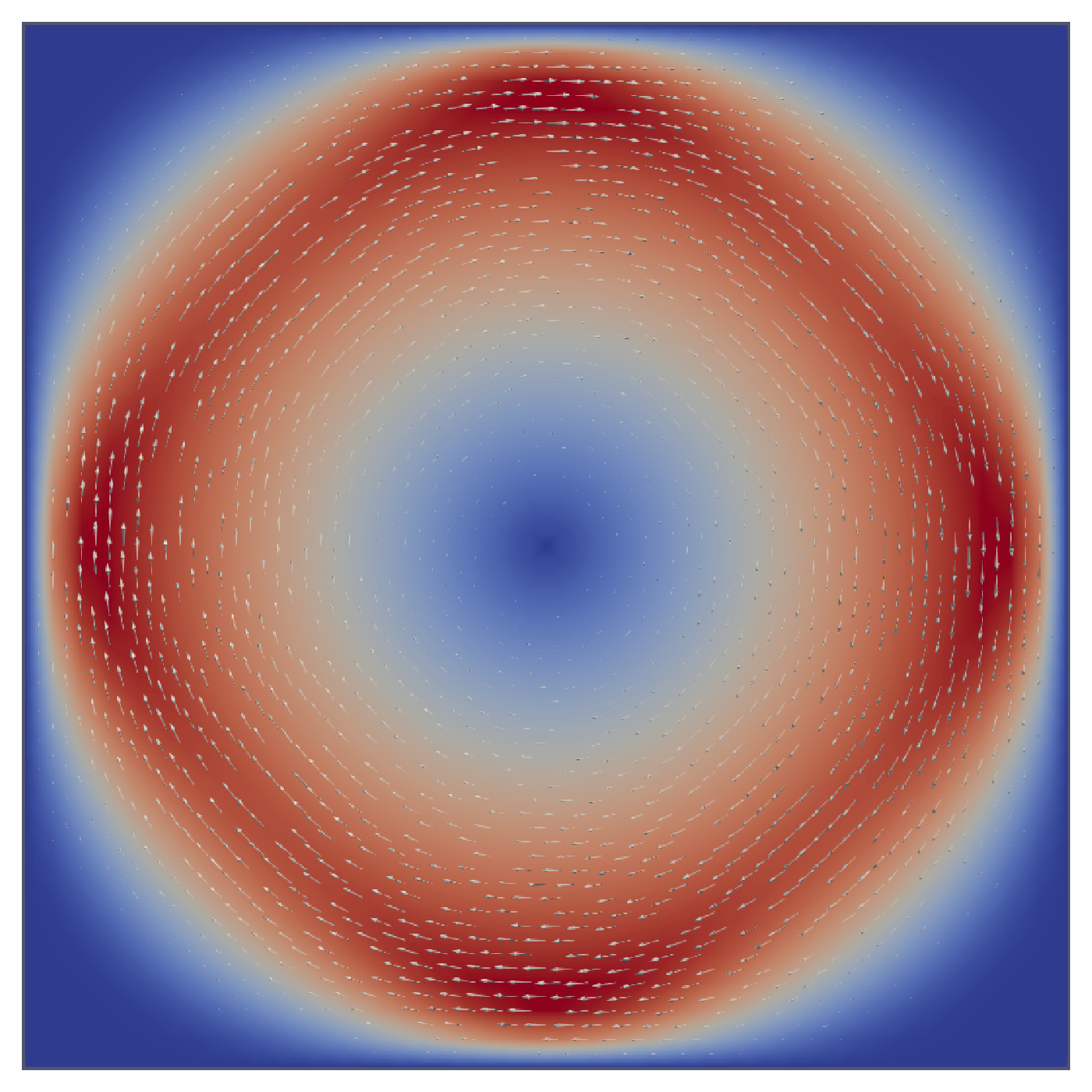}\includegraphics[width=40mm, height=40mm]{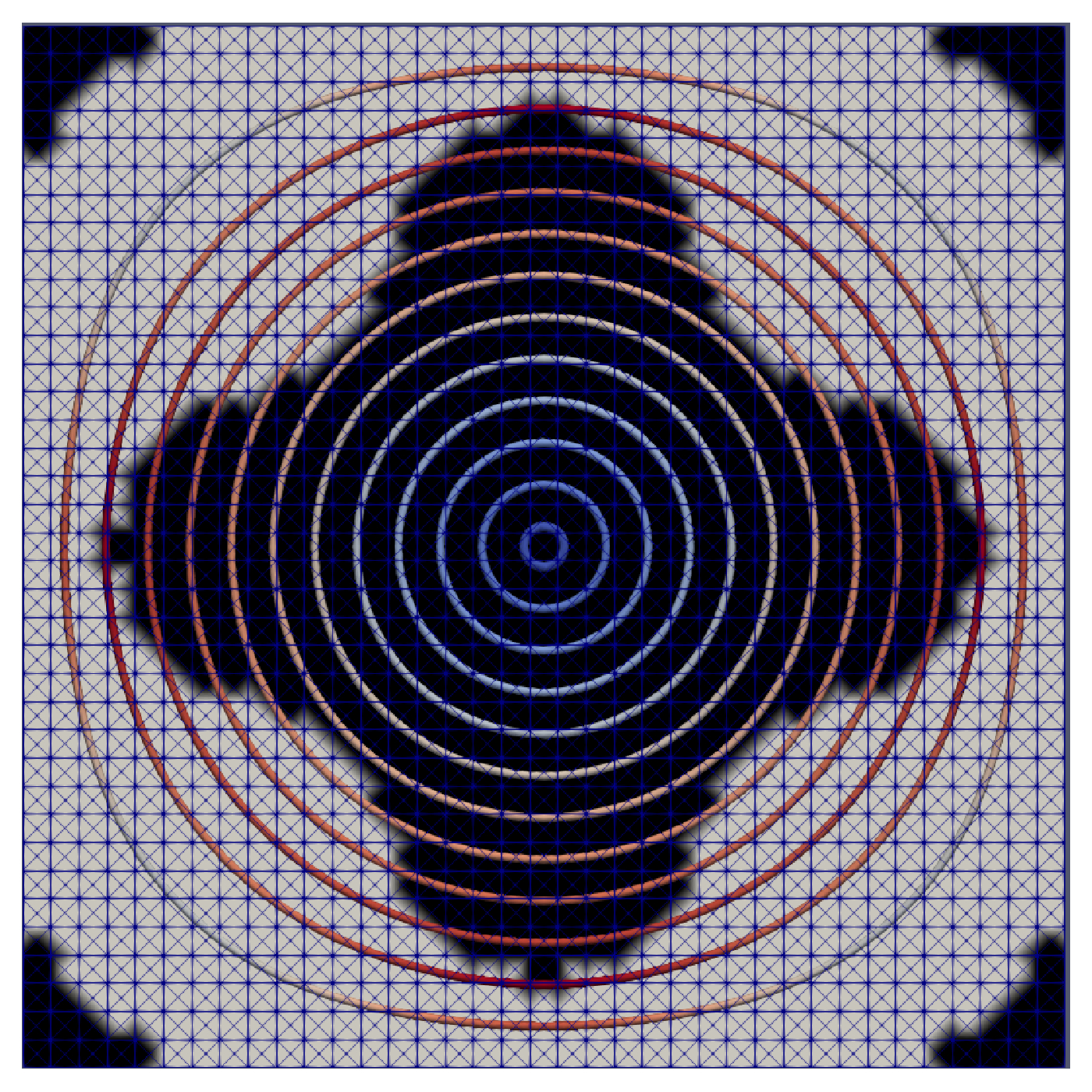}\includegraphics[width=50mm, height=40mm]{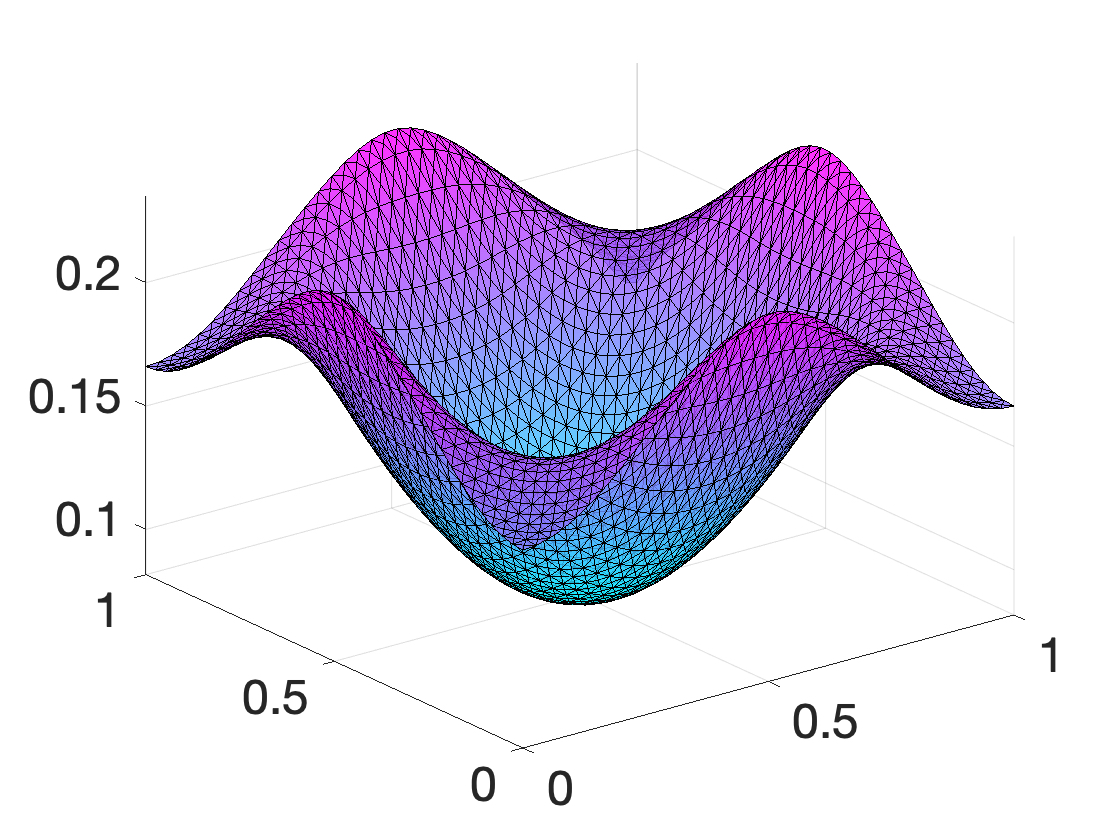}\\
\includegraphics[width=40mm, height=40mm]{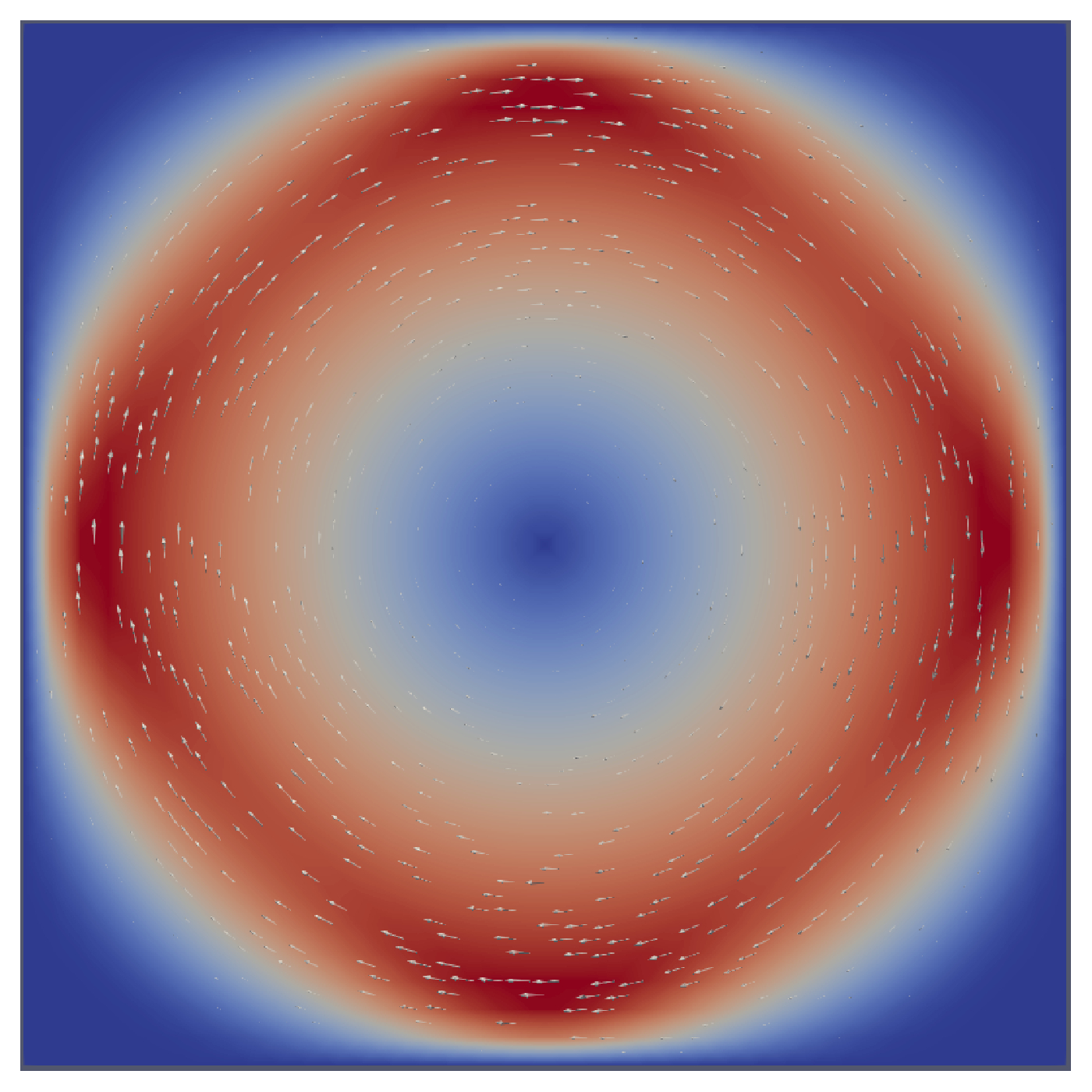}\includegraphics[width=40mm, height=40mm]{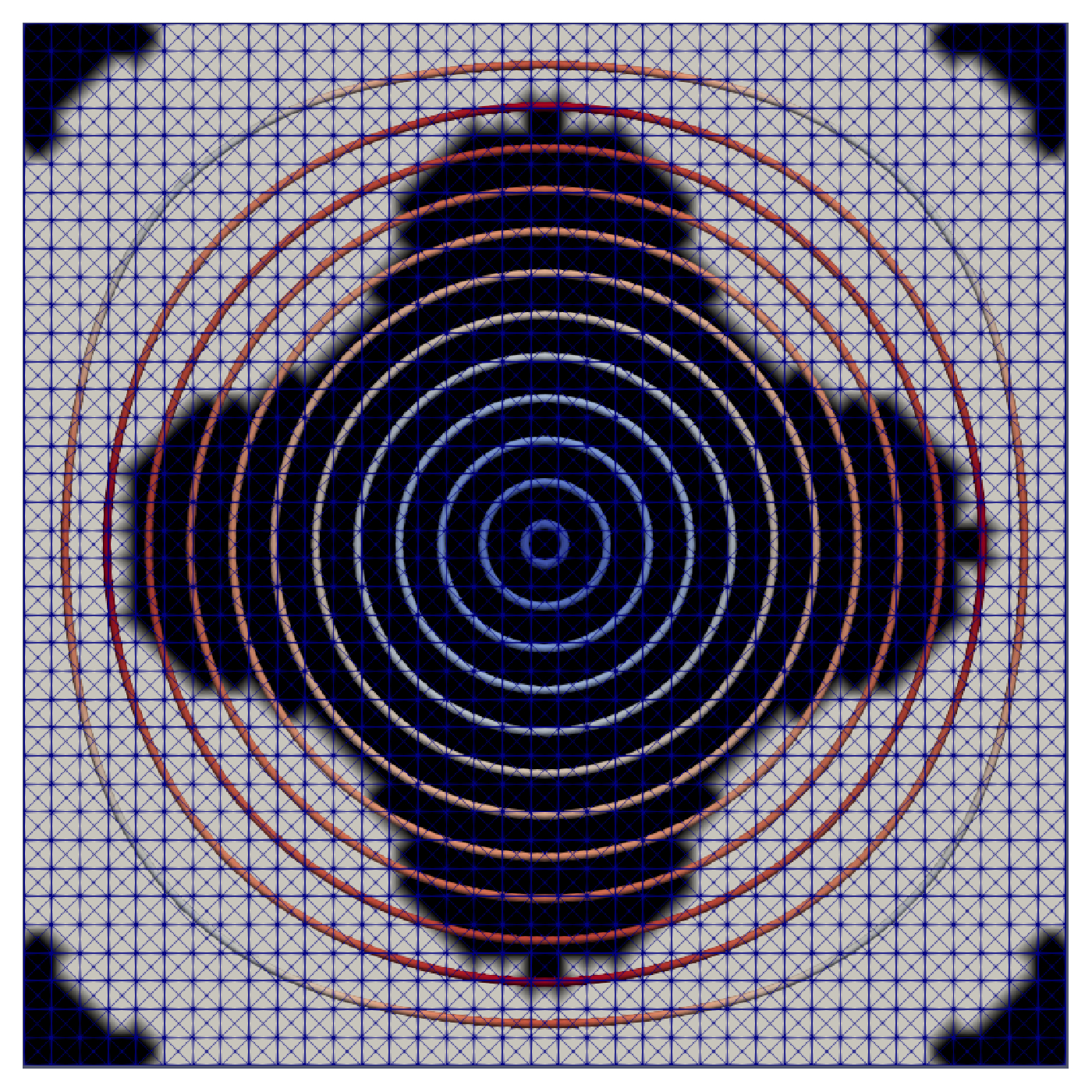}\includegraphics[width=50mm, height=40mm]{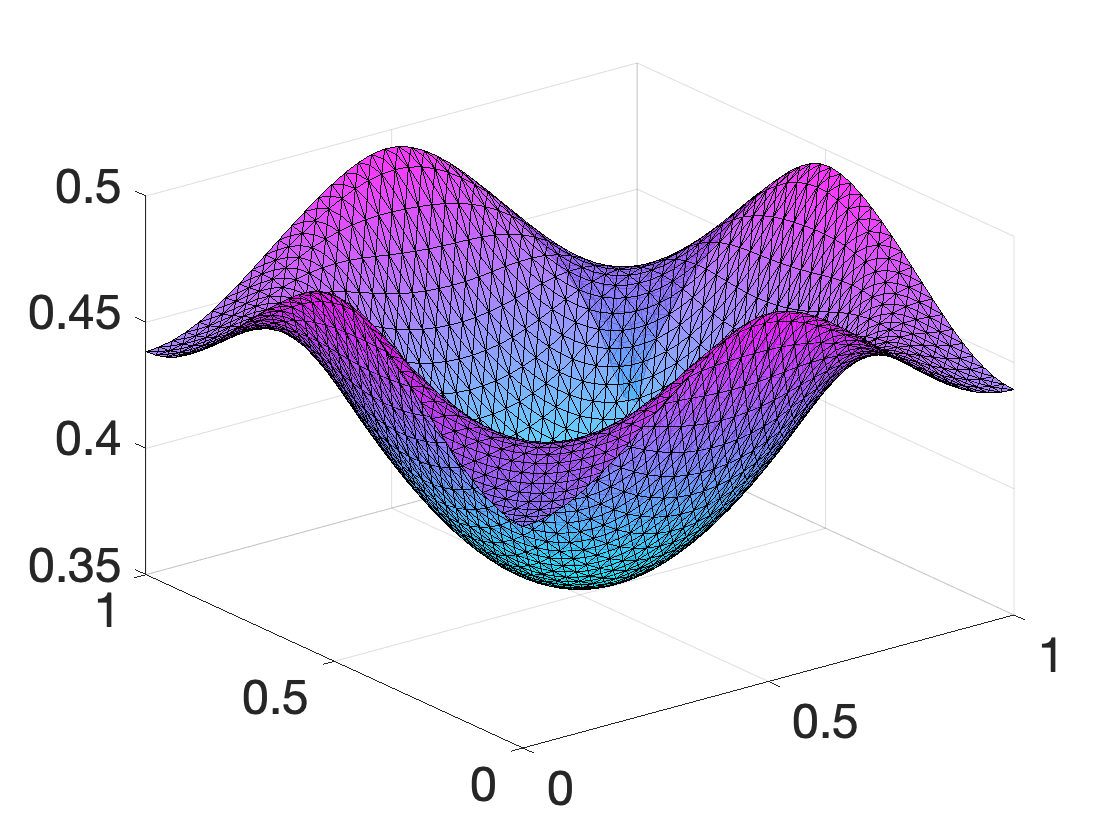}\\
\includegraphics[width=40mm, height=40mm]{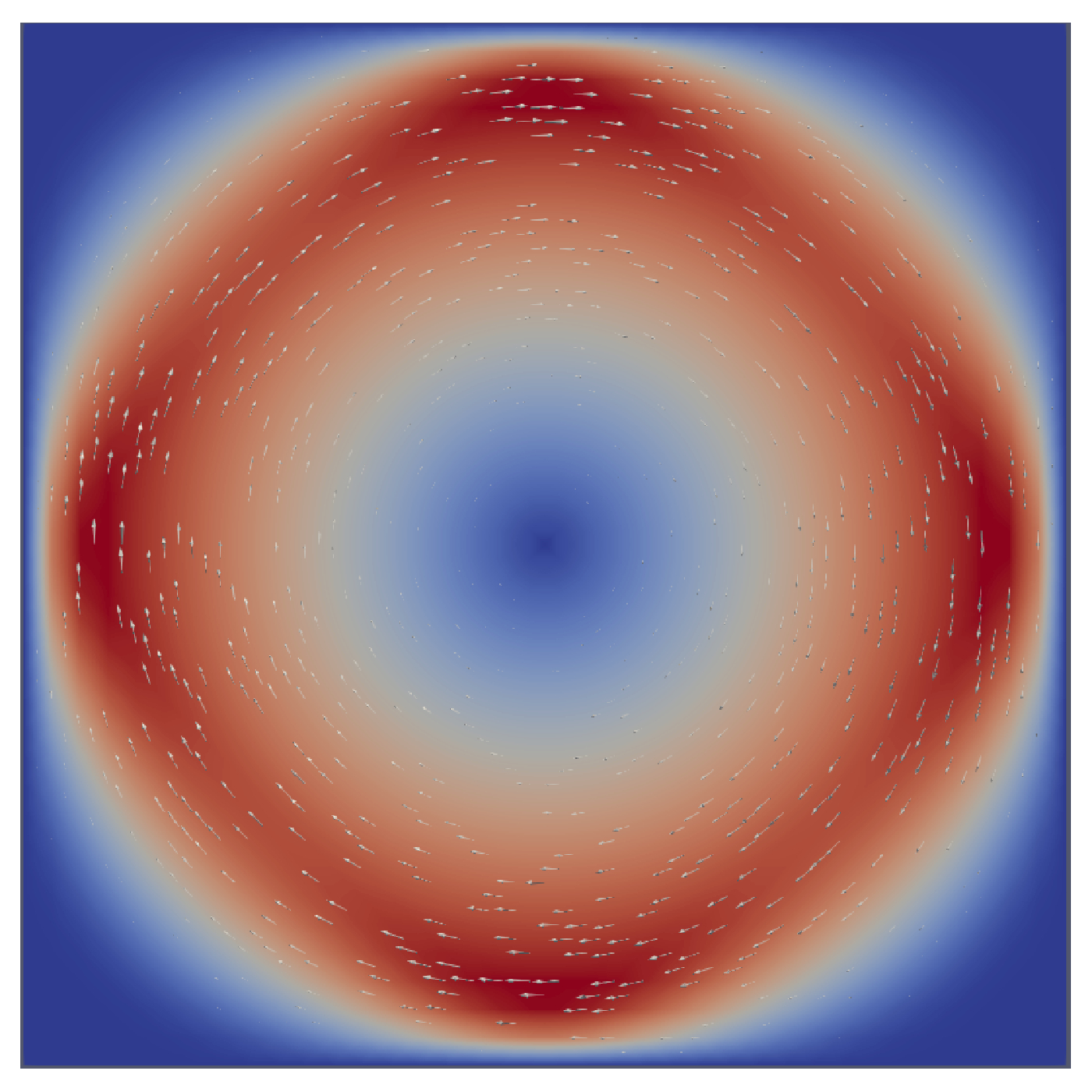}\includegraphics[width=40mm, height=40mm]{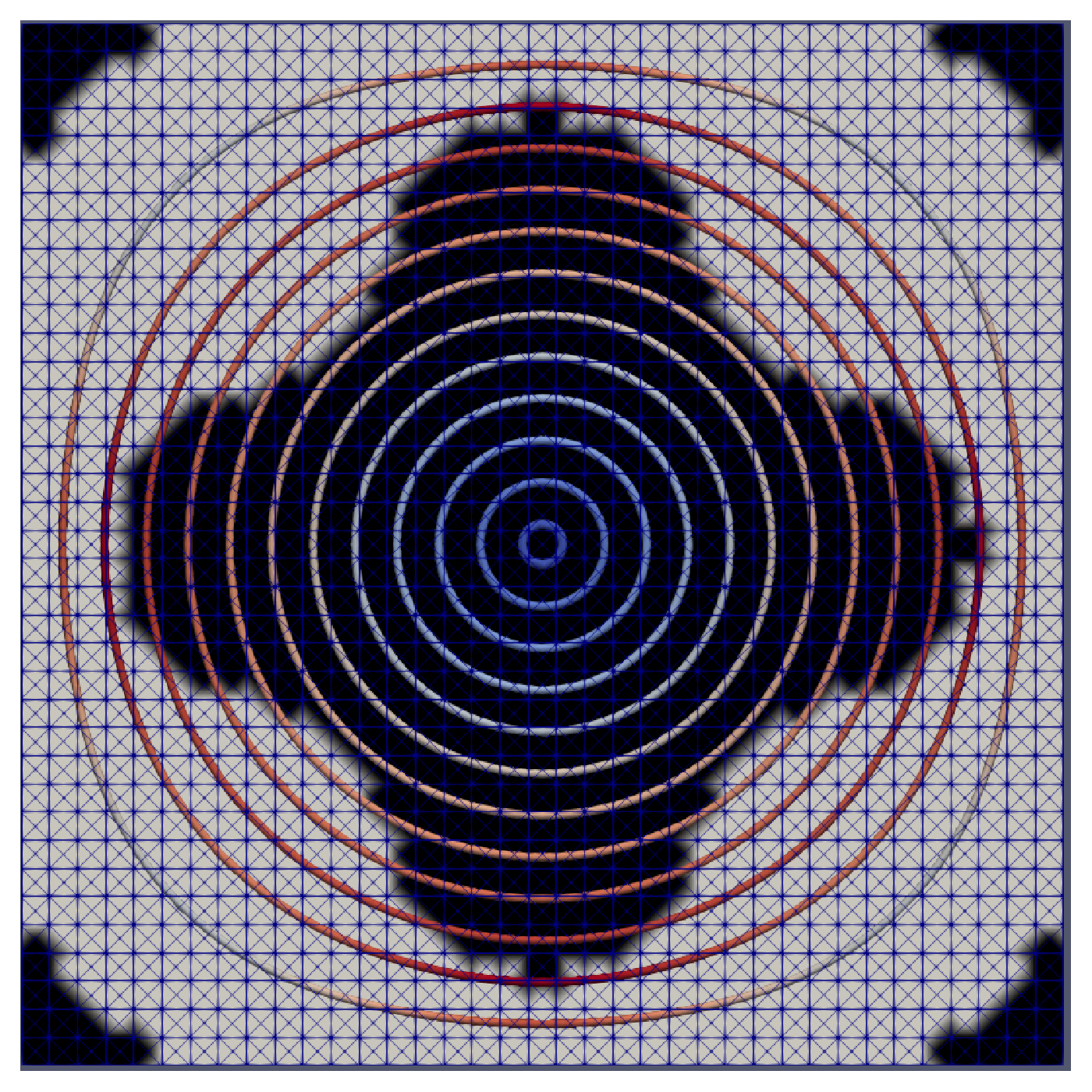}\includegraphics[width=50mm, height=40mm]{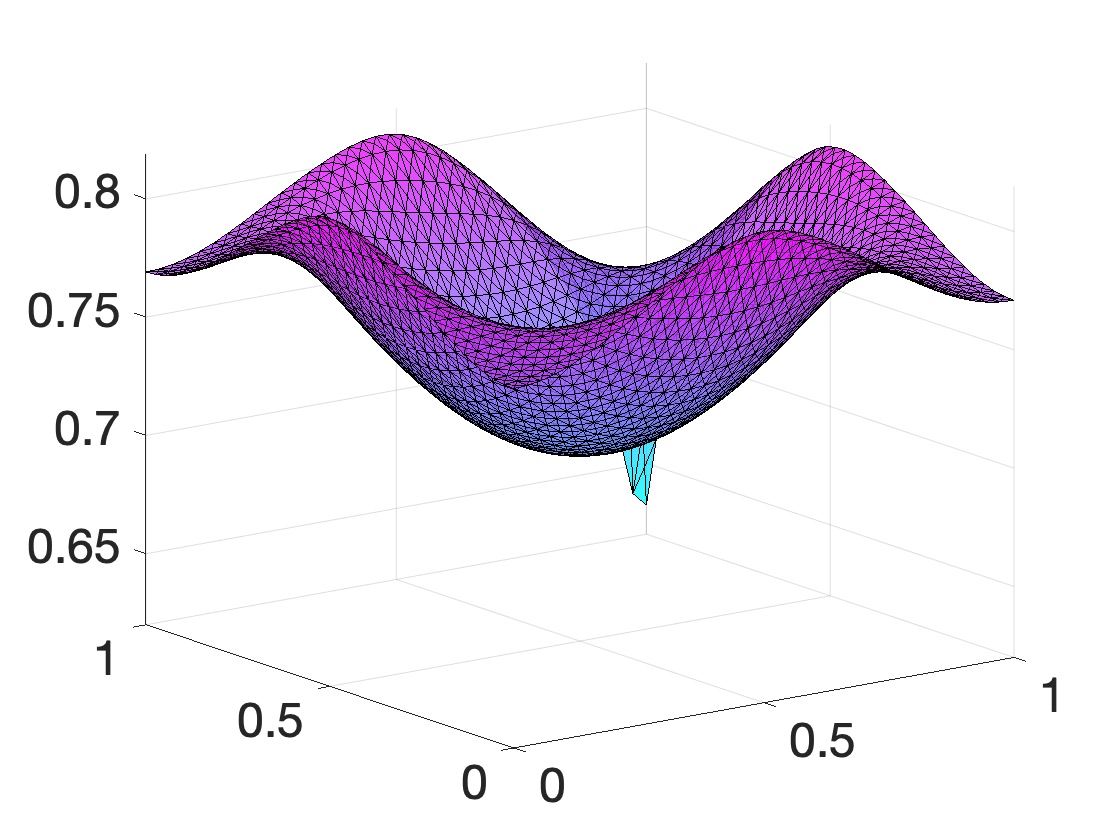}\\
\end{center}
\caption{\small Experiment 1. Left column: velocity flow. Central column: active and inactive sets with streamlines. Right column: temperature field. For $t=0.015$ (up), $t=0.030$, $t=0.060$ and $T_f=0.12$ (bottom).}\label{fig:flowexp2}
\end{figure}

\begin{figure}
\begin{center}
\includegraphics[width=60mm, height=40mm]{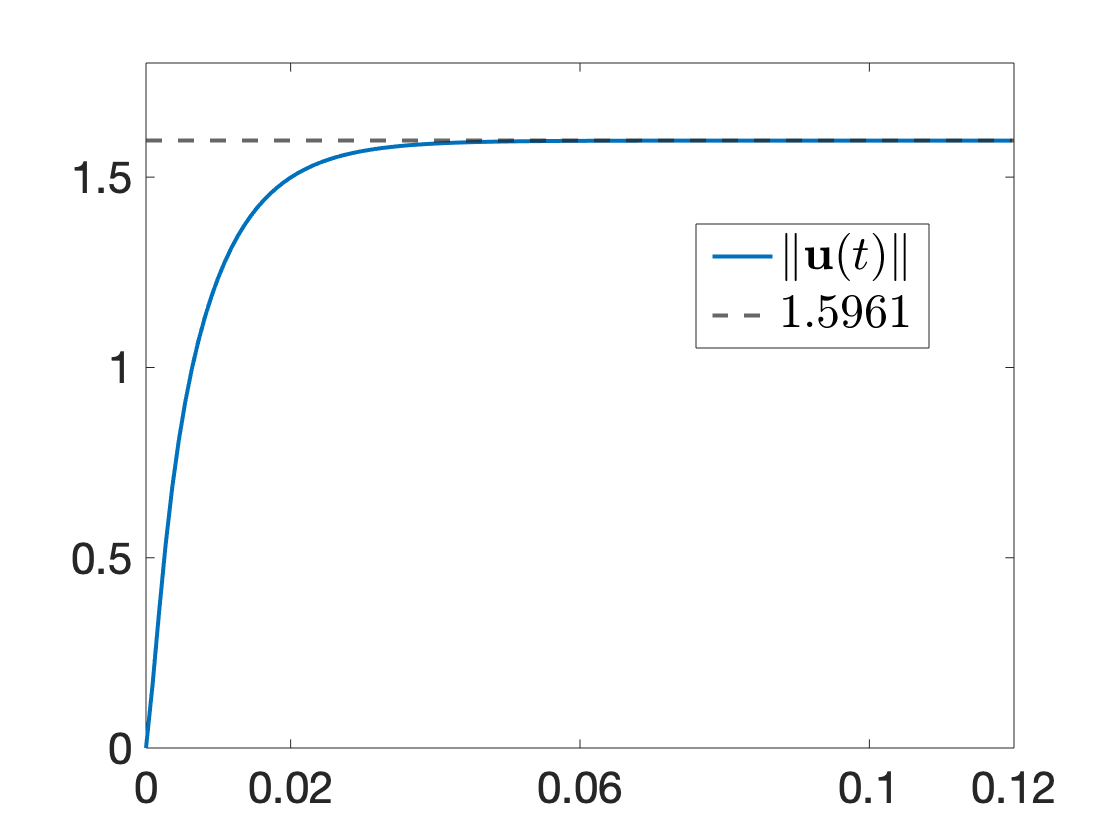}\hspace{1cm}\includegraphics[width=60mm, height=40mm]{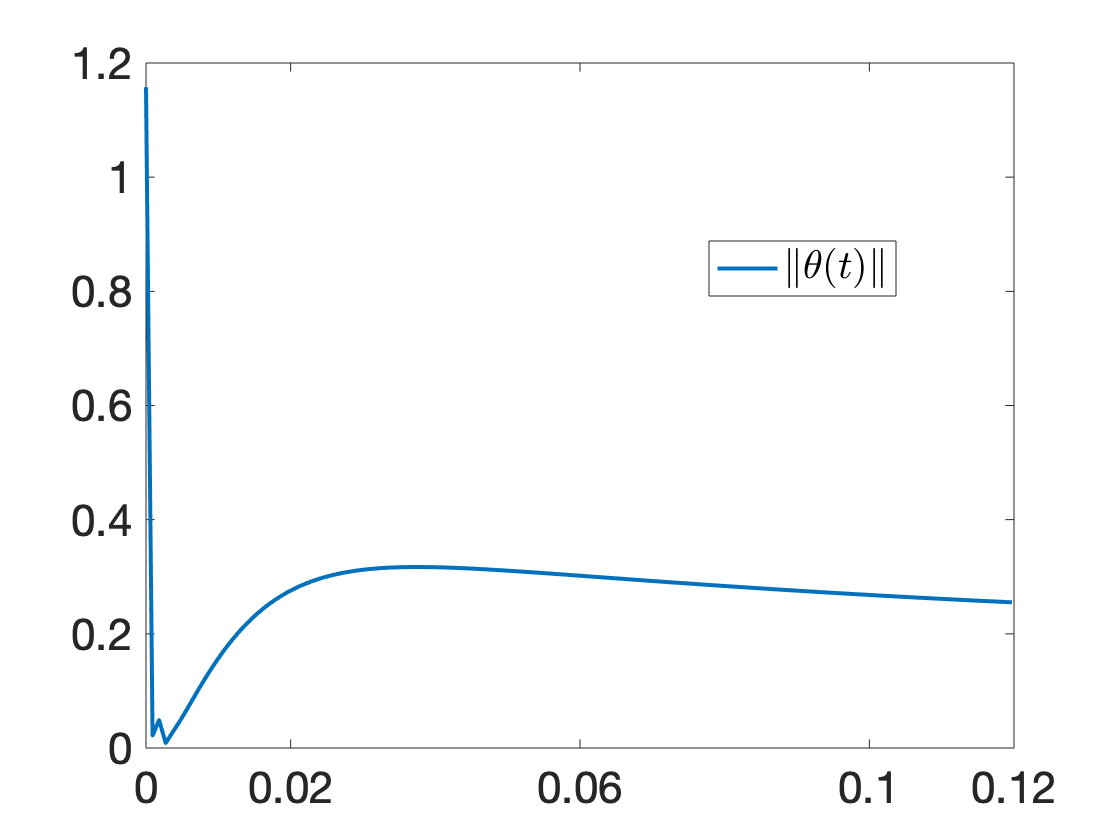}
\end{center}
\caption{\small Experiment 2. Evolution of $\|\mathbf{u}(t)\|_{H^{1,h}}$ (left) and $\|\theta(t)\|_{W^{q,h}}$ (right).}\label{fig:kine1}
\end{figure}

In Figure \ref{fig:flowexp2}, we show the calculated velocity field, the active and inactive sets, representing the rigid and plastic
regions of the material, as well as the streamlines of the flow, and the calculated temperature fields, for several instants. In this case, the flow also behaves as expected: the rigid regions change from the regions given by a constant value of $g$, close to $g_0$ for the initial values, and then evolves. However, the change in the unyielded regions is not as clear as in the former experiment. This can be explained considering that there is not a heat source nor heat sink acting, so the motion is only driven by the body forces and the convective terms. The flow tends to stabilize around $T_f=0.03$. This effect can be appreciated in Figure \ref{fig:kine1}, where the  evolution in time of $\|\mathbf{u}(t)\|_{H^{1,h}}$ and $\|\theta(t)\|_{W^{q,h}}$ is showed. It is possible to see that the velocity norm tends to a constant value as $t$ increases. This effect, actually, is very fast. In constrast, the behavior of the temperature norm, although has the same tendency, it does not clearly reach a limit. This happens even in longer time intervals. This suggests that the heat sink plays an important role in shaping the flow. 

\begin{figure}
\begin{center}
\includegraphics[width=45mm, height=30mm]{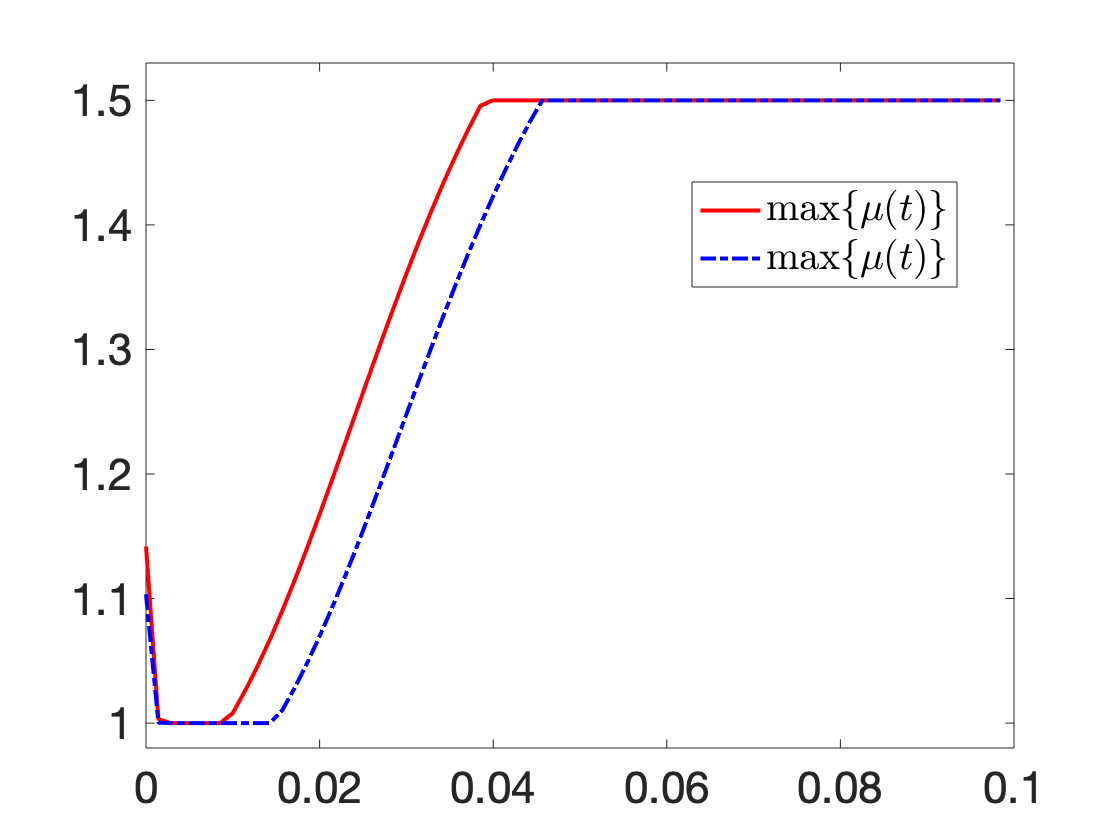}\hspace{1cm}\includegraphics[width=45mm, height=30mm]{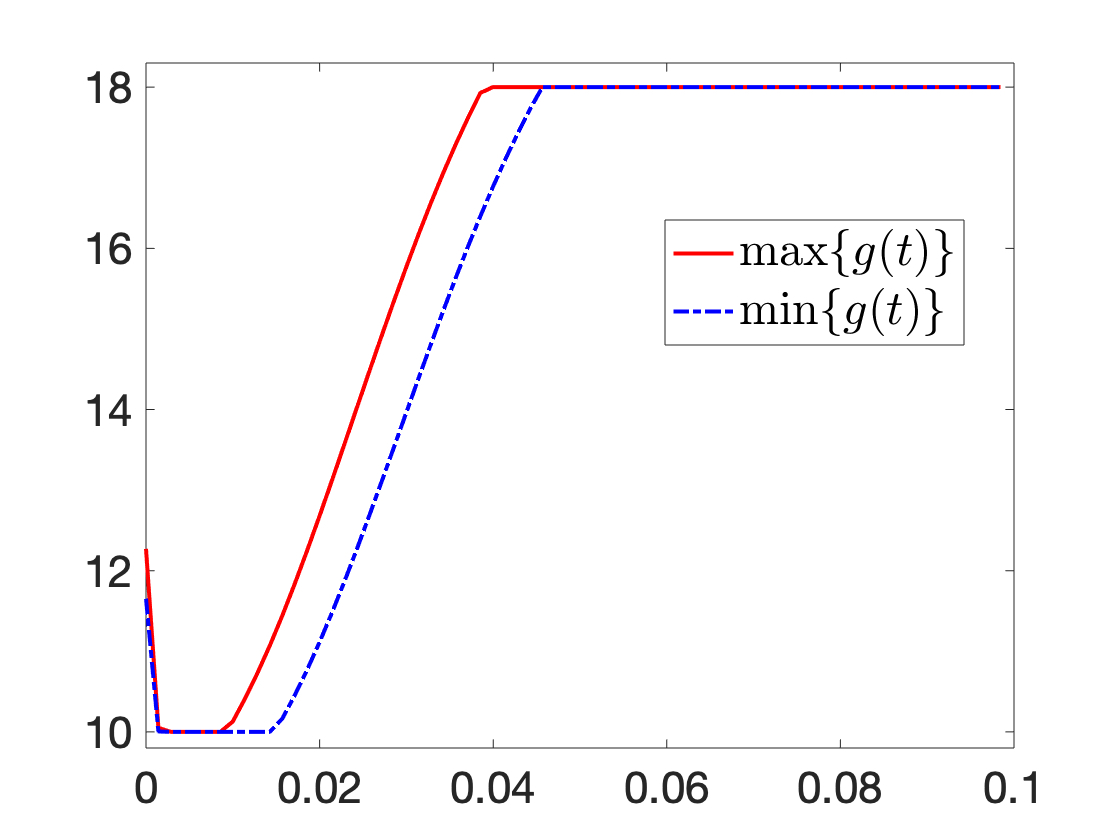}\hspace{1cm}\includegraphics[width=45mm, height=30mm]{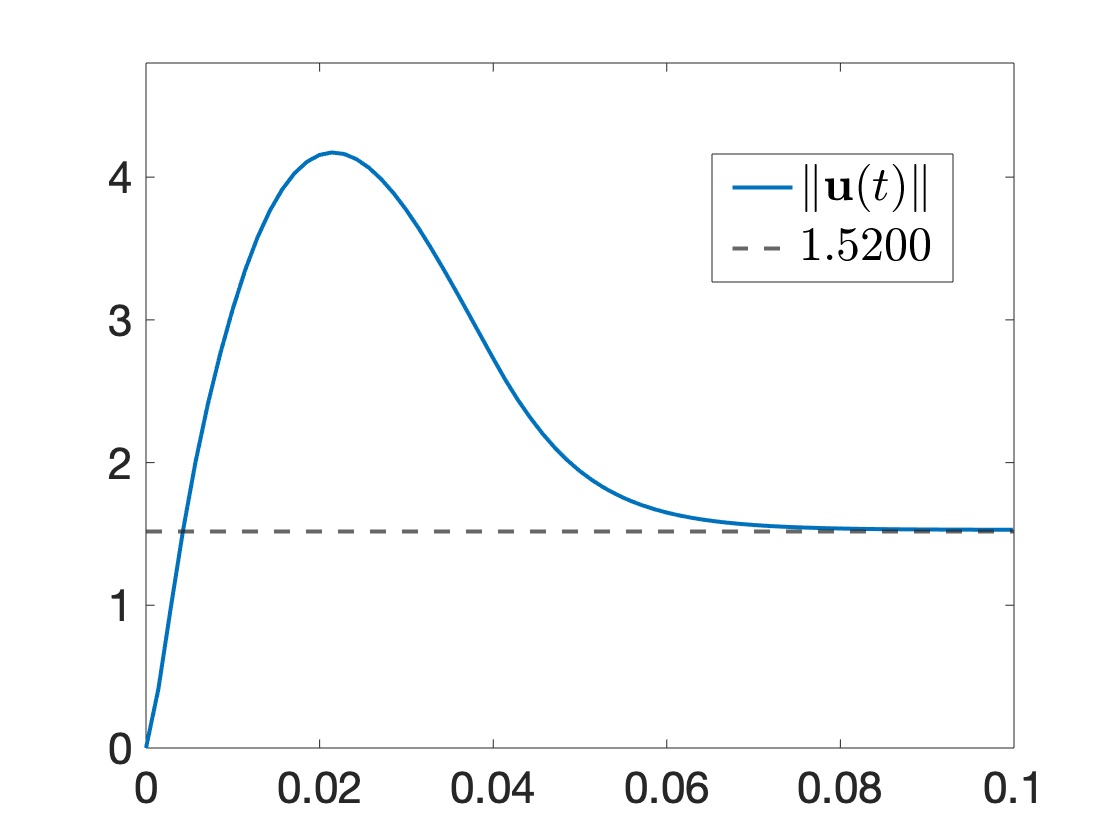}\\
\includegraphics[width=45mm, height=30mm]{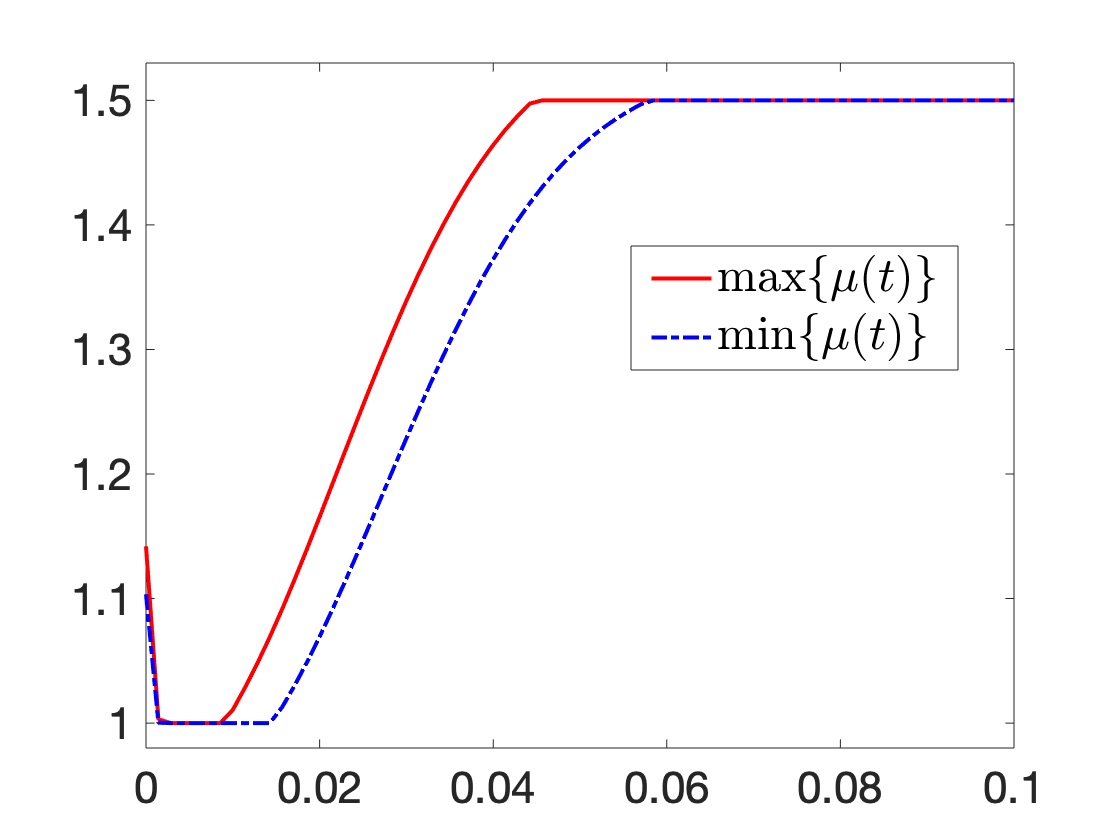}\hspace{1cm}\includegraphics[width=45mm, height=30mm]{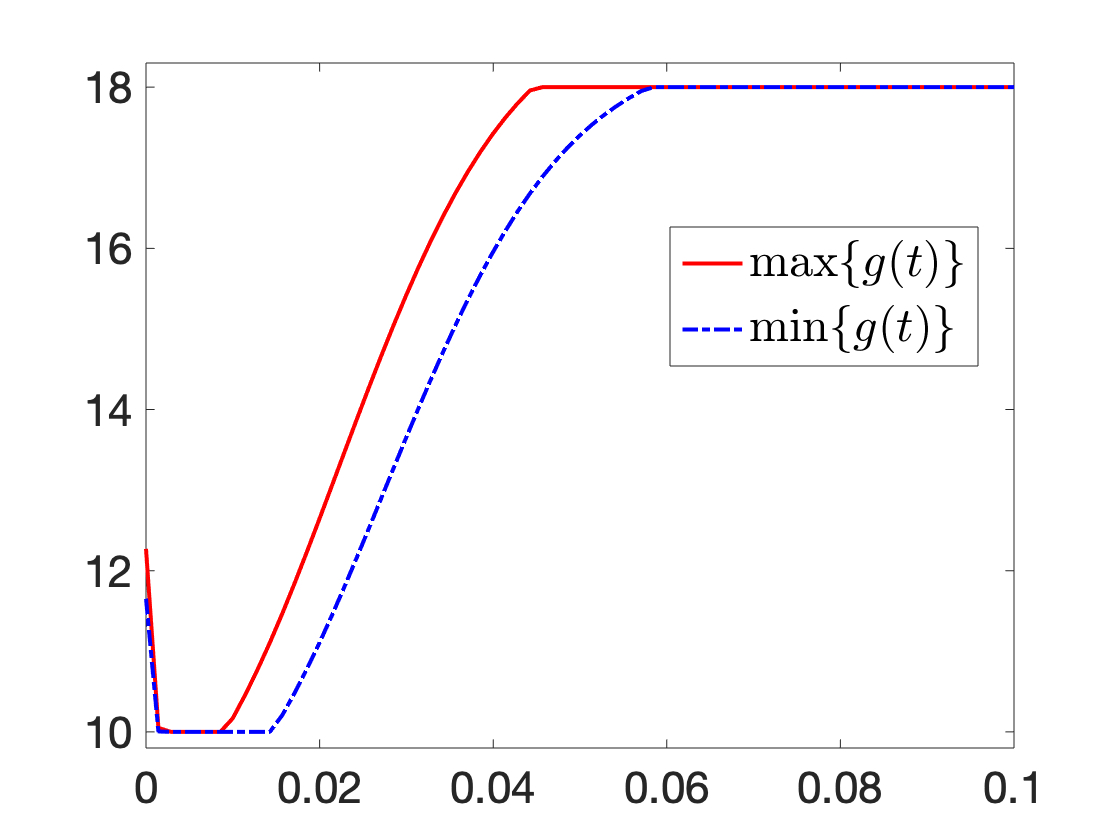}\hspace{1cm}\includegraphics[width=45mm, height=30mm]{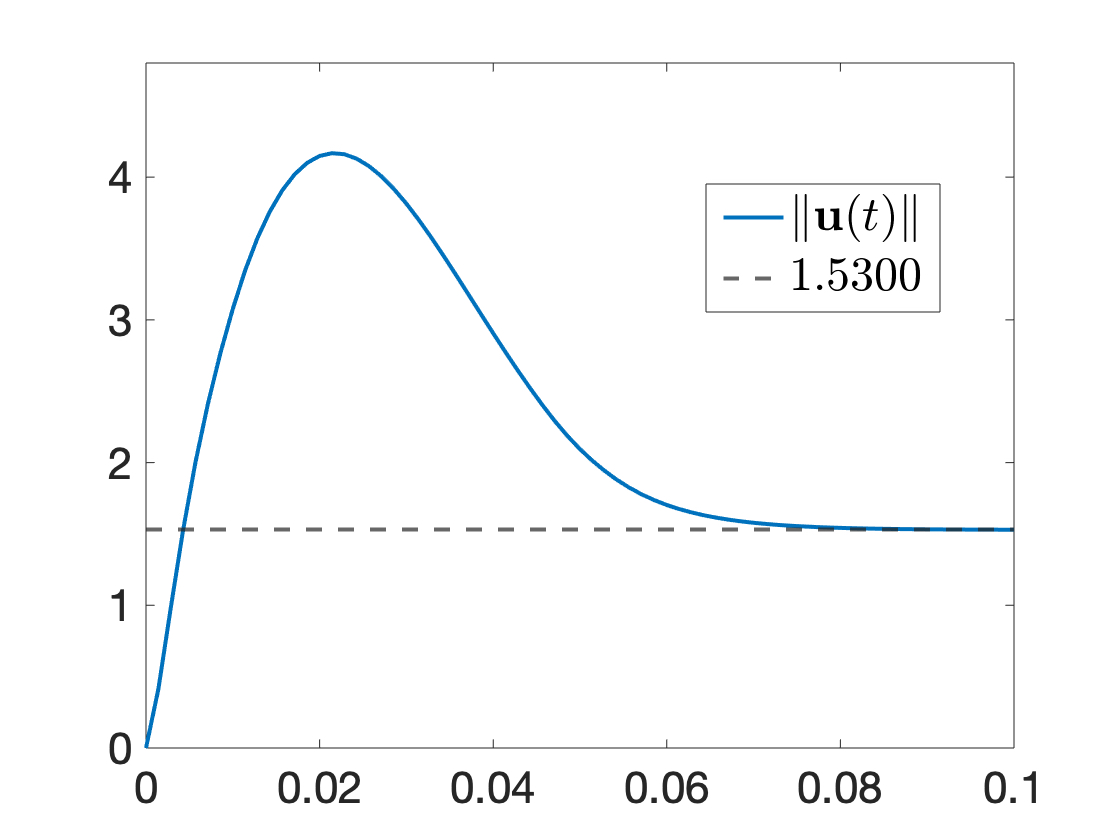}\\
\includegraphics[width=45mm, height=30mm]{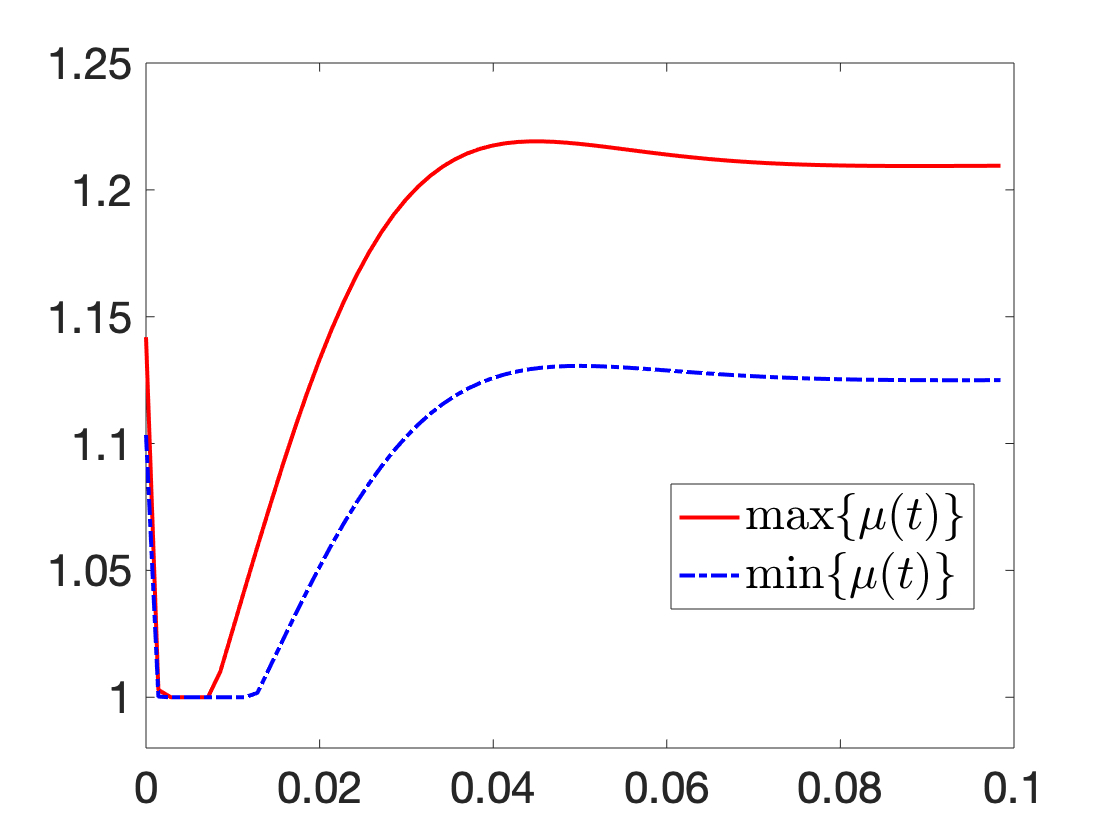}\hspace{1cm}\includegraphics[width=45mm, height=30mm]{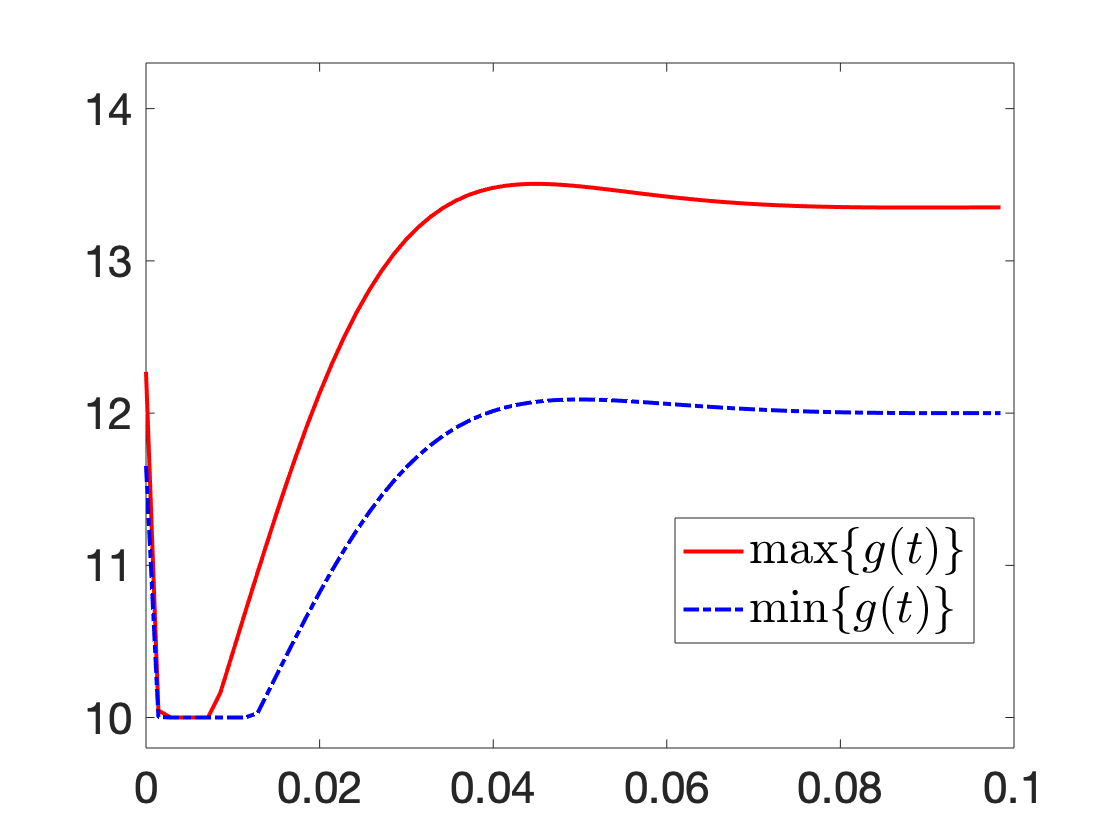}\hspace{1cm}\includegraphics[width=45mm, height=30mm]{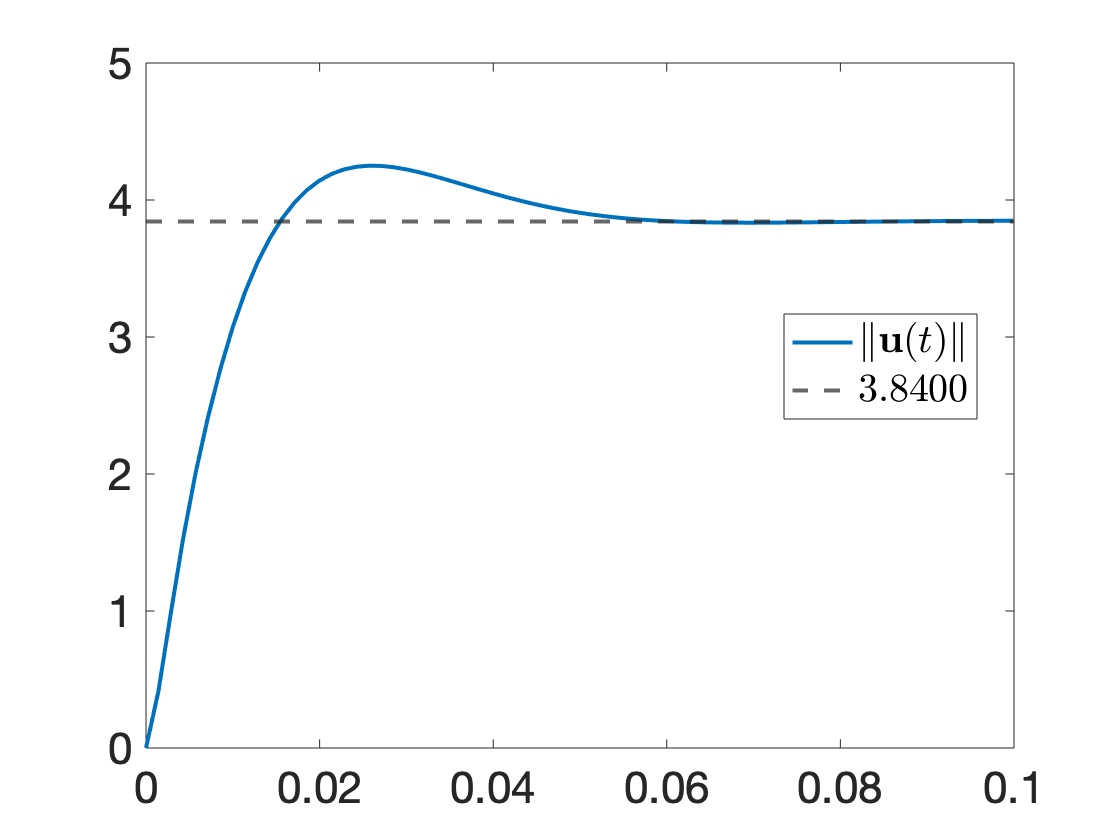}
\end{center}
\caption{\small Evolution of $\min\{\mu(t)\}-\max\{\mu(t)\}$ (left column), $\min\{g(t)\}-\max\{g(t)\}$ (center column) and $\|\mathbf{u}(t)\|_{H^{1,h}}$ (right column), in the interval $[0,T_f]$, for $\alpha=1$ (upper row), $\alpha=10$ (middle row) and $\alpha=100$ (lower row).}\label{fig:alfa1-100}
\end{figure}

In order to show the effect of $\alpha$ in the flow, in Figure \ref{fig:alfa1-100}, we show the evolution of the minimum and maximum values for the parameter functions $\mu$ and $g$, as well as the evolution of the norm of the velocity field, in the interval $[0,T_f]$, for several values of $\alpha$ and the following fixed parameters: $\beta=1$, $\mu_0=1$, $\delta_\mu=0.5$, $g_0$, $\delta_g=8$, $\kappa=10$ and $C_p=1$. It is possible to see a direct impact of the parameter $\alpha$ in the behaviour of the flow. This impact opens the door for a control strategy depending on the heat sink/source. This fact deserves a deeper discussion in future contributions. 

\section[Conclusions]{Conclusions}
Based on a Huber regularized multiplier approach, we proposed and implemented a combined BDF-SSN algorithm for the numerical solution of non-isothermal, time-dependent Bingham flow with temperature dependent parameters. The finite elements based on the so called (cross-grid $\mathbb{P}_1$)-$\mathbb{Q}_0$ elements provided a LBB-stable approximation, which in this case, involved weighted stiffness and mass matrices. This weighted matrices were detailed implemented and provided a nice structure for the resulting systems of equations to be solved. The use of a BDF2 approach for the time stepping of each regularized system leads to a semi-implicit method with a nonsmooth convective free system of equations in each time step. As the resulting system involves slantly differentiable functions, we proposed a SSN algorithm for its numerical solution. At each time step, the algorithm converges locally with a superlinear convergence rate. Moreover, the system to be solved in each SSN step is uncoupled, resulting in an efficient combined technique. The two numerical experiments presented showed the efficiency and accuracy of our numerical approach. Mainly, the computed solution exhibits the main expected physical properties, specially, well-shaped and accurate enough yielded and unyielded regions. There are several issues to be considered in future contributions. For instance, the extension of this methodology to other viscoplastic flows like the Casson or the Herschel-Bulkley models. Also, the design of control strategies, based on the coefficient $\alpha$ in the heat source/sink looks like a promising research field. Finally, the use of different expressions for the viscosity and the yield stress functions is an issue that needs to be tackled.

\section*{Acknowledgement}
The author would like to thank the anonymous referees for many helpful comments which lead to a significant improvement of the article.

\end{document}